\documentclass[12pt]{amsart}
\usepackage{amsfonts,amssymb,amscd,amsmath}
\usepackage{latexsym}
\usepackage{mathrsfs}
\usepackage{tocvsec2}

\def\~{{\rm --}}

\begin{document}

\renewcommand{\tilde}{\widetilde}
\renewcommand{\hat}{\widehat}

\newcommand{\BR}{{\mathbb R}}
\newcommand{\BQ}{{\mathbb Q}}
\newcommand{\BC}{{\mathbb C}}
\newcommand{\BP}{{\mathbb P}}
\newcommand{\BZ}{{\mathbb Z}}
\newcommand{\BN}{{\mathbb N}}
\newcommand{\BS}{{\mathbb S}}

\newcommand{\cH}{{\mathcal H}}
\newcommand{\cA}{{\mathcal A}}
\newcommand{\cB}{{\mathcal B}}
\newcommand{\ccF}{{\mathfrak F}}
\newcommand{\cD}{{\mathcal D}}
\newcommand{\cL}{{\mathcal L}}
\newcommand{\cF}{{\mathcal F}}
\newcommand{\cP}{{\mathcal P}}
\newcommand{\cX}{{\mathcal X}}
\newcommand{\cY}{{\mathcal Y}}
\newcommand{\cS}{{\mathcal S}}
\newcommand{\cSol}{\hbox{$\mathcal Sol$}}
\newcommand{\cT}{\hbox{$\mathcal T$}}

\newcommand{\Z}{{\mathbb Z}}
\newcommand{\Q}{{\mathbb Q}}
\newcommand{\N}{{\mathbb N}}
\newcommand{\C}{{\mathbb C}}
\newcommand{\R}{{\mathbb R}}

\newcommand{\CH}{{\mathcal H}}
\newcommand{\CA}{{\mathcal A}}

\def\HH{\mbox{${\mathcal H}$\kern-5.2pt${\mathcal H}$}}

\newcommand{\binomial}[2]{\genfrac{(}{)}{0pt}{}{ #1 }{ #2 }}
\newcommand{\qbinomial}[2]{\genfrac{[}{]}{0pt}{}{ #1 }{ #2 }_q }
\newcommand{\qbinom}[3]{\genfrac{[}{]}{0pt}{}{ #1 }{ #2 }_{ #3 } }


\def\der{\partial}
\def\tensor{\otimes}
\def\gam{\gamma} \def\Gam{\Gamma}
\def\del{\delta} \def\Del{\Delta}
\def\kap{\kappa}
\def\lam{\lambda} \def\Lam{\Lambda}
\def\Comp{{\mathbb C}}
\def\sM{{\mathcal M}}

\newtheorem{theorem}{Theorem}[section]
\newtheorem{maintheorem}[theorem]{Main Theorem}
\newtheorem{proposition}[theorem]{Proposition}
\newtheorem{definition}[theorem]{Definition}
\newtheorem{lemma}[theorem]{Lemma}
\newtheorem{corollary}[theorem]{Corollary}
\newtheorem{notation}[theorem]{Notation}
\newtheorem{remark}[theorem]{Remark}
\newtheorem{example}[theorem]{Example}

\newtheorem{theorem }{Theorem}[section]
\newtheorem{maintheorem }[theorem]{Main Theorem}
\newtheorem{proposition }[theorem]{Proposition}
\newtheorem{definition }[theorem]{Definition}
\newtheorem{lemma }[theorem]{Lemma}
\newtheorem{corollary }[theorem]{Corollary}
\newtheorem{notation }[theorem]{Notation}
\newtheorem{remark }[theorem]{Remark}
\newtheorem{example }[theorem]{Example}

\newtheorem{ maintheorem }[theorem]{Main Theorem}
\newtheorem{ theorem}{Theorem}[section]
\newtheorem{ proposition}[theorem]{Proposition}
\newtheorem{ definition}[theorem]{Definition}
\newtheorem{ lemma}[theorem]{Lemma}
\newtheorem{ corollary}[theorem]{Corollary}
\newtheorem{ notation}[theorem]{Notation}
\newtheorem{ remark}[theorem]{Remark}
\newtheorem{ example}[theorem]{Example}

\newtheorem{thm}{Theorem}[section]
\newtheorem{prop}[thm]{Proposition}
\newtheorem{lem}[thm]{Lemma}
\newtheorem{cor}[thm]{Corollary}
\newtheorem{conj}[thm]{Conjecture}
\newtheorem{con}[thm]{Conjecture}
\newtheorem{dfn}[thm]{Definition}
\newtheorem{df}[thm]{Definition}
 \newcommand{\rem}{{\bf Comment.\ }}
 \newcommand{\rmk}{{\bf Comment.\ }}
 \newcommand{\exmp}{{\bf Example.\ }}
 \newcommand{\ex}{{\bf Example.\ }}
 \newcommand{\prob}{{\bf Problem.\ }}

\newtheorem{note}{Note} 
\renewcommand{\thenote}{}
\newtheorem*{acka}{Acknowledgments}
\newtheorem{ack}{Acknowledgments}
\renewcommand{\theack}{}
\renewcommand{\appendixname}{\bf Appendix}
\renewcommand{\proof}{{\em Proof.\ }}

\hyphenation{
ap-pen-dix as-ymp-tot-ic at-trib-uted at-trib-ut-able
Bry-li-n-sky com-mu-ta-tion de-ge-ne-rate
de-riv-a-tive dis-trib-ute equi-vari-ant ex-tra-or-di-nary  
geo-met-ric griev-ance griev-ous grad-ed ho-lo-no-my ho-mo-thetic
in-fin-ite-ly in-fin-i-tes-i-mal Ha-rish Cha-n-dra mul-ti-plic-able 
non-euclid-ean non-iso-mor-phic non-smooth par-a-digm 
par-a-bol-ic pa-rab-o-loid pa-ram-e-trize phe-nom-e-non 
post-script pseu-do-dif-fer-en-tial pseu-do-fi-nite 
qua-drat-ics quad-ra-ture Han-kel rec-tan-gle semi-def-i-nite 
set-up wide-spread Euler-ian Feb-ru-ary Gauss-ian Grothen-dieck 
Hamil-ton-ian Her-mi-t-ian her-mi-t-ian Jan-u-ary 
Japan-ese Ka-shi-wa-ra Kor-te-weg Le-gendre No-vem-ber Rie-mann-ian 
Sep-tem-ber Za-mo-lo-d-chi-kov Kni-zh-nik quan-tum Op-dam
Mac-do-nald Ca-lo-ge-ro Su-ther-land Mo-ser 
Ol-sha-net-sky  Pe-re-lo-mov in-de-pen-dent ope-ra-tors 
cy-clo-to-mic ra-tio-nal de-gen-er-a-tion 
in-ter-est-ing de-for-ma-tions de-for-ma-tion pro-ce-dure 
fol-lows ope-ra-tors  pre-serve suf-fices ap-proach 
for-mu-las con-sider its com-ple-tion cor-re-spond-ing 
au-to-mor-phism be-cause pro-por-tional fi-nal-ly let-ting 
equi-v-a-lence ge-n-er-al-ized Mac-do-nald iden-ti-ties 
cor-re-s-pond sub-dia-grams par-ti-tion na-t-u-ral-ly 
or-dered stan-dard de-for-ma-tion ar-gu-ment com-bined 
sphe-r-i-cal rep-re-sen-ta-tions tri-go-no-me-t-ric
ge-n-er-al-ly speak-ing pri-m-it-ive ir-re-du-cible 
sum-ma-tion  rep-re-sen-ta-tives pro-por-ti-o-na-li-ty
ultra-sphe-ri-cal Ro-gers}

\def\ffor{\quad\hbox{ for }\quad}
\def\wwhen{\quad\hbox{ when }\quad}
\def\wwhere{\quad\hbox{ where }\quad}
\def\aand{\quad\hbox{ and }\quad}
\def\for{\  \hbox{ for } \ }
\def\iif{ \ \hbox{ if } \ }
\def\when{ \ \hbox{ when } \ }
\def\where{\  \hbox{ where } \ }
\def\and{\  \hbox{ and } \ }
\def\and{\  \hbox{ and } \ }
\def\oor{\  \hbox{ or } \ }
\def\proof{{\em Proof. \  }}

\def\equal{\stackrel{\,\mathbf{def}}{= \kern-3pt =}}

\def\la{\lambda}
\def\La{\Lambda}
\def\om{\omega}
\def\Om{\Omega}
\def\Th{\Theta}
\def\th{\theta}
\def\al{\alpha}
\def\be{\beta}
\def\ga{\gamma}
\def\ep{\epsilon}
\def\up{\upsilon}
\def\Up{\Upsilon}
\def\de{\delta}
\def\De{\Delta}
\def\ka{\kappa}
\def\kapp{\hbox{\bf \ae}}
\def\si{\sigma}
\def\Si{\Sigma}
\def\Ga{\Gamma}
\def\ze{\zeta}
\def\io{\iota}
\def\bio{b^\iota}
\def\aio{a^\iota}
\def\twio{\tilde{w}^\iota}
\def\hwio{\hat{w}^\iota}
\def\gio{\g^\iota}
\def\Bio{B^\iota}

\def\del{\delta}
\def\pa{\partial}
\def\vp{\varphi}
\def\ve{\varepsilon}
\def\inf{\infty}

\def\vph{\varphi}
\def\vps{\varpsi}
\def\vPh{\varPhi}
\def\vep{\varepsilon}
\def\vpi{{\varpi}}
\def\vth{{\vartheta}}
\def\vsi{{\varsigma}}
\def\vrh{{\varrho}}

\def\bph{\bar{\phi}}
\def\bsi{\bar{\si}}
\def\bvp{\bar{\varphi}}

\newcommand{\bS}{{\mathbf S}}
\newcommand{\bH}{{\mathbf H}}
\newcommand{\bF}{{\mathbf F}}
\newcommand{\bE}{{\mathbf E}}

\def\tal{\tilde{\alpha}}
\def\tbe{\tilde{\beta}}
\def\tde{\tilde{\delta}}
\def\tpi{\tilde{\pi}}
\def\txi{\tilde{\xi}}
\def\tPi{\tilde{\Pi}}
\def\tPhi{\tilde{\Phi}}
\def\tV{\tilde{V}}
\def\tJ{\tilde{J}}
\def\tla{\tilde{\lambda}}
\def\tga{\tilde{\gamma}}
\def\tGa{\tilde{\Gamma}}
\def\tvs{\tilde{{\varsigma}}}
\def\tu{\tilde{u}}
\def\tU{\tilde{U}}
\def\tw{\widetilde w}
\def\tW{\widetilde W}
\def\tB{\tilde B}
\def\tv{\tilde v}
\def\tV{\tilde V}
\def\tz{\tilde z}
\def\tb{\tilde b}
\def\ta{\tilde a}
\def\tih{\tilde h}
\def\trh{\tilde {\rho}}
\def\tx{\tilde x}
\def\tf{\tilde f}
\def\tg{\tilde g}
\def\tG{\tilde G}
\def\tk{\tilde k}
\def\tl{\tilde l}
\def\tL{\tilde L}
\def\tD{\tilde D}
\def\tR{\tilde R}
\def\tP{\tilde P}
\def\tH{\tilde H}
\def\tp{\tilde p}

\def\hH{\hat{H}}
\def\hh{\hat{h}}
\def\hR{\hat{R}}
\def\hY{\hat{Y}}
\def\hX{\hat{X}}
\def\hP{\hat{P}}
\def\hT{\hat{T}}
\def\hV{\hat{V}}
\def\hG{\hat{G}}
\def\hF{\hat{F}}
\def\hw{\widehat{w}}
\def\hW{\widehat{W}}
\def\hu{\hat{u}}
\def\hs{\hat{s}}
\def\hv{\hat{v}}
\def\hb{\hat{b}}
\def\hB{\widehat{B}}
\def\hze{\hat{\zeta}}
\def\hsi{\hat{\sigma}}
\def\hrh{\hat{\rho}}
\def\hth{\hat{\theta}}
\def\hy{\hat{y}}
\def\hx{\hat{x}}
\def\hz{\hat{z}}
\def\hg{\hat{g}}
\def\he{\hat{e}}
\def\hE{\widehat{E}}

\def\B{\mathbf{B}}
\def\I{\mathbf{I}}
\def\P{\mathbf{P}}
\def\G{\mathbf{G}}
\def\S{\mathbf{S}}
\def\F{\mathbf{F}}
\def\one{\mathbf{1}}
\def\Sn{\mathbf{S}_n}
\def\0{\mathbf{0}}
\def\H{\mathbf{H}}
\def\V{\mathbf{V}}

\def\f{\mathcal{F}}
\def\çF{\mathcal{F}}
\def\o{\mathcal{O}}
\def\t{\mathcal{T}}
\def\r{\mathcal{R}}
\def\l{\mathcal{L}}
\def\m{\mathcal{M}}
\def\k{\mathcal{K}}
\def\n{\mathcal{N}}
\def\d{\mathcal{D}}
\def\p{\mathcal{P}}
\def\cP{\mathcal{P}}
\def\a{\mathcal{A}}
\def\h{\mathcal{H}}
\def\c{\mathcal{C}}
\def\y{\mathcal{Y}}
\def\e{\mathcal{E}}
\def\v{\mathcal{V}}
\def\z{\mathcal{Z}}
\def\x{\mathcal{X}}
\def\s{\mathcal{S}}
\def\g{\mathcal{G}}
\def\u{\mathcal{U}}
\def\w{\mathcal{W}}
\def\i{\mathcal{I}}
\def\j{\mathcal{J}}
\def\b{\mathcal{B}}

\def\lan{\langle}
\def\llb{(\!(}
\def\ran{\rangle}
\def\rrb{)\!)}
 \def\dim{{\hbox{\rm dim}}_{\mathbb C}\,}
\def\lng{\hbox{\rm{\tiny lng}}}
\def\sht{\hbox{\rm{\tiny sht}}}
\def\sph{\hbox{\rm{\tiny sph}}}
\def\inv{\hbox{\rm{\tiny inv}}}

\def\br#1{\langle #1 \rangle}

\def\rank{\hbox{rank}}
\def\gl{\mathfrak{gl}_N}

\newcommand{\Aut}{\operatorname{Aut}}
\newcommand{\Hom}{\operatorname{Hom}}
\newcommand{\End}{\operatorname{End}}
\newcommand{\Ind}{\operatorname{Ind}}
\newcommand{\ad}{\operatorname{ad}}
\newcommand{\pr}{\operatorname{pr}}
\newcommand{\aweyl}{\tilde{\mathbb S}_n}
\newcommand{\hec}{{\mathcal H}^t_n}
\newcommand{\Func}{{\mathcal F}({\mathbb C}^n,{\mathcal H}^t_n)}
\newcommand{\tr}{\operatorname{tr}}
\newcommand{\Out}{\operatorname{Out}}
\newcommand{\Rad}{\operatorname{Rad}}
\newcommand{\Spec}{\operatorname{Spec}}
\newcommand{\id}{\operatorname{id}}
\newcommand{\Int}{\operatorname{Int}}
\newcommand{\ct} {\operatorname{ct}}

\newcommand{\rat}{{\mathbb Q}}
\newcommand{\real}{{\mathbb R}}
\newcommand{\cplx}{{\mathbb C}}
\newcommand{\zint}{{\mathbb Z}}

\newcommand{\sq}{\phantom{1}\hfill$\qed$}
\newcommand{\Rea}{\Re}
\newcommand{\Ima}{\Im}

\newcommand{\st}{\bowtie}
\newcommand{\modd}{\mbox{\,mod\,}}
\newcommand{\lr}{\langle}
\newcommand{\rr}{\rangle}
\newcommand{\eps}{\varepsilon}
\newcommand{\phk}{\phi^{(k)}}
\newcommand{\psk}{\psi^{(k)}}
\newcommand{\Res}{\mbox{Res}\;}
\newcommand{\sgn}{\mbox{sgn}}
\newcommand{\mn} {\left\{ \begin{array}{c}m\\
n\end{array}\right\}}

\def\TT{\mathfrak{T}}
\def\JJ{\mathfrak{J}}
\def\HH{\mathfrak{H}}
\def\FF{\mathfrak{F}}
\def\GG{\mathfrak{G}}
\def\CC{\mathfrak{C}}
\def\LL{\mathfrak{L}}

\def\BB{\mathfrak{B}}
\def\AA{\mathfrak{A}}
\def\ZZ{\mathfrak{Z}}
\def\HH{\hbox{${\mathcal H}$\kern-5.2pt${\mathcal H}$}}
\def\tHH{\widetilde{\HH\ }}

\font\smm=msbm10 at 12pt 
\def\symbol#1{\hbox{\smm #1}}
\def\lsmash{{\symbol n}}
\def\rsmash{{\symbol o}}
\def\#{\sharp}

\font\tenbf=cmbx10
\font\tenrm=cmr10
\font\tenit=cmti10
\font\ninebf=cmbx9
\font\ninerm=cmr9
\font\nineit=cmti9
\font\eightbf=cmbx8
\font\eightrm=cmr8
\font\eightit=cmti8
\font\sevenrm=cmr7
\font\sevenbf=cmbx7


\title [On Galois action in rigid DAHA modules]
{On Galois action in rigid DAHA modules}
\author[Ivan Cherednik]{Ivan Cherednik $^\dag$}

\begin{abstract}
Given an elliptic curve over a field $K$ of algebraic
numbers, we associate with it an action of the absolute Galois
group $G_K$ in the type $A_1$ rigid DAHA-modules at
roots of unity $q$ and over the rings $\Z[q^{1/4}]/(p^m)$
for sufficiently general prime $p$. We describe rigid modules
in characteristic zero and for such rings. The main examples of
rigid modules are generalized nonsymmetric Verlinde algebras;
their deformations for arbitrary $q$ are constructed in this
paper, which is of independent interest on its own. The Galois
action preserves the images of the elliptic braid group in the
groups of automorphisms of rigid modules over $\Z[q^{1/4}]/(p^m)$.
If they are finite in characteristic zero, then $G_K$ acts there 
and no reduction modulo $(p^m)$ is needed; we find all such cases.
In the case of $3$-dimensional DAHA-modules, these images are
quotients of equilateral triangle groups directly related to
the Livn\'e groups. Also, this paper can be viewed as a certain
extension of the DAHA theory of refined Jones polynomials of
torus knots (for $A_1$) to $G_K$.
\end{abstract}

\thanks{$^\dag$  February 2, 2014.
\ \ \ Partially supported by NSF grant
DMS--1101535, the Fulbright Program and the Simons Foundation}

\address[I. Cherednik]{Department of Mathematics, UNC
Chapel Hill, North Carolina 27599, USA\\
chered@email.unc.edu}

 \def\bysame{{\bf --- }}
 \def\~{{\bf --}}
\newcommand{\comment}[1]{}
\renewcommand{\tilde}{\widetilde}
\renewcommand{\hat}{\widehat}
\newcommand{\dagx}{\hbox{\tiny\mathversion{bold}$\dag$}}
\newcommand{\ddagx}{\hbox{\tiny\mathversion{bold}$\ddag$}}
\newtheorem{conjecture}[theorem]{Conjecture}
\newcommand*\toeq{
\raisebox{-0.15 em}{\,\ensuremath{
\xrightarrow{\raisebox{-0.3 em}{\ensuremath{\sim}}}}\,}
}

\comment{
"Deligne-Simpson Problem, Quivers and Galois action
on rigid DAHA-modules in rank one."

Given an elliptic curve over a field K of algebraic
numbers, we associate with it an action of the absolute Galois
group G_K in the rigid DAHA-modules of type A_1. The latter
are considered for roots of unity q and over the rings
Z[q]/(p^m) for sufficiently general prime p and any m.
We describe rigid modules in characteristic zero and
those for generic p.

They are directly connected with the Deligne-Simpson Problem
and therefore with the quivers which encode the conjugacy
classes for DSP (Crawley-Boevey and others). This was used
by Oblomkov and Stoica for DAHA of type C-check-C_1 (physics
links are expected here), but we need a more subtle theory
at roots of unity and in finite characteristic.

There are 3 types of rigid DAHA modules at roots of unity
for A_1: the generalized nonsymmetric Verlinde algebras
(the symmetric part of the one for t=q is the classical
Verlinde algebras for sl_2-hat, the so-called non-semisimple
Verlinde algebras (connected with the logarithmic CFT),
and the "massless" series essentially lifted from the rational DAHA.
The generalized Verlinde algebras can be deformed to
arbitrary q, which is of independent interest on its own.

The Galois action preserves the images of the elliptic braid
group in the groups of automorphisms of such modules. This can
be viewed as an extension of the classical theory of Tate modules
of elliptic curves. If the image of the elliptic braid group is
finite in characteristic 0, then G_K acts there and no reduction
modulo (p^m) is necessary; we will describe all such cases,
which is (morally) similar to the classical Schwarz's list in the
theory of hypergeometric functions.
}

\vskip -0.0cm
\par
{\centering
In memory of Andrei Zelevinsky
\medskip
\par}
\vskip -0.0cm
\maketitle
\smallskip

{\tiny
\centerline{{\bf MSC} (2010): 17B45, 20C08, 33D52, 12F12, 14G32,
30F10, 11G05, 20F65, 20H10, 51M10}
}
\smallskip

\vskip -0.0cm
\renewcommand{\baselinestretch}{1.2}
{\textmd
\tableofcontents
}
\renewcommand{\baselinestretch}{1.0}
\vfill\eject

\renewcommand{\natural}{\wr}

\setcounter{section}{-1}
\setcounter{equation}{0}
\section{\sc Introduction}
Given an elliptic curve defined over a field of algebraic numbers
$K$, we defined the corresponding action of the
{\em absolute Galois group\,}
$G_K=$Gal$(\overline{\Q}/K)$ in the {\em rigid DAHA-modules\,}
of type  $A_1$. The rigidity we use is actually that
in the corresponding {\em Deligne-Simpson Problem\,}.
These modules are considered at roots of unity $q$
and over the rings $\Z[q^{1/4}]/(p^m)$ for odd
prime $\,p\,$ such that gcd$(p,N)=1$ for the order $N$ of $q$
and for $\,m\in \N$.
{\em Generalized nonsymmetric Verlinde algebras\,} are the main
examples of
rigid modules; we construct their deformations
for any $q$, which is of interest in its own right.

The {\em elliptic braid group\,},
denoted by $\b_{1}$,  plays the key role in this paper;
$G_K$ acts in its images in the groups of automorphisms of
rigid modules considered over the rings $\Z[q^{1/4}]/(p^m)$.
We describe all rigid DAHA-modules in characteristic $0$;
they are defined over $\Z[q^{1/4}]$ and remain rigid modulo
$(p^m)$ for $p$ as above. We find all cases
when $\b_1$ has finite images in characteristic $0$; then
$G_K$ acts there without further reduction modulo $(p^m)$.
Also, we post a conjectural list of the cases when the images
of $\b_1$ are arithmetic discrete. A link
to equilateral {\em triangle groups\,} and
the Livn\'e lattices in $PU(2,1)$ is established.

Given an elliptic curve $E$ with a puncture at the
origin $\,0$, the elliptic braid group controls its covers
equivariant with respect to the reflection $x\mapsto -x$ in $E$.
Considering such covers is related to \cite{Bel} and
Grothendieck's program of {\em dessins d'enfants\,}. Let us also
mention here \cite{BL}
(the construction of the polylogarithm sheaf for $E\setminus 0$).

In contrast to the {\em Tate modules
of elliptic curves\,} based on the unramified covers of $E$,
covers of $E\setminus 0$ form a huge class.
A much more restrictive system of covers of $E\setminus 0$ can be
obtained via the DAHA modules at roots of unity.
They are similar to the unramified ones, though
the functoriality with respect to the roots of
unity $q$ is unclear. 
Even without this,
such an action of $G_K$ is interesting; for instance,
the DAHA-based theory of {\em refined Jones polynomials\,} of
torus knots from \cite{CJ} has a counterpart for $G_K$ 
(in type $A_1$).

\subsection{\bf Tate modules via rigidity}\label{sec:Tate}
Let $V_{2N}$ be the unique irreducible
nonzero module of the extended Weyl algebra
$\w_{2N}$ generated by the elements $X^{\pm1},Y^{\pm1}$
and $S$ with the relations
\begin{align*}
Y^{-1}X^{-1}YXq^{1/2}\!=\!1,\ \,
SXSX\!=\!1\!=\!SYSY,\  S^2\!=\!1,\ X^{2N}\!=\!1\!=\!Y^{2N},
\end{align*}
where  $q^{1/2}$ is assumed a primitive $(2\!N)$th root of unity.
Switching to the generators
$\,A=XS,\ \tilde{B}=q^{1/4}XSY,\  C=SY\,$, we interpret these
relations as the following (very special) case
of the multiplicative {\em Deligne-Simpson problem,  DSP\,}:
\begin{align*}
A^2=1=\tilde{B}^2=C^2,\, A\tilde{B}C=\tilde{S}\equal q^{1/4}S,\ \
A,\tilde{B},C,\tilde{S}\in \hbox{GL}(2\!N,\C),
\end{align*}
provided that the multiplicities of the eigenvalues $\pm 1$ coincide
(and equal to $N$) for $A,\tilde{B},C$ and are $N+1,N-1$ for the
eigenvalues $q^{1/4},-q^{-1/4}$ of $\tilde{S}$.
This {\em DSP\,} is {\em rigid\,}; it
has a unique solution (up to a conjugation) due to the
uniqueness of $V_{2N}$. Up to a conjugation,
the matrices $A,\tilde{B},C,\tilde{S}$
can be assumed here over $\Q(q^{1/4}).$
\smallskip

For an elliptic curve
$E_{\C}$, let $\,o,o_1,o_2,o_3\in P^1_{\C}\,$ be the images of
$\,0\in E_{\C}\,$ and the remaining
points $\{0_1,0_2,0_3\}$ 
of the $2$nd order in $E_{\C}$ upon 
the identification
$E/\{s\}\toeq P^1$, where $s: z\mapsto -z$ in $E_{\C}\ni z$.

Finite quotients of the group
$\b_1=\,<\!A,\tilde{B},C,\tilde{S}\!>$ result in
(connected) Galois covers of $P^1_{\C}$ ramified 
at $\{o,o_1,o_2,o_3\}$, which can be naturally 
considered as covers of $E_{\C}$ ramified only at 
$0\,$ associated with the
corresponding quotients of\,
$\!<\!\tilde{X}\!\!=\!q^{-1/4}X,\tilde{Y}\!\!=
\!q^{1/4}Y,\tilde{S}^2\!>\!$.
We apply this to the image $\tilde{\mathfrak{B}}$ of 
$\b_1=\,<\!A,\tilde{B},C,\tilde{S}\!>$ in
GL$(2\!N,\Q(q^{1/4}))$.
Dividing $\tilde{\mathfrak{B}}$ by the
center, the corresponding {\em unramified\,} covering 
$\pi_{2N}:E\simeq\tilde{E}\to E$ is the multiplication 
by $2\!N$; 
its kernel $\tilde{E}_{2N}\simeq\Z_{2N}^2$ is the Galois group.

\comment{
The commutator subgroup of  $\tilde{\mathfrak{B}}$ corresponds 
to the extension of $\C(\tilde{E})$ by the function 
$f^{\!\frac{1}{2\!N}}$ with $(f)\!\subset\!\tilde{E}_{2N}$ such that
$\{f^ef^{-1},e\!\in\! \tilde{E}_{2N}\}$ is from 
$Z^1(\tilde{E}_{2N},(\C(\tilde{E})^*)^{2N})$ and maps 
to the $2$-cocycle from $H^2(\tilde{E}_{2N},\mu_{2N})$ 
canonically associated with $\w_{2N}$. 
Here $\mu_{2N}=\{\ze:\ze^{2N}\!=\!1\}$; 
we use the Galois cohomology, Theorem Hilbert-90 and Tsen's theorem.
}

Let $q^{\frac{1}{2}}\!=\!e^{\frac{\pi\imath}{N}},\  
e_1,e_2\!\in\! \tilde{E}_{2N}$ be the generators corresponding
to $X,Y$. The group $<\!X\!,Y\!,q^{1/2}\!>$ corresponds to 
the extension of $\C(\tilde{E})$ by the function 
$f^{\!\frac{1}{2\!N}}$ for $f$ with the divisor 
$(f)\!=\!\sum_{i=1}^{2\!N\!-\!1}\!D^{\,e_2^i}\!-\!(2\!N\!-1)\!D$, 
where
$D\!=\!\sum_{i=0}^{2\!N\!-\!1} e_1^i$ and $e^i_j=ie_j$.
The action of $\tilde{E}_{2N}$ is pointwise.
\smallskip

Assume that $E$ and $0$ are defined over
a field $K\subset \overline{\Q}$,
as well as $o\in P^1$ and the set $\{o_1,o_2,o_3\}$.
The rigidity of {\em DSP\,} above gives that
\vskip 0.2cm

{\em
$(I)$ the group \,Gal$\bigl(\overline{\Q}/K)\bigr)\,$
acts in the group $\tilde{\mathfrak{B}}$ by automorphisms;

$(II)$ Gal$\bigl(\overline{\Q}/K(q^{1/4})\bigr)$
acts via conjugations by matrices in
\,GL$(V_{2N})$.
}
\vskip 0.2cm

\comment{
Switching to the coverings over $E$ and dividing by
the center, the action of Gal$(\overline{\Q}/K)$ is that in
the group $E_{2N}$ of the points of $(2\!N)$th order of $E$.
}
Recall that the classical {\em Tate module} is the projective
limit $\varprojlim E_{\ell^m}$ as
$m\to\infty$
for a prime $\ell$\,; it has a standard action of $G_K$. 
Thus
our constructions for $N=\ell^m$ extends such modules.

Note that $V_{2N}\simeq H^0(\o(D))$ for a divisor $D$ of
degree $2\!N$ in $E$ and $q=e^{\pi\imath/N}$,
where $\mathfrak{B}$ acts due
to $\o(D^e)\simeq \o(D)$ 
(the theory of Kummer-Weil pairing).
Using $A,\tilde{B},C,\tilde{S}$ 
establishes the isomorphisms
$\o(D^e)\simeq \o(D)$ via the trivialization of $\o(D)$ 
at $E\setminus 0$; here $\tilde{S}^2=q^{1/2}$
is the monodromy of $\o(D)$ at $0\,$.
\smallskip

Claims $(I,II)$ are extended in this paper to any rigid
DAHA-modules of type $A_1$ at roots of unity (and in finite
characteristic); these are the boundary case
of our construction from $(\al^\bullet)$ in Section
\ref{sec:Boundary}. However the compatibility 
for different roots of unity is missing in such a generality.
Taking $N=\ell^m$ for odd prime numbers $\ell\,$ together with
making
$\Z[q^{1/4}]/(\ell^m)$ the ring of definition is a possibility here,
but this is beyond the present paper. See also an example
in (\ref{NcheckNxx}) below.
\smallskip

We note that the approach from \cite{Bel} (and other
works) devoted to the (Regular) Inverse Galois Problem
is mainly based on {\em rigid triples}, which are
the triples $a,b,c\in G$ generating $G$ and
satisfying $abc=1$. They are assumed from
the corresponding (given) conjugacy classes in $G$. 
This kind of rigidity means the
uniqueness of such a triple up to a simultaneous
conjugation in $G$. We need $4$ points in $P^1$ and the 
so-called {\em linear rigidity\,} (in matrices) based on
Katz' theory of rigid systems. Enriching Belyi's
approach by means of the multiplicative
{\em DSP\,} (M. Dettweiler and others) was an important
development.

\subsection{\bf The main results of the paper}
We begin the paper with adjusting the classification results from
\cite{C101}, mainly Sections 2.8-9 there, to our objectives.
For a primitive  $N$th root of unity $q$,
we need DAHA-modules over the rings $\Z[q^{1/4}]/(p^m)$ for
prime numbers $p$ such that gcd$(p,2\!N)=1$ and any $m\in \N$.
The gcd-constraint makes such a  modification of \cite{C101}
relatively straightforward. We also need to address
the case when $q^{1/2}$ is not primitive of order $2\!N$
for odd $N$.

$(i)$ Theorem \ref{psl-two-inv} describes
all {\em rigid DAHA-modules\,} $V$ in characteristic $0$
for primitive $q^{1/2}$ of order $2\!N$.
The {\em dim-rigidity\,}, when (by definition)
only finitely many
irreducible DAHA-modules exist of dimension equal to dim$V$,
is sufficient for $A_1$.
Some of them are {\em unirigid\,}, when the dimension
alone entirely determines their isomorphism classes;
they automatically have the projective
action of $PSL_2(\Z)$.

$(ii)$ Proposition \ref{ext-verlinde} generalizes
Section 2.10.5 of \cite{C101} and provides
$q$\~deformations of the {\em generalized
nonsymmetric Verlinde algebras\,}
(in nonsymmetric polynomials and for
arbitrary admissible $t$). The classical Verlinde algebra
of $\hat{\mathfrak{sl}}_2$ is for $t=q$ and upon
the symmetrization.

$(iii)$ Theorem \ref{FINBIMAGE} describes all
dim-rigid modules with finite images of $\b_1$.
We analyze the positivity of the hermitian inner product
at (all) primitive roots of unity $q^{1/2}$.
There are two series of dimension $2$ and $4$ for
all $N$ and two
exceptional cases $N=6,10$ of dim$=8$. We also provide
a conjectural list of arithmetic discrete $Image(\b_1)$.

$(iv)$ The finiteness/discreteness of this kind is a very classical
topic (cf. Schwarz' table for hypergeometric
functions). If the image of $\b_1$ is finite, then
{\em Riemann Existence Theorem\,} provides an action of
$G_K$ in $\, Image(\b_1)$
without modular considerations. Using DAHA in finite
characteristic is needed 
when $\, Image(\b_1)$ is not finite.

$(v)$ Theorem \ref{GALACTION} defines the action of
$G_K$ in the images of
$\b_1$ in rigid modules over $\Z[q^{1/4}]/(p^m)$.
The basic examples are $(a)$ the Tate modules
(see $(I,II)$ above), $(b)$ subgroups
of the Livn\'e groups considered
modulo $(p^m)$ and
$(c)$ the Verlinde algebra ($t=q$)
over $\Z[q^{1/4}]/(p^m).$

$(vi)$ Importantly, we can define the action from
$(c)$ above only
via the nonsymmetric theory. Using the group $\b_1$
and {\em DSP\,} is possible only {\em before} the
symmetrization of the modules $V$. However
after the action of $G_K$
is introduced in the whole $V$, it can be normalized by
the condition $T\mapsto T$. Then the symmetric part of $V$,
defined as $\{v\in V\, \mid\, Tv=t^{1/2}v\}$, becomes
invariant under this action.

$(vii)$ The rigidity in our paper is closely
related to that in \cite{ObS}, where {\em DSP\,}
was used for the classification of rigid
irreducible DAHA modules for $C^\vee C_1$
(apart from the roots of unity), which system is
a natural setting here. Theorem \ref{GALACTION}
generally can be extended to $C^\vee C_1$, though
there is no complete theory at roots
of unity for $C^\vee C_1$ so far.

$(viii)$ The case of $A_1$ is of special interest due
to its connection with covers of elliptic curves ramified
at one point. Also, the images of the elliptic braid group
$\b_1$ in rigid DAHA-modules at roots of unity for $A_1$
are directly related to the {\em triangle groups\,}
of type $(n,n,n;n)$; see Proposition \ref{TRIANGB}.
We establish a connection with the Livn\'e groups (dim$V=3$);
see \cite{Par1}. Links to \cite{LLM} and other papers on finite
quotients of triangle groups are expected.

\comment{
{\bf Acknowledgements.}
The author thanks David Kazhdan, Maxim Kontsevich, Nikita
Nekrasov and Yan Soibelman for useful discussions and Hebrew
University, IHES and SCGP (where this work was reported)
for invitations and hospitality.
\smallskip
}

\setcounter{equation}{0}
\section{\sc DAHA and rigidity}
\subsection{\bf Main definitions}
We will consider only the case of $A_1$ in this paper.
Let $\al=\al_1$, $s=s_1$, and $\om=\om_1$ be the fundamental weight;
then $\alpha=2\omega$ and
$\rho=\om$. The extended affine Weyl group
$\widehat{W}=<\!s,\om\!>$
is a free group generated by the involutions $\,s\,$ and
$\,\pi\equal \om s$.

The generators of double affine Hecke algebra
$\HH=\HH_{q^{1/2},t^{1/2}}$ are
$$
Y=Y_{\om_1}=\pi T,\ \,T=T_1,\ \,X=X_{\om_1}
$$
subject to the quadratic relation $\,(T-t^{1/2})(T+t^{-1/2})=0\,$
and the cross-relations:
\begin{align}\label{dahaone}
&TXT=X^{-1},\ T^{-1}YT^{-1}=Y^{-1},\ Y^{-1}X^{-1}YXT^2q^{1/2}=1.
\end{align}
Using $\pi=YT^{-1}$, the second relation becomes $\pi^2=1$.
This algebra is defined over
$$
\Z_{q;t}\equal\Z[q^{\pm 1/2},t^{\pm 1/2}].
$$
It is important that $\HH_{q^{1/2},1}$ is the
{\em extended Weyl algebra} generated by $X^{\pm1},Y^{\pm1},s\, $
subject to the relation $Y^{-1}X^{-1}YXq^{1/2}=1$ for
the involution $\,s\,=\,T(t^{1/2}=1)\,$ such that
$sXs= X^{-1}$ and $sYs=Y^{-1}$.
\smallskip

{\sf Automorphisms.}
The following maps can be extended to automorphisms
of $\HH$:
\begin{align}\label{tau+def}
&\tau_+(X)=X,\ \tau_+(T)=T,\ \tau_+(Y)=q^{-1/4}XY,\
\tau_+(\pi)=q^{-1/4}X\pi,\\
\label{tau-def}
&\tau_-(Y)=Y,\ \,\tau_-(T)=T,\ \,\tau_-(X)=q^{1/4}YX,\ \
\tau_-(\pi)\,=\,\pi.
\end{align}
They require adding $q^{\pm 1/4}$ to the ring of definition
of $\HH\,$.

The generalized Fourier transform corresponds to the
following automorphism of $\HH$ (it is not an
involution)\,:
\begin{align}
\si(X)= Y^{-1},\ \si(T)=T,\ &\si(Y)=q^{-1/2}Y^{-1}XY=XT^2,\
\si(\pi)=XT,
\notag\\
\si\ &=\  \tau_+\tau_-^{-1}\tau_+\ =\ \tau_-^{-1}\tau_+\tau_-^{-1}.
\label{tautautau}
\end{align}
The last relation is the defining one
of the {\em projective} $PSL_2(\Z)$. Therefor the latter
group acts in $\HH$
by outer automorphisms.
Check that $\,\si\tau_+=\tau_-^{-1}\si,\ \,
\si\tau_+^{-1}=\tau_-\si.$

Due to the group nature of the definition of $\HH$,
we have the inversion anti-involution $\HH\ni
H\mapsto H^\ast$\,:
$$
X^\ast=X^{-1},\, Y^\ast=Y^{-1},\, T^\ast=T^{-1},\,
(q^{1/4})^\ast=
q^{-1/4},\, (t^{1/2})^\ast=t^{-1/2}.
$$
It commutes with all automorphisms of $\HH$.
\smallskip

{\sf Polynomial representation.}
It is defined as
$\mathscr{X}=\Z_{q;t}[X^{\pm 1}]$ over the ring of
definition of $\HH\,$; recall that
$\Z_{q;t}=\Z[q^{\pm 1/2},t^{\pm 1/2}]\,$.
The operator naturally $X$ acts by
the multiplication. The formulas for the
other generators are
\begin{eqnarray*}
T=t^{1/2}s+\frac{t^{1/2}-t^{-1/2}}{X^{2}-1}
\circ (s-1),\ \, Y=\pi T
\end{eqnarray*}
in terms of the (multiplicative) reflection $s(X^n)=X^{-n}$
and $\pi(X^n)=q^{n/2}X^{-n}$ for $n\in\Z$; here $\circ$
stays for the composition.
This module is induced from the one-dimensional
representation $\{T=t^{1/2}=Y\}$ of
the affine Hecke subalgebra $\mathcal{H}_Y=\lan Y, T\ran$.

Let $\lr f \rr$
be the constant term of $f\in \Q(q^{1/2},t^{1/2})[X,X^{-1}]$,
\begin{align}\label{mutildemuone}
&\mu(X;q,t)\equal\prod_{j=0}^\infty
\frac{(1-q^jX^2)(1-q^{j+1}X^{-2})}
{(1-tq^jX^2)(1-tq^{j+1}X^{-2})}.\\
\label{mucontone}
&\mu_\circ\equal\mu/\lr\mu\rr= 1+\frac{t-1}{1-qt}(X^2+qX^{-2})
+\ldots\ ,\\
&\where \lr \mu \rr\ =\ \prod_{j=1}^\infty \frac{(1-tq^j)^2}
{(1-t^2q^j)(1-q^j)}.\notag
\end{align}

We define the (symmetric) inner product in $\mathscr{X}$:
\begin{align}\label{innerp}
&\lan f,g\ran\equal \lan fg^*\mu_\circ\ran,
\hbox{\ \,satisfying\ \ }
\lan f,H(g)\ran =\lan H^*(f),g\ran,
\end{align}
where $\ast$ in $\mathscr{X}$ is the restriction of $\ast$
from $\HH$ and we add the denominators of the coefficients of
$\mu_\circ$ to the ring of definition of $\mathscr{X}$.
We note that the automorphism $\tau_-$ naturally acts in
$\mathscr{X}$.
\smallskip

{\sf The e-polynomials.}
Their coefficients are from the ring $\Z^{loc}_{q;t}$, which we
define as the localization of $\Z_{q;t}$
by the expressions in the form $(1-t^aq^b)$ for integral
$a,b$ assuming that $t^a q^b\neq 1$.
For generic $q,t$, the following relations fix
$e_n\, (n\in \Z)$ uniquely  up to proportionality:
\begin{align}\label{nonsymp}
&Ye_{n} = (t^{1/2}q^{|n|/2})^{-sign(n)}\,e_{n}\for n\in \Z,\,
\text{where}\\
&sign(0)\equal\, -1, \text{\  i.e. $\,0\,$ is treated as negative.} 
\notag
\end{align}
We normalize $\{e_n\}$ as follows:
$\,e_{n}=X^{n}+\text{ ``lower terms'' }\, ,$
where by ``lower terms'', we mean
polynomials in terms of $X^{\pm m}$ as $|m|<n$
and, additionally,  $X^{|n|}$ for negative $n$.

The $e_{n} \,(n\in \Z)$ are called the
{\em nonsymmetric Macdonald polynomials} or simply
$e$\~polynomials; their coefficients actually are
given only in terms of $q,t$ (i.e. no
$q^{1/2},t^{1/2}$ are needed). One has:
\begin{align}\label{fnorme}
&\lr e_m e_n^*\mu_\circ\rr\ =\ \de_{mn}
\prod_{0<j<|\tilde{n}|}\frac{(1-q^j)(1-q^j t^2)}
{(1-q^j t)(1-q^{j}t)},\\
&|\tilde{n}|=|n|+1 \hbox{\, if\, } n\leq 0 \and
|\tilde{n}|=|n| \hbox{\, if\, } n>0, \notag\\
&e_n(t^{-1/2})\ =\ t^{-|n|/2}\prod_{0<j<|\tilde{n}|}
\frac{1-q^jt^2}{1-q^jt} \for n\in \Z. \label{epolval}
\end{align}

\subsection{\bf Finite-dimensional modules}\label{sec:rootsofunity}
We first assume that $q$ is generic. The classification
of finite-dimensional irreducible modules for such $q$ is as
follows.

Let us introduce the following automorphisms
of $\HH$:
\begin{align}\label{fiota}
\iota:&\, T\mapsto -T,\, X\mapsto X,\, Y\mapsto Y,
\ q^{1/2}\mapsto q^{1/2},\, t^{1/2}\mapsto t^{-1/2},\\
\label{fvarsx}
\varsigma_x:&\, T\mapsto T,\, X\mapsto -X,\, Y\mapsto Y,
\ q^{1/2}\mapsto q^{1/2},\, t^{1/2}\mapsto t^{1/2}.\\
\label{fvarsy}
\varsigma_y:&\, T\mapsto T,\, X\mapsto X,\, Y\mapsto -Y,
\ q^{1/2}\mapsto q^{1/2},\, t^{1/2}\mapsto t^{1/2}.
\end{align}
The automorphism $\varsigma_x$ obviously preserves $\mathscr{X}$;
the polynomial $e_m$ is multiplied by $(-1)^m$.
Note that $\,\iota$ changes $\HH_{q^{1/2},t^{1/2}}$ to
$\HH_{q^{1/2},t^{-1/2}}\simeq\HH_{q^{1/2},-t^{1/2}}$ and
commutes with the projective action of
$PSL_2(\Z)$, i.e. with the $\tau_{\pm}$ defined above.
The latter act in an obvious way on
the group generated by $\,\varsigma_x, \varsigma_y$, isomorphic
to $\Z_2^2\,$:
\begin{align}\label{pslvarsigma}
&\tau_+(\varsigma_x)=\varsigma_x,\ \tau_+(\varsigma_y)=
\varsigma_x\varsigma_y,\ \
\tau_-(\varsigma_y)=\varsigma_y,\ \tau_-(\varsigma_x)=
\varsigma_x\varsigma_y.
\end{align}

The following theorem is
from \cite{C101}, Theorem 2.8.1 and Theorem 2.8.2 (see also
the references therein). We set $t=q^k$ for $k\in \C$, assuming
that $q$ is generic in the following theorem. By (ir)reducibility
and $Y$\~semisimplicity we mean those over the field of
fractions of $\Z_{q;t}$, though all $\HH$\~modules discussed below
are over $\Z_{q;t}$.

\begin{theorem} \label{tnegk} 
(i) For generic $q$, irreducible
finite-dimensional $\HH$\~modules
over the field of fractions $\Q(q^{1/2},t^{1/2})$ of $\Z_{q;t}$ are
either $Y$\~spherical,
i.e. quotients of $\mathscr{X},$ or their images under the
automorphism $\,\varsigma_y$, or the images of
$Y$\~spherical modules defined for DAHA\, $\HH_{q^{1/2},t^{-1/2}}$
with the parameters $q^{1/2},t^{-1/2}$ upon the application of
$\,\iota$ or $ \iota\varsigma_y$.

The polynomial representation $\mathscr{X}$ considered
over $\Q(q^{1/2},t^{1/2})$
is not $Y$\~semisimple
if and only if $t=q^{-n}$ for $n\in \N$.
It is reducible if and
only if $\,t=\zeta q^{-n/2}$ for $\zeta=\pm 1$ and
$n\in \N$, where $\zeta=-1$ when
$n\in 2\Z_+$. Therefore if $\mathscr{X}$ is reducible, then
this module is $Y$\~semisimple;
the series $\,\mu_\circ$, the pairing $\,\lr f,g\rr$ on
$\mathscr{X}$
and all polynomials $e_m\, (m\in \Z)\,$ are well defined
in this case.

(ii) For $n\in\Z_+$,
let us fix $q^{1/2}$. We set $k=-1/2-n$ and $\,t=(q^{1/2})^{-1-2n}$.
The scalar squares $\lr e_m,e_m\rr$ are nonzero exactly at
$\{m\}'\equal\{-2n\le m \le 2n+1\}.$
The radical $Rad_0$ of the pairing $\lr f,g \rr$ is
$\,\oplus_{m\not\in \{m\}'\,} \BC e_m.$ Furthermore
$Rad_0=(e_{-2n-1})$ as an ideal in $\mathscr{X},$
and the $\HH\, $\~module
$\,V_{4n+2}\equal\mathscr{X}/Rad_0$
considered over $\Z^{loc}_{q;t}$
is the greatest
 $\Z^{loc}_{q;t}$\~free finite-dimensional quotient of
$\mathscr{X}.$

The coefficients of the
polynomial $e_{-2n-1}\,$ belong to $\Z[q^{\pm 1/2}]$ and the
module $W=W_{4n+2}\equal\mathscr{X}/(e_{-2n-1})$ is
well defined over $\Z_{q;t}=\Z[q^{\pm 1/4}]$
(i.e. without the localization with respect to $(1-q^{j/2})$
for $j\in \N$). It is a direct sum of the following two
$\Z_{q;t}$\~free non-isomorphic $\HH\,$\~submodules of
dimension $2n+1=2|k|$:
\begin{align}\label{epmte}
&V^\pm_{2n+1}\cong \mathscr{X}/\,(e_{n+1}\mp t^{-1/2}e_{-n});\\
&\hbox{see \,Theorem 2.8.1\, from\, \cite{C101}}.
\end{align}
These modules are orthogonal to each other with respect to
$\lr\, ,\, \rr$ restricted to $W$; also,
$\, V^-_{2n+1}\cong \varsigma_x(V^+_{2n+1}).$

The module $V^+_{2n+1}$ is $\tau_{\pm}$\~invariant.
The modules $\,\varsigma_x^{\ep}\varsigma_y^{\de}(V^+_{2n+1})$
are pairwise non-isomorphic for different
choices of $\ep,\de=0,1$; they are transformed by $\tau_{+},$
$\tau_{-}$ according to formulas $(\ref{pslvarsigma})$.

(iii) The remaining case of reducibility of
$\mathscr{X}$  is when
$\,t=-(q^{1/2})^{-n}$ for
$n\in \N$ (i.e. for $k'\equal-n/2$ and $\zeta=-1\,)$, where
we set $\,t=\zeta\,(q^{1/2})^{2k'}$.
Then the maximal $\Z_{q;t}$\~free finite-dimensional
quotient of $\mathscr{X}$ is irreducible
of even dimension $2n=4|k'|$. It is
isomorphic to $V'_{2n}\equal \mathscr{X}/(e_{-n}),$ where
$e_{-n}\in \mathscr{X}$ and
the ideal $(e_{-n})$ is the radical $Rad_0$ of the pairing
$\lr f,g\rr$. This module is linearly generated over
$\Z^{loc}_{q;t}$
by $e_m$ for $1-n\le m \le n$.
The squares $\lr e_m,e_m\rr$ are nonzero exactly at
$m=0,1,-1,2,\ldots, 1-n,n.$
The automorphisms $\,\varsigma_x$, $\,\varsigma_y$ and
$\,\tau_{\pm}$ do not
change the isomorphism class of $V'_{2n}$.\sq
\end{theorem}

The following is the explicit description of
$V^\pm_{2n+1}$ over $\Z^{loc}_{q;t}$ from \cite{C101}:
\begin{align}\label{epmtex}
V^\pm_{2n+1}& =
\oplus_{m=1}^{2n+1}\,\Z^{loc}_{q;t}\,(\ep_{m}\!\pm\!\ep_{-2n-1+m})
\hbox{ mod }
Rad_0, \ \, \ep_m\!\equal\! \ep_m/\ep_m(t^{-1/2}),\notag\\
&e_{n+1}\!\pm\! t^{-\frac{1}{2}}e_{-n}\!=\!
e_{n+1}(t^{-\frac{1}{2}})(\ep_{n+1}\!\pm\!\ep_{-n})\!=\!
X^{-n}\!\!\!\prod_{m=-n}^{n}\!\!(X\!\pm\! q^{1/4+m/2}).
\end{align}
Here $e_m(t^{-1/2})$ is from (\ref{epolval}). Note
that $e_{m}$ and  $\ep_{-2n-1+m}$ have coinciding
$Y$\~eigenvalues, namely $t^{1/2}q^{m/2}=t^{-1/2}q^{m/2-2n-1}$;
see (\ref{nonsymp}).
\smallskip

{\sf Roots of unity.}
We assume now that $q$ is a primitive root of unity of degree
$N\ge 1$
and also {\em always} pick $q^{1/2}$  in primitive $(2\!N)$th
roots of unity. We will actually allow taking
$q^{1/2}$ in $N$th roots of unity for odd $N$, but
will incorporate such a case using
{\em squaring $q,t$\,} or via {\em Little DAHA}
to be considered below.

We set $\,t=(q^{1/2})^{2k}$ for $k\in \Z/2$.
Then $\Z_{q;t}$ becomes $\Z_{\dot{q}}\equal \Z[\dot{q}]$, where
$\dot{q}=q^{1/2}$ for even $2k$ or $\dot{q}=q^{1/4}$ for odd $2k$. 
Accordingly,  $\Z^{loc}_{q;t}=\Z^{(2N)}_{\dot{q}}$, where the
latter is the localization of $\Z_{\dot{q}}$ with respect to
$2\!N$.  Thus coefficients
of $e$\~polynomials, if well defined, are
from $\Z[q]^{(N)}$ or $\Z[q^{1/2}]^{(2N)}$ respectively
for even and odd $2k$.

As above, we call a $\Z_{\dot{q}}$\~free
$\HH$\~module $V$ irreducible
if it is irreducible over the field of fractions
of $\Z_{\dot{q}}$.
In all examples below such irreducibility implies
{\sf loc-irreducibility\,}, which means that $V$ remains
irreducible over
$\Z_{\dot{q}}/\mathfrak{p}$ for any prime ideal $\mathfrak{p}$
in $\Z_{\dot{q}}$ over $(p)\subset \Z$ for prime $p$ such
that gcd$(p,2\!N)=1$. The latter conditions gives that any
$\HH$\~submodule of $V_{\mathfrak{p}}^m\equal
V\otimes_{\Z_{\dot{q}}} (\Z_{\dot{q}}/\mathfrak{p}^m)$ for $m\in \N$
has the form $a V_{\mathfrak{p}}^m$ for proper
$a\in \Z_{\dot{q}}/\mathfrak{p}^m$ (assuming\, gcd$(p,2\!N)=1$).
Accordingly, $V$ is called {\sf loc-semisimple\,} if $Y$ is
semisimple over all such $\Z_{\dot{q}}/\mathfrak{p}^m$.

Following mainly Theorems 2.9.3, 2.9.8 from \cite{C101},
we will begin with the $\HH$\~modules
$\,V^{(2)}_{2N}\equal\mathscr{X}/(X^N+X^{-N})$,
$$
V^{(-2)}_{4N}\equal\mathscr{X}/(X^{2N}+X^{-2N}-2),\
V^{(2)}_{4N}\equal\mathscr{X}/(X^{2N}+X^{-2N}+2).
$$
The module $\,V^{(-2)}_{4N}$ will be used in the case of
$k\in\Z$ in the next theorem. The module $\,V^{(2)}_{2N}$,
which is a natural quotient of $V^{(2)}_{4N}$,
will be needed below for $\,k\in 1/2+\Z.$
Respectively, $t^{N}=q^{kN}=\pm 1$ in these
two cases.

The next theorem is part of the general description
of all irreducible $\HH$\~modules at roots of unity
from \cite{C101}; we aim at the {\em dim-rigidity\,}
to be considered below.

\begin{theorem}\label{perfectposit}
(i) Let  $2k\in \Z$ and $0<2k<N$.
Then $\,V^{(-2)}_{4N}$ for $k\in \Z$
or, correspondingly, $\,V^{(2)}_{2N}$ for $k\in 1/2+\Z$ has a
unique $\Z_{\dot{q}}$\~free irreducible nonzero quotient, which is
$\,V_{2N-4k}=\mathscr{X}/(e_{-n})$ of dimension $2N-4k$ for
$\,n\equal N-2k.$ It is also loc-irreducible.
Its $Y$\~spectrum is simple: $\,V_{2N-4k}=\oplus_{m=2k-N+1}^{N-2k}
\Z_{\dot{q}}e'_m$, where the eigenvectors $e_m'$ are the
images of $e_{m}$ for $-n+1\le m\le n,$
which are all well defined. The spectrum remains simple over
the rings $\Z_{\dot{q}}/\mathfrak{p}^m$ (for gcd$(p,2\!N)=1$ and
any $\,m\in\N$).
Also,
\begin{align}\label{TeY}
T(e'_{-n})\ =\ -t^{-1/2}e'_{-n}\ =\ Y(e'_{-n}).
\end{align}
The isomorphism class of this module is invariant under
the action of  $\,\varsigma_x$ and  $\,\varsigma_y$.

These statements can be extended to $k=0$, i.e for $t^{1/2}=1$,
when  $\,X^N-X^{-N}\,$ is taken instead of $e_{-n}.$

(ii) Let $\,k\equal -1/2-n$ for integral $\,0\le n< (N-1)/2.$ Then
$\,V^{(2)}_{2N}=\mathscr{X}/(X^N+X^{-N})\,$ has two
non-isomorphic $\Z_{\dot{q}}$\~free
and $Y$\~semisimple quotients, namely those
from (\ref{epmte}), which are
\begin{align}\label{epmte1}
&V^\pm_{2n+1}\cong \mathscr{X}/(e_{n+1}\mp t^{-1/2}e_{-n}).
\end{align}
Here $e_{n+1},e_{-n}$ are well defined over $\Z_{\dot{q}}$.
These modules are also loc-irreducible and loc-semisimple.

The binomials $(iX)^{n+1}\pm (iX)^{-n}$, where $i=\sqrt{-1}$,
must be taken instead of
$\,e_{n+1}\pm t^{-1/2}e_{-n}$ to extend these statements to
the case when $\,N=2n+1,\, t^{1/2}=i$.

The kernel $(e_{-2n-1})$ of the homomorphism from
$V^{(2)}_{2N}$ to the direct sum of  $V^+_{2n+1}$ and
$V^-_{2n+1}$ is a free $\Z_{\dot{q}}$\~module of dimension
$2N-4|k|$ and is isomorphic to $V_{2N-4|k|}$
from $(i)$ with $t^{1/2}$ replaced by $t^{-1/2}$  under
the action of involution $\,\iota$ (sending
$\,T\mapsto -T$; see above). The vector $\,e=e_{-2n-1}$ satisfies
the relations $\,T(e)\ =\ -t^{-1/2}e\, $ and  $\,Y(e)= t^{-1/2}e.$
\sq
\end{theorem}

In Part $(ii)$, the polynomials $\ep_m$
are well defined for $\,-2n\le m\le 2n+1$
over $\Z_{\dot{q}}^{(N)}$
and the linear decompositions of $V^\pm_{2n+1}$ in terms
of $\ep'_m\pm \ep'_{m-2n-1}$ from (\ref{epmtex}) holds over
this ring.
Also $\lan \ep_m,\ep_m\ran=\lan \ep_{m-2n+1},\ep_{m-2n+1}\ran$
in this range.

Both families of modules, $\,V=V_{2N-4k}$
and $V=V^\pm_{2n+1}$, have the invariant hermitian inner products
$\lan f\,,\, g\ran$, which can be defined directly by formulas
(\ref{fnorme}) for the images of $e_m$ in the proper range of $m$. 
Here we extend $\Z_{\dot{q}}$ to  $\Z_{\dot{q}}^{(N)}$.  One has
$\lan Hf\,,\, g\ran=\lan f\,,\, H^{-1}g\ran$ for
$H=X,Y,T$ and $f,g\in V$; see (\ref{innerp}). To justify the
existence of such a form,
it suffices to observe that
$\lan e\,,\, e\ran=0$  for $\,e=e_{-n}$ in $(i)$ and
correspondingly for $\,e=e_{-2n-1}$ in $(ii)$
(and then use the intertwiners).

\subsection{\bf Rigidity at roots of unity}
We will begin with the following exact sequences
from \cite{C101}, which clarify
Theorem \ref{perfectposit}.

\begin{proposition}\label{GREATSI}
Let $k\in \Z/2, |k|<N/2.$
There exist two exact sequences over $\Z_{\dot{q}}$ for $k>0$:
\begin{align}\label{eseq1}
&0\to\iota\varsigma_y(V_{2N+4k})\to V^{(-2)}_{4N}\to V_{2N-4k}\to 0
\for k\in \Z_+,\\
&0\to \iota(V^{+}_{2k}\oplus V^{-}_{2k})
\to V^{(2)}_{2N} \to V_{2N-4k}\to 0
\for k\in 1/2+\Z_+.
\end{align}
The arrows must be reversed here for $k<0$:
\begin{align}\label{eseq2}
&0\to \iota\varsigma_y(V_{2N-4|k|})\to V^{(-2)}_{4N}\to
V_{2N+4|k|}\to 0, \  k\in -1-\Z_+,\\
&0\to \iota(V_{2N-4|k|})\to V^{(2)}_{2N}\to
V^{+}_{2|k|}\oplus V^{-}_{2|k|}\to 0
,\ k\in -\frac{1}{2}-\Z_+.
\end{align}
\noindent
 If $k\not\in \Z/2$, equivalently
$t^2\not\in q^{\Z}$, then $V^{(-2)}_{4N}$ and $V^{(2)}_{2N}$
are both loc-irreducible.
\end{proposition}

Here $\varsigma_y$ can be omitted, since
it does not change the classes of isomorphism
of the modules $V_{2N\pm 4|k|}$; though the
composition map $\,\iota\varsigma_y\,$
naturally appears in the corresponding homomorphism.
\smallskip

Given arbitrary $q^{1/2}$ and $t^{1/2}$, a $\Z_{q;t}$\~free
irreducible finite-dimen\-sional module $V$ of $\HH\,$ is
called\, {\em dim-rigid\,} if all isomorphism classes of irreducible
modules of dimension equal to dim$V$ over the field of fractions
of $\Z_{q;t}$ and its arbitrary algebraic extensions $\mathbb{A}$
constitute a finite family.
Such family is obviously projective $PSL_2(\Z)$\~invariant
(where $q^{1/4}$ must be added to the definition field).
Accordingly, for $q,t$ as above, we call $V$ a
{\em loc-dim-rigit\,}
module
if there are only finitely many isomorphism classes of irreducible
free modules of dimension dim$V$ over the rings
$\mathbb{A}/\tilde{\mathfrak{p}}^{\,m}$ for prime ideals
$\tilde{\mathfrak{p}}\subset\mathbb{A}$ over $(p)$
for prime $p$ such that gcd$(p,2\!N)=1$ and any $m\in \N$.

If there is only one module in such family, then
$\tau_{\pm}$ naturally act in this module over
$\Q(q^{1/4})$ or, correspondingly, over
$\mathbb{A}/\tilde{\mathfrak{p}}^{\,m}$ (for $p,
\tilde{\mathfrak{p}}, m$
as above).
We will call
these modules\, {\em unirigid\,} or, correspondingly,
\,{\em loc-unirigid\,}. Recall that $t^{1/2}=\pm (q^{1/2})^k$
and $q^{1/2}$ is assumed primitive $(2\!N)$th root of unity.

The classification of {\em dim-rigid\,} modules is similar
to that for spherical projective $PSL_2(\Z)$\~invariant modules,
(we will begin with in the next theorem), but there are
some deviations for the family $(\breve{\ga})$ below.

\smallskip

\begin{theorem}\label{psl-two-inv}
(i) Given $N>1$, projective $PSL_2(\Z)$\~invariant
irreducible modules, i.e.  where the action of
$\,\tau_{\pm}$ can be defined, exist only if
$\,t=\pm q^{\Z/2}$. They are isomorphic to one of the
following spherical modules (irreducible
quotients of $\,\mathscr{X}$)
or their $\,\iota$\~images:
\begin{align}\label{3types}
&(\al):\ \  V_{2N-4k}\ \,  \for 2k\,\in\, \Z_+,\ \ N/2>k>0,\\
&(\be):\ \  V_{2N+4|k|} \for k\in -\Z_+,\, -N/2<k<0,
\notag\\
&(\ga):\ \  V^{+}_{2|k|}\  (k\!=\!-1/2-m>-N/2,\, m\in Z_+),\notag
\end{align}
where $t=(q^{1/2})^{2k}$ and
$\,\iota:\,k\mapsto -k, T\mapsto -T,\,\HH_{q^{1/2},t^{1/2}}\mapsto
\HH_{q^{1/2},t^{-1/2}}\,.$

(ii)
The modules from $(\al,\be)$ are $\,\varsigma$\~invariant,
namely, their isomorphism classes are invariant
under the action of $\,\varsigma_x$ and $\,\varsigma_y$.
The modules
\begin{align}\label{typegamma}
&(\breve{\ga}):\ \,
\varsigma_x^{\ep}\varsigma_y^{\de}(V^+_{2|k|})
\hbox{ for }\ep,\de=0,1,
\hbox{\, incl. } V^{+}_{2|k|},\,
V^{-}_{2|k|}\!=\!\varsigma_x(V^{+}_{2|k|}),
\end{align}
are non-isomorphic for different
$\ep,\de$;\, $\tau_{+}$ and $\tau_{-}$
transform them as in $(\ref{pslvarsigma})$.
The modules from $(\al,\ga)$
are $Y$\~semisimple and loc-semisimple,
$V_{2N+4|k|}$ from $(\be)$ is not $\,Y$\~semisimple;
these hold upon applying $\iota$.

(iii)
The modules from $(\al,\be)$, the $\,\varsigma$\~orbit of
$V^{+}_{2|k|}$ from $(\breve{\ga})$ and the $\,\iota$\~images of all
these modules constitute all {\sf \,dim-rigid\,} irreducible
modules of $\HH\,$ (for any $t$);\, they are also
{\sf \,loc-dim-rigid}.
Moreover, for a fixed primitive $(2\!N)$th
root $q^{1/2}$, the dimensions of modules in $(\al,\be,\ga)$
are all pairwise distinct. In particular, the modules of type
$(\al,\be)$ are {\sf \,unirigid\,}, which provides a conceptual
(without
explicit formulas) proof of their projective
$PSL_2(\Z)$\~invariance (especially valuable for $(\be)$).
\end{theorem}
{\em Proof.}
This theorem is essentially from \cite{C101}; for instance,
Part $(iii)$ follows from the classification of
Theorem 2.9.2 there. Given $N$ and a primitive $(2\!N)$th
root $q^{1/2}$, the dimensions of modules $(\al,\be,\ga)$
determine uniquely the corresponding $k$ (and therefore $t$).
Indeed, the dimensions are all different inside the union
$(\al)\cup(\be)$ and inside $(\ga)$; also, the dimension is
even in the first group and odd in the second.

Note that the rigidity
over $\Z_{\dot{q}}/\mathfrak{p}^m$ in Part $(iii)$ requires
knowing the dimensions of {\em all\,} free irreducible
$\HH$\~modules over this ring (and its algebraic extensions).
The classification from \cite{C101},
Theorem 2.9.2 is based on the intertwining operators, which are
generally well defined only over
$\Z_{q;t}^{loc}$. This approach can be extended to the rings
$\Z[\dot{q}]/\mathfrak{p}^m$ for $\mathfrak{p}\supset (p)$
under the condition \,gcd$(2\!N,p)=1$;\,
the latter condition is sufficient here, though not
always necessary.
\sq
\medskip

{\sf Squaring parameters for odd $N$.}
Following Sections 2.10.1-2.10.2 from \cite{C101},
let us address the case when $q^{1/2}$ is not
primitive $(2\!N)$th root of unity; it may occur only
if $N$ is odd. Instead of changing $q^{1/2}$ and
the previous considerations, we will incorporate such a case
by performing the substitution
\begin{align}\label{qtoqhat}
q\mapsto q^\checkmark=q^2,
q^{1/2}\mapsto (q^\checkmark)^{1/2}=q, t\mapsto t^\checkmark=t^2,
\end{align}
assuming that $N$ is odd.
I.e. we define the action of
$\HH^\checkmark\equal\HH_{q,t}$ in
representations $V$ of $\HH=\HH_{q^{1/2},t^{1/2}}$
considered above (in their space,
to be exact), by squaring the parameters $q,t$ in the
$e$-polynomials and the corresponding
matrices of the generators; this is closely related
to Section \ref{sec:DEFORMV} below.
In type $(\be)$, the first exact sequence from (\ref{eseq2}) is
used.
We will denote the resulting $\HH^\checkmark$\~modules by
$V^\checkmark$.


The  $\HH^\checkmark$\~modules $V^\checkmark$ for
$V$ in $(\al)$ for integral $k$ and in $(\be)$ are reducible,
namely they are direct sums
of their two $\HH^\checkmark$\~submodules of
half-dimension, which is equal to $N-2k$ and
$N+2k$ correspondingly.
The $Y$\~spectra in these two submodules
coincide and are simple for $(\al)$. In the case of $(\be)$,
the Jordan form of $Y$ is the same in both
irreducible components.

The justification of next proposition follows
case $(\ga)$ from the theorem. See also Proposition
\ref{ext-verlinde-chk} below. We use the exact sequences
from (\ref{eseq1}) and (\ref{eseq2}) and
the evaluation pairings:
\begin{align}\label{evpairing}
\{f,g\}_\pm\,\equal\, \Bigl(f\bigl((Y^\checkmark)^{-1}
\bigr)(g)\Bigr)
(X\mapsto\, \pm\, (t^\checkmark)^{1/2}
=\pm t)
\hbox { \, in \, } \mathscr{X}.
\end{align}

\begin{proposition}\label{ODDN}
For odd $N=2n+1\!>\!1$, let
$\,q, t=q^k\,$ be from Theorem
\ref{psl-two-inv}, including the inequality $0<|k|<N/2$ there.
We denote the polynomials $\ep_m$ upon the substitution $\checkmark$
by $\ep_m^\checkmark$.

(i) Upon the substitution (\ref{qtoqhat}) for {\sf integral}
 $k$ in
the cases of
$(\al,\be)$,
\begin{align}\label{oddnnew}
V_{2N-4k}^\checkmark=V^{+\checkmark}_{N-2k}\oplus
V^{-\checkmark}_{N-2k},\
V_{2N+4k}^\checkmark=V^{+\checkmark}_{N+2k}
\oplus V^{-\checkmark}_{N+2k},
\end{align}
where the $\HH^\checkmark$\~modules $V^{\pm\checkmark}_{N-2k}$,
$V^{\pm\checkmark}_{N+2k}$  are the quotients of
$V_{2N-4k}^\checkmark$ and
$V_{2N+4k}^\checkmark$ by the radicals of the corresponding
evaluation pairings
$\{f,g\}_{\pm}$. One has
$$
\ V^{\pm\checkmark}_{N-2k}\,=\,
\oplus_{m=1}^{N-2k}\,\Z^{loc}_{q;t}\,(\ep^\checkmark_m\pm
\ep^\checkmark_{m-N+2k})=\mathscr{X}/(\ep^\checkmark_{n-k+1}\mp
\ep^\checkmark_{k-n}).
$$

(ii) Continuing (i),
the modules $V^{\pm\checkmark}_{N-2k}$,
$V^{\pm\checkmark}_{N+2k}$ are $\tau_{\pm}$\~invariant. Their images
under the action of $\varsigma_x^\ep\varsigma_y^\de$ for
$\ep,\de=0,1$  are transformed
by $\tau_{\pm}$ as for $(\breve{\ga})$ from the
theorem and are pairwise non-isomorphic.
The corresponding families
\begin{align}\label{brevealbe}
(\breve{\al}^\checkmark):\,
 V^{\ep,\de,\checkmark}_{N-2k}\!=\!
\{\varsigma_x^\ep\varsigma_y^\de(V_{N-2k}^{+\checkmark})\},
\   (\breve{\be}^\checkmark):\,
 V^{\ep,\de,\checkmark}_{N+2k}\!=\!
\{\varsigma_x^\ep\varsigma_y^\de(V_{N+2k}^{+\checkmark})\}
\end{align}
constitute all $\HH^\checkmark$\~modules of dimensions $N-2k$
and $N+2k$.

(iii) The modules $V^\checkmark_{2N-4k}$ for $k=1/2+n$ and
$V_{2|k|}^{\ep,\de,\checkmark}$ for $k=-1/2-n$ in
$(\breve{\ga}^\checkmark)$ are (remain) irreducible
for $n\in \Z_+$ upon the substitution (\ref{qtoqhat}).
The modules $V_{2N-4k}^\checkmark$ here remain
$\tau_{\pm}$\~invariant and become direct sums of
Jordan $2$\~blocks under the action of $Y$ or $X$.
One has
$\ V^{\pm\checkmark}_{2|k|}\,=\,
\oplus_{m=1}^{2|k|}\,\Z^{loc}_{q;t}\,(\ep^\checkmark_m-
\ep^\checkmark_{m-2|k|})$.
Thus
$(\breve{\ga}^\checkmark)\ :
\  V^{\ep,\de,\checkmark}_{2|k|}=
\{\varsigma_x^\ep\varsigma_y^\de\bigl(V_{|k|}^{+\checkmark}\bigr)\}$
consists of all pairwise non-isomorphic $\HH^\checkmark$\~modules
of dimension $2n+1$. The action of $\tau_{\pm}$
in this family remains the same as that for $(\breve{\ga})$.
\end{proposition}
{\em Proof.} Let us focus here on the series
 $(\al)$ for integral $k$.
The key is the following fact.
The $e$\~polynomials
(when exist) and the inner products $\lan\,,\,\ran_{\pm}$
depend only on $\,q,t$,\, though
$T,Y$ and the evaluation pairings
from (\ref{evpairing}) do involve $t^{1/2}$. Thus
changing $t^{1/2}$ to $-t^{1/2}$, which sends $\,T\mapsto -T$ and
$Y\mapsto -Y$,\, does not influence the $e$\~polynomials.

Generally, if the quotient of $\mathscr{X}$ by
the radical of $\lan\,,\,\ran$ is already irreducible,
then using the evaluation pairing on top of it will
produce no further reduction. However upon changing $\,q,t\,$ by
$\,q^2,t^2$,\, the spectrum of $Y^\checkmark$ in
$V_{2N-4k}^\checkmark$ becomes of multiplicity $2$,
the radicals of $\{\,,\,\}_{\pm}$
in $V_{2N-4k}^\checkmark$ become nonzero and they produce the
split from (\ref{oddnnew}). Note that
for odd $N$ and integral $k$, the transform $\,q,t\mapsto q^2,t^2\,$
is a {\em Galois automorphism\,} sending the coefficients
of $e_m$ to those of $e_m^\checkmark$;
so these polynomials exist simultaneously.

One can also use here the exact sequence (\ref{eseq1}).
The $\HH^\checkmark$\~module
$V^{(-2)}_{4N}=\mathscr{X}/(X^{2N}+X^{-2N}-2)$ becomes
the following direct sum of its submodules:
$V^{(-2)}_{4N}=$
$\mathscr{X}/(X^{N}+X^{-N}-2)\oplus
\mathscr{X}/(X^{N}+X^{-N}+2)$. Note that these submodules
are transposed by $\varsigma_x$, which sends $X\mapsto -X$.
Accordingly, (\ref{eseq1})
will be reduced for $k\in \Z_+$ to
\begin{align}\label{eseq1check}
&0\to\iota\varsigma_y(V_{N+2k}^{\pm\checkmark})
\to \mathscr{X}/(X^{N}+X^{-N}\mp 2)\to
V_{N-2k}^{\pm\checkmark}\to 0.
\end{align}

Concerning $(iii)$, let us briefly consider the series
$(\al)$ for half-integral $k$. We determine when
the $e$\~polynomials from $V_{2N-4k}$ become non-existent
upon applying $\checkmark$, which results in the loss of
$Y$\~semisimplicity in this case. This is essentially
sufficient for the remaining claims of the proposition
for $(\al^\checkmark)$. The case of $(\breve{\ga}^\checkmark)$
is straightforward.
\sq

\subsection{\bf Deformations of Verlinde algebras}
\label{sec:DEFORMV}
The subalgebra $V_{2N-4}^{sym}$ of $V_{2n-4k}$ for $k=1$
of symmetric elements $v\in V_{2N-4}$
is isomorphic to the {\em usual Verlinde algebra\,} for
$\widehat{\mathfrak{sl}_2}$ of level $N$
(with the central charge $N-2$). Its dimension is
$N-1$; generally, the
dimensions of the symmetric subalgebras,
defined as $\{v\in V\,:\, Tv=t^{1/2}v\}\,$, of the
modules under consideration are as follows:
\begin{align}\label{symdim}
&\hbox{dim\,} V_{2N-4k}^{sym}=N-2k+1,\
\hbox{dim\,} V_{2N+4|k|}^{sym}=N+2|k|+1,\notag\\
&\hbox{dim\,} (V^{\pm}_{2|k|})^{sym}=\,m+1\, \for k=-1/2-m,\
m\in \Z_+.
\end{align}
It matches the exact sequences from (\ref{eseq1})
and (\ref{eseq2}), since the symmetric parts of the middle
modules $V^{(-2)}_{4N}$, $V^{(2)}_{2N}$ are
of dimension $2N$ and $N$ correspondingly;\,
$\,\iota\,$ transposes
the eigenvalues $\,t^{1/2}$ and $\,-t^{-1/2}$ of $T$.

Here and further by {\em Verlinde algebras\,},
we mean {\em generalized\,}
(any admissible $k$) and {\em nonsymmetric\,} ones (without
taking the symmetrization).
We note that the projective action of $PSL_2(\Z)$ in
$V_{2N-4k}^{sym}$ for integral $0<k<N/2$
was defined in \cite{Ki} using the
interpretation of Macdonald polynomials via
quantum $GL$ (Etingof, Kirillov).
\smallskip

Now let us address the absence of  modules
$V'_{4|k'|}$ for $-2k'\in \N$, defined in $(iii)$
of Theorem \ref{tnegk}, in the classification
of projective $PSL_2(\Z$)\~invariant and {\em dim-rigid\,}
modules from Theorem \ref{psl-two-inv}.

The next proposition
connects the series $(\al)$
with the modules  $V'_{4|k'|}$ for $2k'\in -\N$ defined in $(iii)$
of Theorem \ref{tnegk}. Cf. Section 2.10.5, the remark after
Theorem 2.9.9 from \cite{C101} and (\ref{litVerlinde}) below.
 This explains why $(\al,\be,\ga)$
are sufficient in Theorem \ref{tnegk}, and we do not need
to add there the specializations
of  modules $V'_{4|k'|}$ at roots of unity.

\begin{proposition}\label{ext-verlinde}
Continuing to assume that $\,q^{1/2}$ is a primitive $(2\!N)$th root
of unity and considering the modules $V_{2N-4k}$ from $(\al)$,
let $k'=k-N/2$.
Then
\begin{align}\label{add-roots}
&t=(q^{1/2})^{2k}=(q^{-1/2})^{N+2k'}=-(q^{1/2})^{2k'}
\and\notag\\
&V_{2N-4k}\,\cong\, V'_{4|k'|}, \hbox{\ setting\ \,} k'=k-N/2,
\end{align}
where the latter module is well defined for primitive
$q^{1/2}$ of order $2\!N$.
Since $V'_{4|k'|}$ are well defined for any $q\neg 0$
and are irreducible for sufficiently general $q$,
they are flat
$q$\~deformations of modules $V_{2N-4k}$ with all
their structures. Namely, the dimension remains the same
and the action of projective
$PSL_2(\Z)$ and the multiplication in $V'_{4|k'|}$
(it is $Y$\~spherical) deform those in
$V_{2N-4k}$. Then, $V'_{4|k'|}$ is
$X$\~ and $Y$\~semisimple for generic $q$ and the standard
Hermitian inner product in $V_{2N-4k}$ is the specialization
of that in $V'_{4|k'|}$, which is direct from
(\ref{fnorme}). Moreover, the inner product in $V'_{4|k'|}$
remain positive definite if $|q|=1$ provided\,
$|\hbox{arg}(q)|\le 2\pi/N$. \sq
\end{proposition}
\smallskip

For instance, the $q$\~deformation
of the standard $\widehat{\mathfrak{sl}}_2$\~Verlinde algebra
$V_{2N-4}^{sym}$ of (for $k=1$  and of dimension $N-1$)
is $(V'_{4|k'|})^{sym}$ for $k'=1-N/2$. We note that
Verlinde algebras of type $(\al)$ and their deformations
can be defined for arbitrary (reduced) root systems.
\smallskip

The deformations of nonsymmetric and symmetric
generalized Verlinde algebras (of type $(\al)$) were
constructed in \cite{C101} only for odd $N$
and integral $0<k<N/2$. Also, it was done there upon the
restriction to
the nonsymmetric {\em Little Verlinde algebra\,}, defined as
the image $V_{N-2k}^{even}$ of even polynomials
from $\mathscr{X}$ in  $V_{2N-4k}$. It is a module over
{\em Little DAHA\,} $\HH^{ev}$,
defined as the span of
$X^{\pm2},Y^{\pm2},T\,$ in $\HH\,.$


The polynomials $e_m$ are of the same parity as $m$, so
the subspace (and a subalgebra) $V_{N-2k}^{even}$ is
linearly generated
by $e'_m$ over $\Z_{\dot{q}}^{(N)}$ for even $m$ taken from the set 
$\{2k-N+1,\ldots,N-2k\}$
(note that the endpoints here have different parity).
It is irreducible for odd $N$ and has two irreducible components
of dimension $N/2-k$ for even $N$.
\smallskip

Accordingly, $V_{N-2k}^{odd}$  is the $\HH^{ev}$\~module
generated by the images of odd polynomials from $\mathscr{X}$;
it has the same dimension $N-2k$.

The claim from \cite{C101} is as follows:
\begin{align}\label{litVerlinde}
&V_{N-2k}^{even}\simeq V_{2|k'|}\,, \hbox{\,\, where \,\,}
k'=k-N/2 \hbox{\,\, for odd\,\, } N,
\end{align}
which is upon the restriction to $\HH^{ev}$ and when
$q^{1/2}=-e^{\nu\pi \imath/N}$ for \,gcd$(\nu,2N)=1$.
 Therefore   $V_{2|k'|}$
for generic $q$ are flat deformations of  $V_{N-2k}^{even}$
(with all their structures). Note that the choice of the sign
of $q^{1/2}$
influences $t'=q^{k'}$ ($k'$ is a half-integer) and
makes the inner product in $V_{2|k'|}$
coinciding with that in $V_{N-2k}^{even}$.
The inner product in the latter module
involves only $q$; see (\ref{fnorme}). Picking
$q=e^{2\pi \imath/N}$, the ``smallest" primitive $N$th root of
unity, makes the inner product in $V_{N-2k}$
 positive definite for any $2k\in \N$ here and above.

The following proposition is related to
(\ref{litVerlinde}), however is not its reformulation.
\begin{proposition}\label{ext-verlinde-chk}
For odd \,$N$ and \,$\Z_+\!\ni k\!< N/2$, let
$q=e^{2\nu\pi \imath/N}$, where \,gcd$(\nu,N)=1$.
Then the module \,$V^{+\checkmark}_{N-2k}$\! of type
 $(\al^\checkmark)$ from Proposition \ref{ODDN} is
isomorphic to
\!$V^{+\checkmark}_{2|k'|}$ for $k'=k-N/2$.\sq
\end{proposition}

Note the operator $X^\checkmark$ {\em concides\,}
with $X^2$ in $V_{2N-4k}$ for the basis of
$X$\~eigenfunctions. However the corresponding
$Y^\checkmark$ (in the same basis) is not connected with $Y$
in any direct way.


\subsection{\bf Boundary cases}\label{sec:Boundary}
With minor reservations,
the statements of Theorem \ref{psl-two-inv} can be extended to the 
boundary cases, namely when $t=\pm 1$ and, respectively,
$k=0,-N/2$.
The list of irreducible $\tau_{\pm}$\~invariant
quotients of $\mathscr{X}$ from $(ii)$ then becomes:
\begin{align}\label{boundcasea}
&(\al^\bullet)\ \,  t=1, k=0\,: \ V_{2N}=\mathscr{X}/(X^N-X^{-N}),\\
(\be^\bullet)\ \, t=-1, k&=-\frac{N}{2},\ \hbox{even\,} N\,: \
V_{4N}=\mathscr{X}/(X^{2N}+X^{-2N}-2),\notag\\
(\ga^\bullet)\ \, t=-1, k&=-\frac{N}{2}, \ \hbox{odd\ } N\,:\
V^+_{N}=\mathscr{X}/((\sqrt{-1}X)^{N}-1).\notag
\end{align}

When $k=0$ and $t=1$, the
$\HH$\~module $V_{2N}$ is irreducible
$\tau_{\pm}$\~invariant and {\em dim-rigid\,}.
It is not {\em unirigit\,},
but is unique of dimension $2N$
up to applying $\iota:T\mapsto -T$.
We note that applying $\iota$  does not change the
quadratic equation for $T$ when $t=1$, however
$\iota(V_{2N})$ is not a $Y$\~spherical module.

In this case, we can set $T=s$, where  $s^2=1,\ sXs=X^{-1}$ and
$\,sYs=Y^{-1}$. The restriction of $V_{2N}$ to $X,Y$ is
the classical (unique)
finite-dimensional module of the Weyl subalgebra $\w\subset \HH\,$ 
which is generated by $X^{\pm1},Y^{\pm1}$ subject to the relation
$Y^{-1}X^{-1}YXq^{1/2}=1$.
The action of $s$ in this module is determined uniquely
up to the multiplication by $\pm 1$, which is
exactly applying $\iota$.
The isomorphism class of this
$\HH$\~module is invariant under $\varsigma_x$ and $\varsigma_y$;
the latter automorphism corresponds to the conjugation by the
polynomial $X^N+X^{-N}$.
\smallskip

Now let $t=-1$. Then the eigenvalues of $T$ coincide and
$T$ becomes non-semisimple in $\mathscr{X}$
and its quotients $V$ from  $(\be^\bullet,\ga^\bullet)$
unless $N=1$. Indeed, the semisimplicity of $T$ would
result in $T=t^{1/2}$ in such $V$ (but it is not a scalar).

For even $N$, the module $V_{4N}=V^{-2}$ is
irreducible; it is a direct sum of Jordan $2$\~blocks
for $X$, $Y$. This module is $\tau_{\pm}$\~invariant
and invariant with respect to
$\varsigma_x$ and $\varsigma_y$. The invariance with respect
to $\varsigma_y$  can be
seen by extending the formulas for $t\neq -1$ to
this special case. Here and for the description
of $V^+_{N}$ for odd $N$ to be discussed next,
see $(iii,iv)$ in Theorem 2.9.3 from \cite{C101}.

Let us (re)establish the projective $PSL_2(\Z)$\~invariance
of $V_{4N}$ using the rigidity approach.
The problem is that the dimension $4N$ of $V_{4N}$ is maximal among
irreducible representations for such $q,t$, so it is not
{\em dim-rigid\,}. However the Jordan
decomposition of $T$ in
$V_{4N}$\, is unique among all modules $V^{C}=
\mathscr{X}/(X^{2N}+X^{-2N}+C)$. One has
$$
T(P\De_-)\,=\,\imath P\De_-\, -2\imath P\De_+
\for \De_{\pm}=X\pm X^{-1}
\hbox{ and any symmetric } P.
$$
Let $P=(X^{2N}-X^{-2N})^2/(X+X^{-1})$, which is a
symmetric Laurent polynomials since $N$ is even. Then
\begin{align}\label{TPDe}
T(P\De_-)\,=\, \imath P\De_-\, -2\imath(X^{2N}-X^{-2N})^2\, =\,
\imath P\De_- \hbox{ \, in\, } V^{-2}.
\end{align}
For  $C\neq -2$, there will be an extra term $-(C+2)$ in the
right-hand side.
Also, there can be no nonsymmetric
polynomials of degree smaller than $P$ satisfying $T(P)=\imath P$
in any $V^C$.

Thus $V^{-2}$
is a direct sum of $(2N-1)$ Jordan $2$-blocks with respect
to $T$ and
$2$ one-dimensional $T$\~eigenspaces. The modules $V^{C}$
for $C\neq -2$ are direct sums of $2N$ Jordan
$2$-blocks under the action of $T$.
Therefore $V_{4N}=V^{-2}$ for $t$ from $(\be^\bullet)$
is, indeed, projective $PSL_2(\Z)$\~invariant.
\smallskip

The $\bullet$-module $V_{N}=V^{+}_{N}$ for odd $N$ is
{\em dim-rigid\,} and
$Y$\~semisimple; its $Y$\~spectrum
is simple. It is $\tau_{\pm}$\~invariant, but not invariant
with respect to $\varsigma_x$ and $\varsigma_y$.
Let $(\breve{\ga}^\bullet)$ denote the natural boundary
of the family $(\breve{\ga})$ from (\ref{typegamma})
for such $t=-1,k=N/2$.

The simplest (but instructional)
example of odd $N$ is the module
$V^{+}_1$, i.e. when $N=1$ and $k=-1/2$\,:
\begin{align*}
&t=q^{-1/2}=-1 \hbox { and } \ t^{1/2}=q^{-1/4}=\imath=\sqrt{-1},\\
&V_1=\mathscr{X}/(\imath X-1),\hbox{ where } T=Y=\imath=X^{-1}.
\end{align*}
Note that $T(\imath X-1)=X^{-1}-\imath,\ \pi(\imath X-1)=
-\imath X^{-1}-1$, so
the ideal $(\imath X-1)$ is obviously $\HH$\~invariant
in this case.

The existence of the action
of $\tau_+$ in $V^+_1$ is equivalent to the relation
$\tau_+(Y)=Y$,
which obviously holds since
$\tau_+(Y)\equal q^{-1/4}XY= \imath XY=\imath (-\imath)
\imath=\imath=Y$ in $V^+_1$. Similarly,
$\tau_-$ acts in $V_1^+$ too and
$\tau_+\tau_-=id\,$ in this module.
Note that  $\tau_+(Y)=-Y$
in $\varsigma_x(V_1^+)$ and therefore $\tau_+$ does not
act in $\varsigma_x(V_1^+)$ (it is not a conjugation by a matrix).


\subsection{\bf Finite and discrete images}
Let us define the elliptic $q$\~braid group as follows:
\begin{align}\label{braidg}
&\b_q\equal \lan X,T,Y,q^{1/2}\ran
\hbox{\, for a primitive $(2\!N)$th root\, } q,\\
\hbox { subject to }
&TXT=X^{-1},\ T^{-1}YT^{-1}=Y^{-1},\ Y^{-1}X^{-1}YXT^2q^{1/2}=1.
\notag
\end{align}
We will also denote by $\b_q^{ev}$ its subgroups generated
by $X^2, T,Y^2$ associated with $\HH^{ev}$. Recall that the
latter algebra acts in $V_{N-2k}^{even}$ and
$V_{N-2k}^{odd}$ (of dimension $N-2k$ over $\Z_{\dot{q}}$) generated
by the images of even and odd polynomials in $V_{2N-4k}$
for the modules from $(\al)$ in (\ref{3types}).

The {\em elliptic braid group\,}
$\b_1$ is obtained by omitting $q^{1/2}$ in $\b_q$.
It is the following renormalization of
$\b_q$ (which will be used later):

\begin{align}\label{braidnorm}
&\hbox{setting\ \,} \tilde{T}\,=\,q^{+1/4}T,\
\tilde{Y}\,=\,q^{+1/4}Y,\
\tilde{X}\,=\,q^{-1/4}X,\\
&\tilde{T}\tilde{X}\tilde{T}=\tilde{X}^{-1},\
\tilde{T}^{-1}\tilde{Y}\tilde{T}^{-1}=\tilde{Y}^{-1},\
\tilde{Y}^{-1}\tilde{X}^{-1}\tilde{Y}\tilde{X}\tilde{T}^2=1.
\notag
\end{align}

We will describe in this section all cases of {\em dim-rigid\,}
modules such that the image $\,Image(\b_q)$ of $\b_q$ in the
corresponding
matrix groups is finite. Without loss of generality,
we can restrict ourselves with
$(\al,\be,\ga)$ from Theorem \ref{psl-two-inv}
and $(\al^\checkmark,\be^\checkmark,\ga^\checkmark)$
from Proposition \ref{ODDN}; indeed, applying
$\iota,\varsigma_x,\varsigma_y$ does not change
the finiteness.

We analyze the positivity of the
inner product $\lan \cdot\,,\,\cdot \ran$.
The same approach provides examples when $\,Image(\b_q)$ is
 {\em discrete\,} to be discussed next
(though not a classification of all such cases).

\begin{theorem}\label{FINBIMAGE}
We omit one-dimensional modules, those from
$(\al^\bullet)$ in Section \ref{sec:Boundary}
and $(\ga^\checkmark)$ for $t^\checkmark=1$
from Proposition \ref{ODDN}; in these
cases the images are finite.
Then among the modules $(\al,\be,\ga)$
from Theorem \ref{psl-two-inv} and those from
Proposition \ref{ODDN}, the group $Image(\b_q)$ can be
finite only for $(\al)$ 
and  $(\al^\checkmark)$ for the following $k$.

There are two series of such
modules $V=V_{2N-4k}$ of type $(\al)$
correspondingly for $N>1$ and for $N>2$\,:
\begin{align}
&(\al^{-1}_N):\, k=(N-1)/2, \hbox{dim}V=2,\ \,
(\al^{-2}_N):\, k=(N-2)/2,
\hbox{dim}V=4.\notag
\end{align}
Apart from these series, there are only two exceptional cases:
\begin{align}
&(\al_{6}^1):\,  N=6,\, k=1,\, \hbox{dim}V=8,\ \ \
(\al_{10}^3):\,  N=10,\, k=3,\, \hbox{dim}V=8.\notag
\end{align}

Also, upon $q\mapsto q^\checkmark=q^2, t\mapsto
t^\checkmark=t^2$ as in Proposition \ref{ODDN},
there is the following additional series
of such modules $V_{N-2k}^{+\checkmark}$ for even $N\!>\!2$:
\begin{align}
&(\al^{-2\checkmark}_N):\, k=(N-2)/2,\,
\hbox{dim}V^\checkmark=2.\notag
\end{align}
\end{theorem}
{\em Proof.}
The operators $T,X,Y$ are
represented by matrices with the entries in
$\Z[q^{1/4}]$ (over $\Q(q^{1/4})$ in the basis of
$e$\~polynomials) for
the modules from $(\al,\ga)$.
Actually $\Z[q^{1/2}]$ or $\Q(q^{1/2})$
is sufficient here if we switch
from $T$ to $t^{-1/2}T$ and $Y$ to $t^{-1/2}Y$, which
holds for $(\be)$ as well. Accordingly, the field
 $\Q(q)$ is sufficient for modules $V^\checkmark$.

Now assume that $Image(\b_q)$ in GL$(V)$ is finite for all
$(2\!N)$th primitive root of unity $q^{1/2}$ and the corresponding
$t^{1/2}$.
Then for every such $q$ (embedded into $\C$), there exists a positive
definite hermitian inner product $(\cdot\,, \,\cdot)$ in $V$
such that
$(Hu,v)=(u,H^{-1}v)$ for $H=X,T,Y$ and any vectors $u,v\in V$. This
results in the semisimplicity of $X,Y$ in $V$; therefore
the cases
$(\be),(\be^\checkmark)$ and $(\al^\checkmark)$ for
odd $2k$ can be excluded from the consideration.
Moreover, such an inner product must be proportional to the
standard one $\lan \cdot,\,,\, \cdot\ran$ in the cases
$(\al,\ga)$ due to the irreducibility of $V$.

We conclude that the finiteness of $Image(\b_q)$ may occur
only for $(\al,\ga)$ and for $(\al^\checkmark)$
when $k$ is integral. This is equivalent
to the positivity of $\lan \cdot\,,\, \cdot\ran$ for any primitive
$q^{1/2}$, i.e. for any embedding
$\Z[q^{1/2}]\hookrightarrow \C$. It is straightforward to show
that only $(\al_N^{-1}), (\al_N^{-2}), (\al_6^2),
(\al_{10}^3), (\al^{-2\checkmark}_N)$ posses such
positivity, though we use computers to find that there
only $2$ exceptional cases  (apart from the infinite series).
\sq

\smallskip

\rmk
One can examine the total
positivity of the form $\lan \cdot\,,\,\cdot\ran$
upon the restriction to $V^{sym}$. The only new case
of total positivity (for all $\nu$) is the
{\em Verlinde series} $k=1$; in this case,
$\lan p_n\,,\,p_n\ran=1$ for any Rogers-Macdonald polynomials $p_n$.
Furthermore, there are no new cases of total positivity upon the
restriction to the even part
$V^{sym+even}$ of $V^{sym}$ from $(\al,\ga)$.
However  we have the following additional cases of total
positivity for  $\tilde{V}=V_{N-2k}^{sym+odd}=
V_{N-2k}^{sym}\cap V_{N-2k}^{odd}$\,:
\begin{align*}
&k=(N-4)/2,\, \hbox{dim}\tilde{V}=4\hbox{ for all } N> 4,\ \,
N=12,\, k=2,\, \hbox{dim}\tilde{V}=8.
\end{align*}

The groups $\b_q$, $\b_q^{ev}$ do not preserve
 $V_{N-2k}^{sym}$, though the projective
$PSL_2(\Z)$ and the absolute Galois group (see below)
act there. Note that
$\b_q^{ev}$ preserves $V^{odd}$ and $V^{even}$,
but the (total) positivity property remains the same in
these modules
as for $V$. Indeed,  $\lan e_n,e_n\ran=\lan e_{1-n},e_{1-n}\ran$
for any $n$ (provided the existence of $e_n$).
\sq

\begin{proposition}\label{bqcenter}
The center of the group $Image(\b_q)$ for $V$ of types
$(\al,\be,\ga)$
belongs to the group
generated by the scalar matrix $\dot{q}$, where
$\dot{q}=q^{1/2}$ for $2k\in \Z$
and $\dot{q}=q^{1/4}$ for $2k\in 1/2+\Z$.
\end{proposition}
{\em Proof.} Due to the irreducibility of $V$, the center
elements are scalars $z\in \Z[\dot{q}]$.
Following the proof of Theorem
\ref{FINBIMAGE}, $z$ is unimodular for types $(\al,\ga)$ for any
embeddings $\Z[\dot{q}]\hookrightarrow \C$.
The corresponding hermitian forms can be non-positive, but
this is not a problem since we consider only central elements.
Thus $z$ is
a root of unity, which can be only a power of $\dot{q}$.
Modules of type $(\be)$ posses hermitian invariant
(non-positive) forms too, the leading terms of those in the
polynomial representations covering $V$.
\sq


\smallskip
{\sf Discrete  braid-images.} 
It is of interest to apply the method used in the proof
of the theorem, to construct examples of {\em discrete\,}
groups $Image(\b_q)$.  We can use that upon the
switch from $T$ to $t^{-1/2}T$ and $Y$ to $t^{-1/2}Y$,
the images of $T,Y,X$ can be made with entries in
$\Z[q^{1/2}]$; in the basis of $e_m$ for
$m=2k-N+1,\ldots,N-2k$, they are matrices over $\Q(q^{1/2})$,
which field is sufficient here.

The discreteness of non-finite $\,Image(\b_q)$ holds if there
is exactly
one complex place (valuation) of $\Q(q^{1/2})$ where the form
$\lan \cdot\,,\,\cdot\ran$ is non-positive; in this case
the groups $Image(\b_q)$ are {\em arithmetic discrete\,}.
Respectively, $\Q(q)$ must be considered for the
$\checkmark$-series. Omitting $(\be)$,  $(\be^\checkmark)$
and the series $(\al^\checkmark)$ for odd $2k$
from Proposition \ref{ODDN},
our considerations show
that there are only finitely many such cases.


Let us list what was found. We set
$q^{1/2}=\exp(\pi i \nu/N)$ when $1\le \nu \le
2\!N, \  (\nu,2\!N)=1$.
Only the modules of type $(\al,\ga)$ and those
of type $(\al^\checkmark,\be^\checkmark)$
from Proposition \ref{ODDN}
will be considered
below. The results of our analysis and
computer calculations are as follows.
\smallskip

There exists exactly one pair
$q^{\pm 1/2}=e^{\pm \frac{\pi i \nu}{N}}$
of primitive $(2\!N)$th
roots of unity such that the inner product
$\lan\cdot\,,\,\cdot\ran$
is non-positive only (presumably)
for the following modules of
type $(\al)$:
\begin{align}\label{aldiscr}
&\hbox{\small $\{N=4,\,k=1/2,\hbox{dim}V=6, \nu=3\}$},\
&&\hbox{\small $\{N=5,\,k=1/2,\hbox{dim}V=8, \nu=3\}$},\notag\\
&\hbox{\small $\{N=5,\,k=1,\hbox{dim}V=6, \nu=3\}$},\
&&\hbox{\small $\{N=6,\,k=1/2,\hbox{dim}V=10, \nu=5\}$},\notag\\
&\hbox{\small $\{N=6,\,k=3/2,\hbox{dim}V=6,  \nu=5\}$},\
&&\hbox{\small $\{N=7,\,k=3/2,\hbox{dim}V=8, \nu=3\}$},\notag\\
&\hbox{\small $\{N=7,\,k=2,\hbox{dim}V=6, \nu=5\}$},\
&&\hbox{\small $\{N=9,\,k=2,\hbox{dim}V=10, \nu=7\}$},\notag\\
&\hbox{\small $\{N=9,\,k=5/2,\hbox{dim}V=8, \nu=5\}$},\
&&\hbox{\small $\{N=15,k=4,\hbox{dim}V=14, \nu=13\}$},\notag\\
& &&\hbox{\small  $\{N=15,k=\frac{11}{2},\hbox{dim}V=8, \nu=7\}$}.
\end{align}

Let us now discuss the type $(\al^\checkmark)$ for integral $k$
and $(\ga^\checkmark)$.
In type $(\al^\checkmark)$, we have the following
$4$ additional cases (all with $\nu=1$):
\begin{align}
&\{N=5,k=1,\, N=7,k=2\, \hbox{dim}V=3\},\label{aldisodd}\\
\{N=9,\, &k=2,\, \hbox{dim}V=5\},\,\
\{N=15,\, k=4,\, \hbox{dim}V=7\}.\notag
\end{align}

\comment{
\begin{align}\label{aldisodd}
&\hbox{\small $\{N=5,\,k=1/2,\hbox{dim}V=4, \nu=2\}$},\ \ \
&&\hbox{\small $\{N=5,\,k=1,\hbox{dim}V=3, \nu=1\}$},\notag\\
&\hbox{\small $\{N=7,\,k=2,\hbox{dim}V=3,  \nu=1\}$},\ \ \
&&\hbox{\small $\{N=9,\,k=2,\hbox{dim}V=5, \nu=1\}$},\notag\\
& &&\hbox{\small $\{N=15,k=4,\hbox{dim}V=7, \nu=1\}$}.
\end{align}
}

These additional cases exactly correspond to
the ones in (\ref{aldiscr}) for odd $N$ and integral $k$,
but they are of half-dimension. They also
{\em coincide\,} with those obtained in (\ref{gadischk})
below from $(\ga^\checkmark)$. The coincidence is due to the
isomorphism from Proposition \ref{ext-verlinde-chk}.

\smallskip

In type $(\ga)$, the inner product is always indefinite
for $\nu=\pm 1$; it is  positive definite at all other places
in the following cases:
\begin{align}
&\{N=3,4,5,6,9,\, k=-3/2,\, \hbox{dim}V=3\},\label{gadiscr}\\
&\{N=6,\, k=-5/2,\, \hbox{dim}V=5\}.\notag
\end{align}

For $(\ga^\checkmark)$, we also have the following
$4$ additional cases:
\begin{align}
&\{N=5,7,\, k=-3/2,\, \hbox{dim}V=3\},\label{gadischk}\\
\{N=9,\, k=-5/2,\, &\hbox{dim}V=5\},\,\,
\{N=15,\, k=-7/2,\, \hbox{dim}V=7\}.\notag
\end{align}
This list is actually that in
$(\ref{aldisodd})$.

Summarizing,
the groups $Image(\b_q)$
from (\ref{aldiscr}-\ref{gadiscr}), are discrete
and non-finite in GL${}_{\C}(V)$ with respect to the embedding
$\Q(q^{1/2})\hookrightarrow \C$
such that  $q^{1/2}\mapsto e^{\frac{\pi i \nu}{N}}$, where
$\nu$ is from these lists.

The sequences (\ref{gadiscr}), (\ref{gadischk}) for dim$V=3$
are connected in Proposition \ref{LIVNE} below with the arithmetic
Livn\'e discrete groups. Moreover, the theory of the latter
provides an interesting example, when $\b_q$ is discrete
but {\em not\,} arithmetic. It occurs for dim$V$=3, $N=9$;
this case is missing in (\ref{gadischk}) because there are
{\em two\,} places (and four $\nu$) such that the inner product
$\lan\,,\,\ran$ is non-positive. Nevertheless the corresponding
group is discrete! We mention that the Livn\'e groups are
examples of {\em Mostow's groups\,} \cite{Mos}. It is interesting
to compare the cases with dim$V>3$ listed
above with the Mostow list.

\begin{conjecture}\label{DISCCONJ}
Up to isomorphisms and changing the signs
of the images of \, $T,Y,X$, \, arithmetic discrete
non-finite groups
$\,Image(\b_q)$ for irreducible $\HH\,$\~modules
and irreducible $\HH^\checkmark$\~modules with semisimple
$X$ and $Y$ can be only those from the lists
(\ref{aldiscr}-\ref{gadischk}).
\end{conjecture}
According to R.~Schwartz, discrete images cannot appear
if $T$ is elliptic and $X$ or $Y$ is parabolic when dim$V=3$;
see Theorem 1.1 from \cite{Par1}. Thus the semisimplicity
of $X,Y$ does not seem really necessary.

\subsection{\bf Series of dim=\texorpdfstring{$2,4$}{2,4}}
\label{sec:dim24}
Let us give a complete description of the series
$(\al_N^{-1}): \ k=(N-1)/2,\ \hbox{dim}V=2$. Using
Section \ref{sec:DEFORMV}, it suffices to consider
$V'_2$ defined for
$t=-q^{-1/2}$ for any $q$; then we will make
$t=q^{\frac{N-1}{2}}$ for primitive $(2\!N)$th
roots of unity $q^{1/2}$.

One has
$V'_2=\mathscr{X}/(e_{-1})=\Z_{q;t}\,e'_0\oplus\Z_{q;t}\,e'_1$,
where here and below $e'_i$ is the image of $e_i$ in $V'$.
Explicitly,
\begin{align*}
&e_0=1,\ e_1=X, \ e_{-1}\,= X^{-1}+\frac{1-t}{1-tq}X\,=X^{-1}-tX,\\
&Y(e_0)=t^{1/2},\, Y(e_1)=-t^{1/2}e_1,\ T(e'_i)=t^{1/2}e'_i\,(i=0,1)
\end{align*}
Thus this representation is as follows:
\begin{align}\label{dim2fin}
X=\begin{pmatrix} 0 & t^{-1}  \\ 1 & 0 \\ \end{pmatrix},\
Y=\begin{pmatrix} t^{1/2} & 0  \\ 0 & -t^{1/2} \\ \end{pmatrix},\
T=\begin{pmatrix} t^{1/2} & 0  \\ 0 & t^{1/2} \\ \end{pmatrix}.
\end{align}
The relations $Y^{-1}X^{-1}YX=-1$,\, $Y^2=X^{-2}=t=-q^{-1/2}$\,
and $T=t^{1/2}$
determine $\,Image (\b_q)$ up to isomorphisms. Equivalently,
$$\tilde{Y}^{-1}\tilde{X}^{-1}\tilde{Y}\tilde{X}=-1=  \tilde{T}^2=
\tilde{Y}^2=\tilde{X}^{2},\hbox{\, where\, } \tilde{T}=\pm \imath,
$$
in the normalization of (\ref{braidnorm}) from $\b_q$ to $\b_1$.

If $q^{1/2}$ is a primitive $(2\!N)$th
root of unity, the group  $\,Image (\b_q)$
is finite of order $16\!N$ unless $(t^{1/2})^N=1$
for odd $N$ and the order is $8\!N$. Accordingly,
$\,Image (\b_1)$ divided by its center, generated by $\imath$,
is $\, \Z_2^2$.
\medskip

{\sf Dimension $4$}.
The series $(\al_N^{-2})$ of dimension $4$ is as follows.
We will begin with the image of $\b_q$ in
(the automorphisms of) $V'_4$. Recall that the latter module is
defined for any $q$ and $t=-q^{-1}$. The module in $(\al_N^{-2})$
is the specialization of $V'_4$ for $t=q^{\frac{N-2}{2}}$ and a
primitive $(2\!N)$th root of unity $q^{1/2}$.
We set $t^{1/2}=\imath q^{-1/2}$ for $\imath=\sqrt{-1}$.

One has $V'_4=\mathscr{X}/(e_{-2})=\oplus_{i=0}^3
\Z_{q;t}\,e_i'$, where $e_{-2}=(X^2+1)(X^{-2}-t)$,
\begin{align}\label{v4prime}
&e_0\,=1,\ \,e_1\,=X, \ \,e_{-1}\,=\, X^{-1}+\frac{1-t}{2}X,\ \,
e_2\,=X^{2}+\frac{1-t^{-1}}{2},\notag\\
&Y(e_0)=t^{1/2},\, Y(e_1)=-\imath e_1,\ Y(e_{-1})=\imath e_{-1},\
Y(e_2)=-t^{1/2} e_2,\notag\\
&\ \ \ \ \ T(e_0)=t^{1/2}e_0,\
T(e_1)=\frac{t^{1/2}-t^{-1/2}}{2}e_{1}+t^{-1/2}e_{-1},\\
&T(e_{-1})=\frac{t+t^{-1}+2}{4t^{-1/2}}e_{1} +
\frac{t^{1/2}-t^{-1/2}}{2}e_{-1},\ \, T(e'_2)=t^{1/2}e'_2, \notag\\
& T(e'_{-2})=-t^{-1/2}e'_{-2},\ \,Xe_0=e_1,\
Xe_{1}\,=\,-\frac{1-t^{-1}}{2}e_0+e_2,\notag\\
&Xe_{-1}\,=\,\frac{t+t^{-1}+2}{4}e_0+\frac{1-t}{2}e_2,\ \,
Xe_2'=t^{-1}e_{-1}'\ (\hbox{in \,} V_4'\,)\,.\notag
\end{align}

We see that $Y^2$ is scalar in the even and
odd subspaces of $V'_4$, which are
$\Z_{q;t}\,e_0'\oplus \Z_{q;t}\,e_2'$ and
$\Z_{q;t}\,e_1'\oplus \Z_{q;t}\,e_{-1}'$.
Since $T$ and $X^2$ preserve these subspaces,
we obtain that $TY^2=Y^2 T$ and $X^2 Y^2=Y^2 X^2$.
Applying $\si^{-1}$ from (\ref{tautautau})
in $V'_4$ (this module is projective $PSL_2(\Z)$\~invariant),
we conclude that $TX^2=X^2T$. Therefore
$\,Image(\b_q)$ contains the commutative
subgroup $\n\equal<\!T,X^2,Y^2,q^{1/2}\!>$. It is normal due
to the following relations:
\begin{align}\label{xtxminus4}
&XTX^{-1}=T^{-1}X^{-2},\ YT^{-1}Y^{-1}=TY^{-2},\\
&Y^{-1} X^2 Y=q X^{-2},\ X^{-1}Y^2 X=q^{-1}Y^{-2}.\notag
\end{align}
The latter two relations result
from $Y^{-1}XY=q^{1/2}XT^2$ and commutativity of $X^2,Y^2$
with $T$.  For instance,
\begin{align*}
&Y^{-1}X^2Y\,=\,(q^{1/2}XT^2)^2\,=\, qX(T^2 X)T^2\\
=\ &qX(TXX^{-2}T^{-1})T^2\,=\,qX(X^{-3}T^{-2})T^2\,=\,qX^{-2}.
\end{align*}

Due to $X^{-1}Y^{-1}XY=q^{1/2}T^2$ and relations
from (\ref{xtxminus4}), the commutator subgroup of  $\,Image(\b_q)$ 
(the span of all group commutators) equals
$\c\,=\,<\!q^{1/2}T^2, q^{1/2}Y^2, q^{-1/2} X^2\!>\,=\,
<\!\tilde{T}^2, \tilde{Y}^2, \tilde{X}^2\!>\,$
in notations from (\ref{braidnorm}).

Note that the switch from $\,Image(\b_q)$ to $\,Image(\b_1)$
and correspondingly from $X,Y,T$ to
$\tilde{X},\tilde{Y},\tilde{T}$ will make the relations
from (\ref{xtxminus4}) without $q$. For instance, one has
$$
\ \tilde{Y}^{-1} \tilde{X}^2 \tilde{Y}=
\tilde{X}^{-2},\ \tilde{X}^{-1}\tilde{Y}^2 \tilde{X}=\tilde{Y}^{-2}.
$$

We conclude that $\,Image(\b_q)$ is an extension
of $\n$ by $\Z_2^2$. Accordingly,
$\,Image(\b_1)$ is an extension of $\c$ by $\Z_2^3$.
Also, the center of $\,Image(\b_q)$ belongs to
$\n$ (actually, it is trivial in $\, Image(\b_1)$; see below).
 Indeed, any product of $X,Y,T$ can be written
as $X^\ep Y^\de Z$ for $Z\in \n$, $\ep\in \{0,1\}\ni \de$;
thus, $\ep=0=\de$ for central elements, since $X^\ep Y^\de$
commutes with $\n$ only when  $\ep=0=\de$.

Let $V^{even}=\Z_{q;t}\,e'_0
\oplus\Z_{q;t}\,e'_{2},\
V^{odd}=\Z_{q;t}\,e'_1
\oplus\Z_{q;t}\,e'_{-1}$.
One can check that in this basis
\begin{align*}
X^2_{even}=\begin{pmatrix} \frac{t^{-1}-1}{2}
& \frac{t+t^{-1}+2}{4t}  \\
1 & \frac{t^{-1}-1}{2}\ \   \\ \end{pmatrix},\
X^2_{odd}=\begin{pmatrix} \frac{t^{-1}-1}{2}
& \frac{t+t^{-1}+2}{4}  \\
t^{-1} & \frac{t^{-1}-1}{2} \ \ \\ \end{pmatrix}.
\end{align*}
Therefore $(TX^2)_{odd}=t^{-1/2}=(X^2T)_{odd}$.

Since $T=t^{1/2},Y^2=t$ in $V^{even}$
and $X^2$ has eigenvalues $-1,t^{-1}$ there,
$X^{2M}$ in $V^{even}$ can be a product of powers
of $T$ and $Y^2$ there if and only if
$(-t)^{M}=1=q^M$. This results in $X^{2M}=(-1)^M$.
Similarly, $T$ in $V^{odd}$ has eigenvalues $t^{1/2},-t^{-1/2}$
and $T^M$ in $V^{odd}$
can be a product of powers of $TX^2=t^{-1/2}$ and $Y^2=-1$
there
if and only if $(-t)^{M}=1=q^M$. Thus $T^M=t^{-M/2}$.
Finally, the relation $X^{2a}Y^{2b}T^c=1$, where at least one
of $a$,$b$,$c$ is nonzero,
implies $(-t)^N=1=q^N$ for some $N\in \N$.

Continuing, let $q$ be a primitive $N$th root of unity. We conclude
that $N\mid a$ and $N\mid c$, which makes $X,T$ scalars in $V$.
Since the eigenvalues of $Y^2$ in $V$ are $-1,t$, we obtain that
$N\mid b$ and, finally, the group $\,Image(\n)$ modulo scalars is
isomorphic to $\Z_N^3$. Thus the order of $\,Image(\b_q)$
modulo scalars
is $4N^3$ in this case. If there is no such $N>0$,
then $\n$ is isomorphic to $\Z^3$.

Finally, the scalars form the center of $\, Image(\b_q)$
due to the irreducibility of $V'_4$. The center is
always $\{q^{\Z/2}\}$ (either for $N>0$ or if $q$ is not
a root of unity); it is trivial for the group
$\, Image(\b_1)$. This is checked as above.

\smallskip
{\sf Compatibility for different $N$.}
\comment{
Let us discuss a
connection of $V'_4$  considered for a primitive $(2\!N)$th
root of unity $q^{1/2}$  and the module $\check{V}'_4 $
of $\check{\HH}\equal \HH_{q^{d/2},t^{d/2}}$ obtained upon
the substitution
$\check{q}^{1/2}\equal q^{d/2}$.
}
Let $d|N$ and $\check{N}=N/d>2\,$ for
odd $d\,$. Then $\check{V}'_4$
is irreducible for $\check{q}^{1/2}\equal q^{d/2}$ and
$\check{\HH}\equal \HH_{q^{d/2},t^{d/2}}$.
The quotient of $\,Image(\b_q)$ by the center for
a primitive $(2\!N)$th root of unity $q^{1/2}$ and that
for $\check{q}^{1/2}$ will be denoted by
$\mathfrak{B}^\dag_N$  and $\mathfrak{B}^\dag_{\check{N}}$.

Then $\mathfrak{B}^\dag_N$
is the extension of the commutative
subgroup
\begin{align}\label{nNdag}
\n_N^\dag=<\!T,X^2,Y^2\!>\!/\!<\!T^N,X^{2N},Y^{2N}\!>\simeq \Z_N^3
\end{align}
by $\Z_2^2=<\!X,Y\!>/<\!X^2,Y^2,XYX^{-1}Y^{-1}\!>$ subject
to the relations
\begin{align}\label{xtxminus4dag}
Y^{-1}XY=X&T^2=YXY^{-1},\ \ \,
XTX^{-1}=T^{-1}X^{-2}=X^{-1}TX,\notag\\
YT^{-1}Y^{-1}=T&Y^{-2}=Y^{-1}T^{-1}Y,\ \
Y^{-1} X^2 Y=X^{-2}=Y X^2 Y^{-1},\notag\\
&X^{-1}Y^2 X=Y^{-2}=X Y^2 X^{-1}.
\end{align}

Furthermore,  the following two
equivalent relations hold in $\b_q^\dag$
for any $q\in \C^*$ and odd $m=2l+1$:
\begin{align}\label{Cformula}
(Y^{m}\!X^{-m})^2\!=\!Z\equal Y^2T^{-2}\!X^{-2},\
Y^{-m}\!X^{-m}Y^{m}\!X^{m}\!=\!
ZY^{-2m}\!X^{2m}.
\end{align}
Indeed, using (\ref{xtxminus4dag}) and the commutativity of $\n$,
\begin{align*}
Y^{-m}&X^{-m}Y^{m}X^{m}=(Y^{-m}X^{-2l}Y^{m})(Y^{-m}X^{-1}Y^{m})X^m\\
&=X^{2l}(T^{-2}Y^{-4l}X^{-1})X^{2l+1}=X^{4l}T^{-2}Y^{-4l},
\hbox{ \, where\, }\\
Y^{-m}&X^{-1}Y^m\!=\!Y^{-2l-1}X^{-1}Y^{2l+1}\!=\!
Y^{-1}(Y^{-4l}X^{-1})Y\!=\!Y^{-4l}T^{-2}X^{-1}.
\end{align*}

We note that $XZX^{-1}=Z^{-1}=YZY^{-1},\ TZ=ZT$ in $\b_q^\dag$.
Also, $Z=q^{-1/2}\tau_-(X^2)$;\,
$Z=X^{-2}$ in $V^{even}$  and $Z=-t^{1/2}T^{-1}$
in $V^{odd}$. The eigenvalues of $Z$ in each of these two spaces
coincide and are $\{-1,t\}=\{-1,-q^{-1}\}$; recall that
$t=-q^{-1}$. Actually, the  eigenvalues of $Z$ are
important here only up to a common factor.


Next, $Y^mT^{-m}Y^{-m}=YT^{-m}Y^{-1}=T^mY^{-2m}$ and
$Y^{-m}X^{2m}Y^{m}=Y^{-1}X^{2m}Y=X^{-2m}$ for odd $m$, which
follows from (\ref{xtxminus4dag}) and the commutativity
of $\n$. The same holds for $X^{-1}$ instead of $Y$.
Finally, setting $\check{Z}=Y^{2d}T^{-2d}X^{-2}$ (see
(\ref{Cformula})),
we obtain the
homomorphism
\begin{align}\label{NcheckN}
&\mathfrak{B}^\dag_N/\!<\!Z\!>\,\to\,
\mathfrak{B}^\dag_{\check{N}}/\!<\!\check{Z}\!>,\ \,
X\mapsto X^d, Y\mapsto Y^d,
T\mapsto T^d.
\end{align}


\setcounter{equation}{0}
\section{\sc The Galois action}
\subsection{\bf Deligne-Simpson Problem}
Following our approach to the
Tate modules from the Introduction
and Section 2.7.3 from \cite{C101}, let as
switch to the $A,B,C$\~generators of $\b_q$.

Setting $A=XT,\, B=XTY,\,  C=T^{-1}Y$,
the relations of $\b_q$ and the action of $\tau_{\pm}$
become
\begin{align}\label{ABCinter}
A^2\,=\,1&\,=\,C^2=q^{1/2}B^2,\, \hbox{\, where \,}
ABC=A^2 YT^{-1}Y=
YY^{-1}T=T,\notag\\
&\tau_+:\ A\mapsto A,\ B\mapsto q^{-1/4}C,\ C\mapsto
q^{1/4}C^{-1}BC,
\\
&\tau_-:\  A\mapsto q^{1/4}ABA^{-1}, \
B\mapsto q^{-1/4}A,\  C\mapsto C. \notag
\end{align}
Also, the automorphisms $\varsigma_x,
\varsigma_y$ from (\ref{fvarsx}),(\ref{fvarsy})
 act as follows:
\begin{align}\label{varsiabc}
&\varsigma_x: A\!\mapsto\! -A,B\!\mapsto\! -B,C\!\mapsto\! C,\ \,
\varsigma_y: A\!\mapsto\! A,B\!\mapsto\! -B,C\!\mapsto\! -C.
\end{align}

Let us now fix a primitive root of unity $q^{1/2}$ of order $2\!N$.
We take an $\HH\,$\~modules
$V$ that is  {\em dim-rigid\,} and {\em loc-dim-rigid\,}
from Theorem \ref{psl-two-inv} and  pick prime $p$
such that \,gcd$(p,2\!N)=1$. Then we choose
a prime ideal $\mathfrak{p}$ in  $\Z_{\ddot{q}}$
over $(p)$ for $\ddot{q}=q^{1/4}=(q^{1/2})^{1/2}$
for even $2k$ and  $\ddot{q}=q^{1/2}$ if\, $k\in 1/2+\Z$.
Also, all {\em boundary\,} modules from
(\ref{boundcasea}) will be allowed (though not all of
them are {\em dim-rigid}).

Similarly, we can take a {\em dim-rigid\,}
$\HH^\checkmark\,$\~module $V^\checkmark$ for odd $N$
from Proposition \ref{ODDN}. In this case the condition
\,gcd$(p,N)=1$ is sufficient and $\,\ddot{q}\,$ must be
modified as follows:
$$
\ddot{q}^{\,\checkmark}=q^{1/2}
\hbox{\, for even } 2k \hbox{\ \, and\, } \ddot{q}^{\,\checkmark}=q
\hbox{\, if } k\in 1/2+\Z.
$$

Let $\mathfrak{B}_{\mathfrak{p},m}=\mathfrak{B}(q,t,\mathfrak{p},m)$
be \,$Image(\b_q)$ and $\tilde{\mathfrak{B}}_{\mathfrak{p},m}=
\,Image(\b_1)$ from (\ref{braidnorm}). The groups
$\mathfrak{B}$ and
$\tilde{\mathfrak{B}}$ are the images of $\b_q$ and $\b_1$
without the reduction modulo
$\mathfrak{p}^m$. The group $\b_1$ is obtained from $\b_q$ by
sending $q^{1/4}\mapsto 1$ in all relations. It is generated by
\begin{align}\label{brdnrm}
&\tilde{T}\,=\,q^{+1/4}T,\
\tilde{Y}\,=\,q^{+1/4}Y,\
\tilde{X}\,=\,q^{-1/4}X\hbox{ or }\notag\\
&\tilde{A}=XT,\,\ \tilde{B}=q^{+1/4}XTY,\,\  \tilde{C}=T^{-1}Y.
\end{align}
For instance, relations (\ref{dim2fin}) for $V_2$ and
$t^{1/2}q^{1/4}=\pm\imath\,$ read
\begin{align}\label{dim2finx}
&\tilde{X}=\begin{pmatrix} 0 & q^{1/4} \\ q^{-1/4} & 0 \\
\end{pmatrix},\ \
\tilde{Y}=\begin{pmatrix} \pm \imath & 0  \\ 0 &
\mp \imath q^{1/2} \\
\end{pmatrix}, \ \ \, \tilde{T}=\pm \imath\,,\\
&\tilde{A}\,=XT=\begin{pmatrix} 0 & t^{-1/2} \notag\\
t^{1/2} & 0 \\ \end{pmatrix},\ \,
\tilde{C}\,=T^{-1}Y=\begin{pmatrix} 1 & 0 \\
0 & -1 \\ \end{pmatrix},\\
&\tilde{B}=q^{+1/4}XTY=q^{+1/4}
\begin{pmatrix}  0 & -1  \\ t & 0 \\ \end{pmatrix}=
\begin{pmatrix}  0 & \mp \imath t^{-1/2}
\\ \pm \imath t^{1/2} & 0 \\ \end{pmatrix}.\notag
\end{align}
Thus the matrices $\tilde{A},\tilde{B},\tilde{C}$ for $V_2$
are conjugated to the classical Pauli matrices
$\si_1,\si_2,\si_3\,$ for $\,\pm\!=\!+$\, (plug in here
$t^{1/2}=1$).
\smallskip

Thus $\tilde{\mathfrak{B}}_{\mathfrak{p},m}$ is represented by
invertible matrices
of dimension dim$V$ or dim$V^\checkmark$ with entries in
$\Z_{\ddot{q}}=\Z[\ddot{q}]/\mathfrak{p}^m$ for $p>0$ and
with entries in $\Q(\ddot{q})$ if $p=0=m$ (by definition).
Here $\ddot{q}^{\,\checkmark}$ must be used for $V^\checkmark$,
which is the square of $\,\ddot{q}\,$.
\smallskip

Accordingly, $\mathfrak{B}_{(p),m}$ and
$\tilde{\mathfrak{B}}_{(p),m}$ are defined for the image
of $\b_q,\b_1$ in
$$\hbox{GL}\bigl(V\otimes_{Z_{q'}}
\Z_{q'}/(p^m)\bigr)\!=\!
\prod_{\mathfrak{p}\supset (p)}\hbox{GL}
\bigl(V\otimes_{\Z_{q'}}
\Z_{q'}/\mathfrak{p}^m\bigr),\  q'=\dot{q},\ddot{q}
\hbox{\, for\, } \mathfrak{B},\tilde{\mathfrak{B}};
$$
where $\mathfrak{p}$ are prime ideals in $\Z_{q'}$ over
$(p)$,\, gcd$(p,2\!N)=1$. This is unless for $p=0$, where the
definition field is $\Q(q')$. One can take here $V^\checkmark$
for odd $N$ provided \,gcd$(p,N)=1$; then $q'$ must be replaced by
$(q')^2$.
\smallskip

From now on, we will mainly switch from $V$ to $\tilde{V}$
and from $V^\checkmark$ to $\tilde{V}^\checkmark$; the
tilde-modules are (by definition) the same vector spaces
supplied with the action of $\b_1$ via the generators
$\tilde{A},\tilde{B},\tilde{C},\tilde{T}$. The matrix
entries of these generators can be assumed to
belong to $\Z_{\ddot{q}}$  or $\Z_{\ddot{q}^{\,\checkmark}}$.

Also, unless stated otherwise, we will always perform the
reduction
modulo $\mathfrak{p}^m$ or
$(p^m)$ except for the finite $\,Image(\b_q)$
(see Theorem \ref{FINBIMAGE}).
In the latter case, we set $p=0=m$
and consider $\tilde{V}$ and $\tilde{V}^\checkmark$
over  $\Q(\ddot{q})$ and  $\Q(\ddot{q}^{\,\checkmark})$.
Thus $\tilde{A},\tilde{B},\tilde{C},\tilde{T}$
will be considered below over $\Z_{\ddot{q}}/\mathfrak{p}^m$ for
$m\in \N$ or over $\Z_{\ddot{q}}/(p^m)$ (and with
$\ddot{q}\mapsto \ddot{q}^{\,\checkmark}$ for
$V^\checkmark$).
\smallskip

{\sf Deligne-Simpson problem.}
Using (\ref{ABCinter}), let us reformulate Theorem \ref{psl-two-inv}
as the multiplicative (irreducible) {\em Deligne-Simpson
problem\,}, abbreviated as {\em\,DSP\,}.
See \cite{Kos} for a survey on {\em DSP\,}.

We will mainly focus in this section on the case without
$\checkmark$
(when $q^{1/2}$ is primitive of order $2\!N$).
Though see Part $(ii)$ of
Theorem \ref{ABC-DSP}) below on the modules $V^\checkmark$ for odd 
$N$; recall that one must change $q,t$ to $q^2,t^2$ and
square $\ddot{q}$ in this case.

We will use $\,\sim\,$ for the
equivalence of matrices over the field $\Q(\ddot{q})$ for
$\ddot{q}=q^{1/4}$ (even $2k$) and
$\ddot{q}=q^{1/2}$ (odd $2k$)
or, correspondingly, over
$\Z_{\ddot{q}}/\mathfrak{p}^m$ for prime $\mathfrak{p}\supset (p)$
when $p>0$ and  gcd$(p,2\!N)=1$. We will denote by
$D_{\la,\mu}^{(l)}$,
the diagonal matrix with $\, l$ diagonal entries $\la$ and the
diagonal entries $\mu$ otherwise. We take $q^{1/2},
t^{1/2}=\pm (q^{1/2})^k$ as
in Theorem \ref{psl-two-inv} for the series $(\al,\be,\breve{\ga})$.

\begin{theorem}\label{ABC-DSP}
(i) Let $M=$dim$V=2N\pm 4|k|, 2|k|$ for $V$ from  $(\al,\be,
\breve{\ga})$ of
Theorem \ref{psl-two-inv} (without applying $\iota$ for the sake
of definiteness); we
set $\,n=[(M+1)/2]$, $n'=[M/2]+1$. We will assume that
$\,t\neq -1$ to ensure the semisimplicity of $T$.
The following irreducible {\em DSP\,} in invertible
$M\times M$\~matrices
over $\Q(\ddot{q})$ or $\Z_{\ddot{q}}/\mathfrak{p}^m$ (as
above) has a unique solution up conjugation by a (common)
invertible matrix:
\begin{align}\label{ABCT}
\tilde{A}\sim& D_{1,-1}^{(l_a)},\ \tilde{B}\sim D_{1,-1}^{(l_b)},
\,\tilde{C}\sim D_{1,-1}^{(l_c)},\
\tilde{T}=\tilde{A}\tilde{B}\tilde{C}\sim D_{q^{1/4}t^{1/2},
-q^{1/4}t^{-1/2}}^{(n')}\,,\notag\\
&\hbox{where\, }\ l_{a,b,c}=n\hbox{\, for the modules from\, }
(\al),(\be) \hbox{\, and \,} V_{2|k|}^+\,, \hbox{ or}\\
&l_a\!=\!n\!-\!1\!=\!l_b\!=\!l_c\!-\!1,\ l_b\!=\!n\!-\!1\!=\!l_c\!=
\!l_a\!-\!1,\
l_a\!=\!n\!-\!1\!=\!l_c\!=\!l_b\!-\!1\notag\\
&\hbox{correspondingly for the modules }\,
\varsigma_x(V_{2|k|}^+),\
\varsigma_y(V_{2|k|}^+),\
\varsigma_x\varsigma_y(V_{2|k|}^+).
\notag
\end{align}
Furthermore, there are no irreducible solutions of (\ref{ABCT})
with such $M$ and $k$
for any other combinations of $\,l_{a,b,c},n,n'.$

(ii) In the cases
$(\al^\checkmark,\be^\checkmark,\ga^\checkmark)$ for
odd $N$, the dimension $M$ is correspondingly $N\pm 2|k|$
and $2|k|$, except for $(\al^\checkmark)$ with odd $2k$,
when dim$=\!2N\!-\!4k$ and the action of $X$,$Y$ in
$V_{2N-4k}^{\checkmark}$ becomes non-semisimple.
The eigenvalues of $\tilde{T}^\checkmark$
are obtained by squaring $\,t$ ($t^2\neq -1$ since $N$ is
odd). In parallel with
$(i)$, $\,l_a=l_b=l_c=n\,$ for the modules
$V^{+\checkmark}_M$. The transformations
of $l$\~multiplicities under applying
$\varsigma_x,\varsigma_y$ to
$V_{2N-4k}^{+\checkmark}$ for even $2k$ and
to $V_{2N+4k}^{+\checkmark}$ are described exactly
as for $\varsigma_x^\ep\varsigma_y^\de(V_{2|k|})$ in
(\ref{ABCT}).
Accordingly, the $\checkmark$\~variant of the DSP from (\ref{ABCT}) 
has solutions only for such $\,l_{a,b,c},n,n'$.
\sq
\end{theorem}

Let us briefly comment on $V_{2|k|}^+$ from $(\ga)$;
its dimension is $M=2n-1$ and $l_a,l_b,l_c=n$.
The fact that the outer automorphisms $\varsigma_x$,
$\varsigma_y$, and $\varsigma_x\varsigma_y$
(when applying to $V_{2|k|}^+$) correspondingly
diminish $l_a\!=\!l_b$, $l_b\!=\!l_c$ and
$l_a\!=\!l_c$ by $1$
follows directly from the definition of
$\varsigma_x,\varsigma_y$.
The same holds for the  $\checkmark$-modules
except for $(\al^\checkmark)$ and odd $2k$.

Note that the DAHA-parameters $q^{1/2},t^{1/2}$
are not uniquely recoverable from this {\em DSP\,}. By altering
the sign of $q^{1/4}$, one can change the sign of
$t^{1/2}$, which exactly corresponds to applying $\iota$.
Recall that $\,\iota$ can be interpreted as
the substitution $k\mapsto -k$, sending
$q^{1/4}t^{1/2}\mapsto q^{1/4}t^{-1/2}$. This transformation
changes $n'$ to $M-n'$ and does not influence
$l_{a,b,c}$.

Also, the irreducibility of {\em DSP\,} in (\ref{ABCT})
can be omitted
for $q,t=q^k$ and $M$ as in $(\al,\ga)$. The minimality
of dimensions of modules $V_{2n-4k}$ and $V^+_{2|k|}$
gives that $\tilde{A},\tilde{B},\tilde{C}$ have no nontrivial
invariant subspaces.
\smallskip

We see that, generally, irreducible $\HH$\~modules
provide  solutions of {\,\em DSP\,} for the quiver
of type $D_4^{(1)}$.
The other way around, one can use this {\em DSP\,} for finding
irreducible DAHA-modules. See Theorem 1 from \cite{ObS} devoted
to the $C^\vee C_1$\~ case
(apart from the roots of unity), which is directly
related to Theorem 1.5 from \cite{CrB}.
This paper is based on the roots of unity.
Any rigid modules for generic $q$ remain rigid at
sufficiently general roots of unity, but we need more
than these.

The cases of non-semisimple $T$ for $t=-1$,
including the boundary case $(\be^\bullet)$, can be
readily managed as well; this requires
using Jordan normal forms in $DSP\,$ as in \cite{CrB}.
Note that $t^\checkmark=t^2=q^{2k}$ cannot be $-1$,
since $N$ is odd; so $T^\checkmark$ is always semisimple.
Generally, there is no connection with non-semisimplicity
of $T$ and that for $X$ and $Y$, which occurs for
$(\be),(\be^\checkmark)$ and for $(\al^\checkmark)$ when
$2k$ is odd.


\subsection{\bf Absolute Galois group}
The passage to $\tilde{A},
\tilde{B},\tilde{C}$  identifies
$\b_1$, considered as the orbifold fundamental group of
$(E_{\C}\setminus\{0\})/\S_2$
for
an elliptic curve $E_{\C}$ (over $\C$) punctured at $0$, 
with the standard
fundamental group
$\pi_1(P_{\C}^1\setminus \{\al_0,\al_1,\al_2,\al_3\},b_o)$. See
\cite{C101}, Proposition 2.7.6.
The base point $b_o$ here and below is sufficiently general;
$s\in \S_2$ is the reflection $x\mapsto -x$ in $E_{\C}$. The
normal subgroup $<\!X,Y,T^2\!>\subset\b_1$ of index $2$ is
associated
then with the corresponding covers of $E_{\C}$.

In the main theorem, $E\ni 0$ will be a (smooth projective)
elliptic curve and its zero defined over the field
$K\subset \overline\Q \subset \C$. The points
$\{\al_0,\al_1,\al_2,\al_3\}$ will be the images
of the points $\{0, 0_1, 0_2, 0_3\}$ of the $2$nd
order in $E$ under the isomorphism
$E/\S_2 \to P^1$ over $K$.

The {\em Riemann Existence Theorem, RET\,},
provides the existence of a (connected) Riemann surface
$F_{\C}$ such that it is a {\em Galois cover} of $P^1_{\C}$
ramified at points $\al_0,\al_1,\al_2,\al_3$ and
with the {\em Galois group} Aut$(F_{\C}/P^1_{\C})$
isomorphic to
$\tilde{\mathfrak{B}}_{\mathfrak{p},m}$. Furthermore,
$\tilde{T}\!=\!\tilde{A}\tilde{B}\tilde{C},
\tilde{A},\tilde{B},\tilde{C}$ can be identified
with the generators of the (cyclic) stabilizers
of certain ramification points $P_0,P_1,P_2,P_3\in F_{\C}$
over $\al_0,\al_1,\al_2,\al_3$. Galois covers
are those with transitive action of  Aut$(F_{\C}/P^1_{\C})$
in the fibers.
See e.g. Theorem 2.13 from
\cite{Vol} and also Proposition 4.23 there.

The isomorphism classes of such
covers together with the identification
of \,Aut$(F_{\C}/P^1_{\C})$ with the group
$\tilde{\mathfrak{B}}_{\mathfrak{p},m}$
are in one-to-one correspondence with
the triples of generators $\{A',B',C'\}$ of
$\tilde{\mathfrak{B}}_{\mathfrak{p},m}$
considered up to a conjugation and satisfying
the following. They must be taken from
the corresponding conjugacy classes of
$\tilde{A},\tilde{B},\tilde{C}$ and
$T'=A'B'C'$ must belong to the conjugacy class
of $\tilde{T}$ in $\tilde{\mathfrak{B}}_{\mathfrak{p},m}$.
\smallskip

\rmk
The set of
such isomorphism classes
is denoted by $\e_4^{in}$ in
\cite{FrV}; see also Section 3.1 in \cite{Det}.
Actually  $\e_4^{ab}$ is sufficient in this paper, where
the triples $\{A',B',C'\}$ are considered modulo the
action of the group of all automorphisms
of $\tilde{\mathfrak{B}}_{\mathfrak{p},m}$, not only inner.
\sq
\smallskip

The key here is the following. Let us choose
an ordered {\em standard homotopy basis\,}
in $\pi_1(P_{\C}^1\setminus \{\al_i\},b_o)$,
where $\ga_i$ is
represented by a loop at $b_o$ that winds once around $\al_i$
counterclockwise for $i=1,2,3$ and clockwise for
$i=0$ and winds around no other $\al_j$. It is of course not unique;
see e.g. \cite{CoH}. Then
{\em \, RET\,} provides that there is a canonical
Galois cover of $P^1_{\C}$ together with the
identification of its Galois group with
$\tilde{\mathfrak{B}}_{\mathfrak{p},m}$
such that the generators
$A',B',C'$ in the latter group correspond to
$\ga_1,\ga_2,\ga_3$ and $T'=A'B'C'$ corresponds
to $\ga_0$.
\smallskip

Switching to the fields of rational functions,
$\C(F)/\C(P^1)$ is the Galois extension of the fields of rational
functions on $F$ and $P^{1}$. Therefore\,
Gal$(\C(F)/\C(x))\toeq
\tilde{\mathfrak{B}}_{\mathfrak{p},m}$, where we
set  $\C(P^1)=\C(x)$.
\smallskip

Since $\C(F)/\C(P^1)$ is Galois, there exists the minimal field
of definition of this extension
of fields, which contains the field of definition of
the set of points $\{\al_0,\al_1,\al_2,\al_3\}$.
See Proposition 2.5 from \cite{CoH},
Section 1.5 from \cite{FrV} (and references therein)
and Proposition 7.12 from \cite{Vol}.
Let us assume that the definition field is a field of algebraic
numbers $K\subset \overline{\Q}$. We will also assume that
$\al_0$ and the base point $b_o$ are defined over $K$,
though $\al_1,\al_2,\al_3$ may be permuted by
Gal$(\overline{\Q}/K)$.

\subsection{\bf Main Theorem}
We  fix a primitive root of unity $q^{1/2}$ of
order $2\!N$ and
a {\em dim-rigid\,}  module $\,V\,$
from Theorem \ref{psl-two-inv}, which is also {\em loc-dim-rigid\,}.
Recall that
$\tilde{V}$ is $V$ supplied with the action of $\b_1$ via
$\tilde{A},\tilde{B},\tilde{C}$, which
can be assumed with the entries in $\Z_{\ddot{q}}$
for $\ddot{q}=q^{1/4}=(q^{1/2})^{1/2}$
for even $2k$ and  $\ddot{q}=q^{1/2}$ if $k\in 1/2+\Z$.
We will consider
$\tilde{V}$ over $\Z_{\ddot{q}}/\mathfrak{p}^m$ for $m\in \N$
unless the image $\,Image(\b_q)$ of $\b_q$ in GL$(V)$ is finite
(see Theorem \ref{FINBIMAGE}). In the latter case, we set $p=0=m$
and consider $\tilde{V}$ over  $\Q(\ddot{q})$.
Also, see Part $(v)$
of Theorem \ref{GALACTION} concerning the case of modules
$V^\checkmark$ for odd $N$.
\smallskip

Let $G_K=$Gal$(\overline{\Q}/K)$ be the absolute Galois
group of $K$; we also set $G_K^\circ=$Gal$(\overline{\Q}/K^\circ)$
for $K^\circ=K(\ddot{q})$.
\smallskip

The main application of the rigidity of the
DAHA-modules $\,V\,$ under consideration in this paper
is the action of $G_K$
and $G_K^\circ$ in the groups
$\tilde{\mathfrak{B}}_{(p),m}$ and
$\tilde{\mathfrak{B}}_{\mathfrak{p},m}$, which are the images
of $\b_1$ in the groups of automorphisms
of $\tilde{V}$ over the rings $\Z_{\ddot{q}}/(p^m)$ and,
correspondingly, over $\Z_{\ddot{q}}/\mathfrak{p}^m$.
The ring of definition is the field $\Q(\ddot{q})$ if $p=0$;
this is the case of finite $\tilde{\mathfrak{B}}$, which will
denote the image of $\b_1$ for $V$ in characteristic $0$.

The justification of the following theorem is actually
similar to the proof of projective $PSL_2(\Z)$\~invariance of
{\em unirigid\,} modules. Recall that
$\tilde{V}^{sym}=\{v\in V:\, T(v)=t^{1/2}v\}=$
$\{\tilde{v}\in \tilde{V}:\, \tilde{T}(v)=q^{1/4}t^{1/2}\tilde{v}\}$.
However the $G_K$\~invariance holds even for the modules
$V^+_{2|k|}, V^{+\checkmark}_{2|k|}$
and for $V^{+\checkmark}_{N\pm 2|k|}$, which are not
{\em dim-rigid\,}. It even works in the case
of $V_{4N}$ from (\ref{boundcasea}), which is not
dim-rigid. In the non-unirigid cases
we need to use the {\em DSP} at greater potential, which
always appears sufficient. For instance, the
Jordan form of $\tilde{T}$ is used for $V_{4N}$.

\begin{theorem}\label{GALACTION}
Let us consider any module $\tilde{V}$ over $\Z_{\ddot{q}}$
of type  $(\al,\be,\breve{\ga})$ from Theorem \ref{psl-two-inv}
(including applying $\iota$ there)
and pick prime $p$ such that\, gcd$(p,2\!N)=1$ and a prime
ideal $\,\mathfrak{p}\subset \Z_{\ddot{q}}$ over $(p)$.
For $m\in \N$ (we set $m=0$ for $p=0$), let
$f:F_{\C}\to P_{\C}^1$ be the Galois cover such that
\,Gal$(\C(F)/\C(P^1))\toeq
\tilde{\mathfrak{B}}_{\mathfrak{p},m}$, where
$\tilde{T}\!=\!\tilde{A}\tilde{B}\tilde{C},\,
\tilde{A},\tilde{B},\tilde{C}$ generate
the ramification (cyclic) subgroups at certain $P_0,P_1,P_2,P_3$
in the fibers of $\,f\,$ over
$\al_0,\al_1,\al_2,\al_3\in P^1$.

(i) The cover $F_{\C}\to P_{\C}^1$ and the
field extension $\C(F)/\C(P^1)$
can be defined over $K^\circ=K(\ddot{q})$.
Switching to the group $\tilde{\mathfrak{B}}_{(p),m}$, which
has natural projections onto all groups
$\tilde{\mathfrak{B}}_{\mathfrak{p},m}$ for
$(p)\subset \mathfrak{p}$,
the cover  $f^\diamond:F^\diamond_{\C}\to P_{\C}^1$ introduced
for this group in the same way as above and the corresponding
field extension $\C(F^\diamond)/\C(P^1)$ can be
defined over $K$. Accordingly, one has homomorphisms
$G_{K}^\circ\to
\hbox{Aut}\bigl(\tilde{\mathfrak{B}}_{\mathfrak{p},m}\bigr)$
for any $\mathfrak{p}\supset (p)$
and a homomorphism
$G_K\to \hbox{Aut}\bigl(\tilde{\mathfrak{B}}_{(p),m}\bigr).$

(ii) For each prime ideal $\mathfrak{p}\subset \Z_{\ddot{q}}$
over $(p)$ as above,
the action of $G_K^\circ$ in
$\tilde{\mathfrak{B}}_{\mathfrak{p},m}$ is by inner automorphisms
of $\tilde{V}$ considered over $\Z_{\ddot{q}}/\mathfrak{p}^m$. The
action of any $\,g\in G_K$ in
$\tilde{\mathfrak{B}}_{(p),m}$ is by inner automorphisms
of $\tilde{V}$ considered over $\Z_{\ddot{q}}/(p^m)$
(or $\Q(\ddot{q})$ for $p=0$) composed with the Galois conjugation
$\si_g\in \hbox{Gal}(\Q(\ddot{q})/\Q)$ acting naturally in
\,GL$(\tilde{V})$ considered over $\Z_{\ddot{q}}/(p^m)$.
Then the map $g\mapsto \si_g$ coincides with the natural restriction
homomorphism $G_K\to \hbox{Gal}(\Q(\ddot{q})/\Q)$.
\smallskip

(iii) Continuing $(ii)$, there
exists a group homomorphism
$\hat{g}\mapsto \phi_{\hat{g}}=h_{\hat{g}}\si_g\,$ for $\,\hat{g}\,$
from a proper central extension $\varrho:\hat{G}_K\to G_K$ and
$g=\varrho(\hat{g})$
such that the action of
$\,g$ in $\tilde{\mathfrak{B}}_{(p),m}$ is the
action of $\si_{g}$ followed by the conjugation by
$h_{\hat{g}}\in$GL$(\tilde{V})$
for $\tilde{V}$
over $\Z_{\ddot{q}}/(p^m)$ (over $\Q(\ddot{q})$ when $p=0$).
Moreover, applying a proper automorphism
from Aut$(F^\diamond\to P^1)$,
one can assume that
$g(\tilde{T})=\tilde{T}^M$, where $\si_g(\ddot{q})=\ddot{q}^M$,
which is equivalent to the relations
$\si_g(T)=T^M,\, h_{\hat{g}}\,T\,h_{\hat{g}}^{-1}=T$.
Then $G_K$ preserves $\tilde{V}^{sym}$ over
$\Z_{\ddot{q}}/(p^m)$; accordingly, $G_K^\circ$
preserves $\tilde{V}^{sym}$ over
$\Z_{\ddot{q}}/\mathfrak{p}^m$.
\smallskip

(iv) Furthermore, $h_{\hat{g}} \,(\hat{g}\in \hat{G}_K)\,$
are unitary with
respect to the inner product $\lan\cdot\,,\,\cdot\ran$
considered over $\Z_{\ddot{q}}/(p^m)$ for types
$(\al,\breve{\ga})$. In the case $(\breve{\ga})$, this claim and
those from $(i,ii,iii)$ hold for the
modules $\tilde{V}=\,\varsigma_x^{\ep}
\varsigma_y^{\de}(V^+_{2n+1})$
for $\{\ep,\de\}\in \!\bigl\{\{1,0\},\{0,1\},\{1,1\}\bigr\}$
associated correspondingly with the pairs  $\{\al_i,\al_j\}\in\!
\bigl\{\{\al_1,\al_2\},\{\al_2,\al_3\}, \{\al_1,\al_3\}\bigr\}$
subject to the following adjustment.
The automorphism $\si_g$ must be multiplied by
$\varsigma_x^{\ep'-\ep}\varsigma_y^{\de'-\de}$
for  the pair of indices $\{\ep',\de'\}$ corresponding to
the pair $g\{\al_i,\al_j\}$ of ramification points in $P^1$.
\smallskip

(v) All claims above hold for the boundary cases
from (\ref{boundcasea}) and for the
$\HH^\checkmark\,$\~modules $V^\checkmark$ for odd $N$
from Proposition \ref{ODDN}, namely  for
$ (\breve{\al}^\checkmark),(\breve{\be}^\checkmark),
(\breve{\ga}^\checkmark)$.
The component $\,\si_g$ in $\,\phi_{\hat{g}}$
must be adjusted here by
$\varsigma_x^{\ep'-\ep}\varsigma_y^{\de'-\de}$
in the same way as it was done in Part $(iv)$ for
$(\breve{\ga})$
in the following cases:\, $(\breve{\al}^\checkmark)$
for even $2k$, $(\breve{\be}^\checkmark)$ and
$(\breve{\ga}^\checkmark)$.
This is not needed for
$(\al^\checkmark)$
for odd $\,2k$, since $\varsigma_x^\ep\varsigma_y^\de$
(and $\tau_{\pm}$) do not change the isomorphism class
of $V_{2N-4k}^\checkmark$, which remains irreducible.
\end{theorem}
\smallskip

{\em Proof.} We will mainly
consider here the cases $(\al,\be,\breve{\ga})$.
The adjustments to the boundary cases from (\ref{boundcasea})
and the cases $(\breve{\al}^\checkmark,\breve{\be}^\checkmark,
\breve{\ga}^\checkmark)$ from Proposition \ref{ODDN}
are straightforward
with one reservation. The boundary module $V_{4N}$
of type $(\be^\bullet)$ is not {\em dim-rigid\,}.
However one can use here formula (\ref{TPDe}) instead
of {dim-rigidity\,}. Note that
when $\tilde{T}$ becomes non-semisimple, the rigidity
argument continues to work.

Let us take a {\em standard homotopy basis\,}
in $\b_1=\pi_1(P_{\C}^1\setminus \{\al_i\},b_o)$.
Then $F_{\C}$ corresponds to a normal subgroup
$U\subset \b_1$ such that
$\b_1/U=\hbox{Aut}(F_{\C}\!\to\! P^1_{\C})\toeq
\tilde{\mathfrak{B}}_{\mathfrak{p},m},
$
where the standard generators $\ga_0,\ga_1,\ga_2,\ga_3$
map to $\tilde{T},\tilde{A},\tilde{B},\tilde{C}.$
The latter four elements are generators of the
(cyclic) stabilizers at certain ramification points
$P_0,P_1,P_2,P_3\in  F_{\C}$ in the
fibers of $\,f'\,$ over the points $\al_0,\al_1,\al_2,\al_3$.
\smallskip

The covering $F/P^1$ can be defined over
$\hat{K}\subset \overline{\Q}$ for a finite Galois extension
$\hat{K}/K^\circ$
(since there are only finitely many covers with the same
Galois group and given ramification points).
For any element $g\in \hbox{Gal}(\hat{K}/K^\circ)$, let
$f'\!:\!F'\!\to\! P^1$ be the covering obtained by
applying $g$  to $F/P^1$.
It has the same ramification
points in $P^1$ and corresponds to the following homomorphism
$\b_1/U'=$Aut$(F_{\C}\to P^1_{\C})\toeq
\tilde{\mathfrak{B}}_{\mathfrak{p},m}$.

It sends the standard generators to another systems
of generators of $\tilde{\mathfrak{B}}_{\mathfrak{p},m}$,
elements $\,\tilde{T}',
\tilde{A}',\tilde{B}',\tilde{C}'$ generating the stabilizers
at $P'_0,\,P'_{\hat 1}\,,P'_{\hat 2}\,,P'_{\hat 3}\,$,
which are the $\,g$\~images of
$P_0,P_1,P_2,P_3$. These points belong to the
fibers of $\,f'\,$  over $\,\al_0,\al_{\hat 1},
\al_{\hat 2},\al_{\hat 3}$\, for the
permutation $\,\{123\}\mapsto
\{{\hat 1}\,{\hat 2}\,{\hat 3}\}$
induced by $\, g$.
Recall that the base point $b_o$, $\al_0$ and the
set $\{\al_1,\al_2,\al_3\}$ are defined over $K$;
so one can set $g(\al_i)=\al_{\,\hat i}$ for $i\ge 1$.

Furthermore,  there exists $M\in \N$ such that
$\tilde{T}',
\tilde{A}',\tilde{B}',\tilde{C}'$ belong to the same
conjugacy classes as
to  $\,\tilde{T}^M, \tilde{A}^M,\tilde{B}^M,\tilde{C}^M$ upon
the permutation of $\,\tilde{A},\tilde{B},\tilde{C}$
induced by action of $\,g$ on $\,\al_1,\al_2,\al_3$.

Namely, $M$ is determined from the relation $g(\zeta)=\zeta^M$,
where $\zeta$ is a primitive root of unity of order
$\tilde{N}$ that is \,{\em gcd\,} of the orders of $\tilde{T},
\tilde{A},\tilde{B},\tilde{C}$. Thus $M$ is
relatively prime to $\tilde{N}$ and
$\tilde{T}^M, \tilde{A}^M,\tilde{B}^M,\tilde{C}^M$ are
generators of the stabilizers at $P_0,P_1,P_2,P_3$.

Note that applying proper Galois automorphisms in the fibers
gives that
$\tilde{T}'$ is conjugated to a generator
of the stabilizer at $P_0$ and so on for $P_1,P_2,P_3$.
See e.g. Section 2.1 from \cite{FrV}. Which gives the
existence of individual $M$ above, for each
of the ramification points.
We need the {\em branch cycle argument\,} to obtain that
a single $M$ serves all $4$ ramification points and that
it can be uniquely determined modulo $\tilde{N}\,$
from the relation $g(\zeta)=\zeta^M$. This argument uses
a well-known presentation of ramified
covers in the form $z=x^m$ for suitable local parameters.
See Lemma 2.8 from \cite{Vol},
\cite{Bel} and also Section 3.1 from \cite{Det}
and references therein.

Since $M$ is odd, we conclude that
$\tilde{A}',\tilde{B}',\tilde{C}'$ are
conjugated to $\tilde{A},\tilde{B},\tilde{C}$
upon a permutation of $o_1,o_2,o_3$ induced
by $g$. For $\,\tilde{T}'$, we know that it is
conjugated to $\tilde{T}^M$, where
gcd$(M,\hbox{ord}(\tilde{T}))=1$.

Using (\ref{ABCT}),\ \,
ord$(\tilde{T}$\,mod\,$\C^\ast)=$
ord$(-t)=$$2\!N/\hbox{gcd}(2\!N,N-2|k|)$. Therefore
$\,\hbox{ord}(q^{1/2})=2\!N\,$
divides $\hbox{ord}(\tilde{T})$ and
$\hbox{ord}(\tilde{T})=2\!Nl$ for $l\in \N$.
Moreover, $q^{lN/2}t^{lN}$ must be $1$, which gives
that $l=2$ for $k\in \Z$ and $l=1$ for $k\in
\frac{1}{2}\Z\setminus \Z$.
Finally, $\,\ddot{q}\mapsto \ddot{q}^M\,$  can be extended
to $\si\in$Gal$(\Q(\ddot{q})/\Q)$ for this $M$. Using the
{\em branch cycle argument\,}
we conclude that $\si=\si_g$ is the restriction of $g$ to
$\Q(\ddot{q})$, which is claimed in $(ii)$.
\smallskip

Now, on the DAHA end, applying $\si$
to the entries of the matrices representing $\tilde{T},
\tilde{A},\tilde{B},\tilde{C}$, one obtains an irreducible
$\HH$\~module $\tilde{V}^\si$ of the same dimension dim$V$, for
the same $k$ and with $q^{M}$ instead of $q$.
If $\tilde{V}$ is considered over
$\Z_{\ddot{q}}/\mathfrak{p}^m$, then
the module $\tilde{V}^\si$ is naturally defined over
$\Z_{\ddot{q}}/\bigl(\si(\mathfrak{p})^m\bigr)$.

Using the {\em unirigidity\,} in types $(\al,\be)$,
we see that the
collection of matrix images of $\tilde{T}',
\tilde{A}',\tilde{B}',\tilde{C}'$ in
$\tilde{\mathfrak{B}}_{\mathfrak{p},m}$
can be extended to an action of  $\HH\,$ and the resulting module
must be isomorphic to $\tilde{V}^\si$.
Therefore the kernel $U'$ of the homomorphism
\begin{align*}
&\b_1\to \b_1/U'=\hbox{Aut}(F_{\C}\to P^1_{\C})\toeq
\tilde{\mathfrak{B}}_{\mathfrak{p},m}, \hbox{ where } \\
&\ga_0\mapsto \tilde{T}',
\ga_1\mapsto\tilde{A}',\ \ga_2\mapsto\tilde{B}',\
\ga_3\mapsto\tilde{C}',
\end{align*}
coincides with the initial $U$.

Thus the cover $\,f\!:\!F_{\C}\!\to\! P^1_{\C}\,$
and the field $K^\circ(F)$ has a natural
action of the group  $G_K^\circ$.
This gives the claims involving the ideal
$\mathfrak{p}\supset (p)$ from
Parts $(i),$$(ii)$ for types $(\al,\be)$.

In the case
of type $(\breve{\ga})$ (where we have {\em dim-rigidity}\,
but no {\em unirigidity\,}), these
claims hold as well. Indeed,
$\tilde{A}',\tilde{B}',\tilde{C}'$
must be conjugated  in
$\tilde{\mathfrak{B}}_{\mathfrak{p},m}$
to $\tilde{A},\tilde{B},\tilde{C}$ up to a permutation
in this triple.
This excludes
applying $\varsigma_x$, $\varsigma_y$ or their product
for $\tilde{V}$ of type $(\breve{\ga})$;
see (\ref{ABCT}) and (\ref{varsiabc}).

Actually applying $\varsigma_{x,y}$ to $\tilde{V}$ of
type $(\breve{\ga})$
does not change the kernel of the homomorphism from
\,Aut$(f'\!:\! F'\!\to\! P^1)$\, to\,
Aut$(f\!:\!F\!\to\! P^1)$\,,
so the $G_K^\circ$ invariance of $K^\circ(F)$ from Part $(i)$
follows directly from
the {\em dim-rigidity\,} for $(\breve{\ga})$. However this does
not prove that $G_K^\circ$ acts via conjugations for such modules, 
as stated in Part $(ii)$.

Thus we established Parts $(i,ii)$ for
$\tilde{\mathfrak{B}}_{\mathfrak{p},m}$.
This formally results in the corresponding claims there concerning 
the action of complete $G_K$ in the cover
$F^\diamond\to P^1$ corresponding to the group
$\tilde{\mathfrak{B}}_{(p),m}$ from $(i,ii)$.

The justification of
$(iii,iv)$ is similar. We use the uniqueness of
$\lan\cdot\,,\,\cdot\ran$ in $\tilde{V}$ up to proportionality
in Part $(iv)$. Note that $\lan \hat{g}(v)\,,\,\hat{g}(v)\ran$
for $\hat{g}\in \hat{G}_K$ and $v\in \tilde{V}$
is defined over $\Z_{\ddot{q}}/(p^m)$ (not over
$\Z_{\ddot{q}}/\mathfrak{p}^m)$)
and  remains unchanged if $\hat{g}$ is multiplied by any root
of unity.

Generally, it can happen that
$\lan \hat{g}(v)\,,\,\hat{g}(v)\ran=
\ze_g \lan \hat{g}(v)\,,\,\hat{g}(v)\ran$ for
some $\hat{g}\in \hat{G}_K$ and a root of unity $\ze_g$
(serving all $v\in V$). When $p=0$, then
$\ze_g$ can be only $1$ in the case of $(\al)$ due to the positivity
of $\lan\cdot\,,\,\cdot\ran$ at $q=e^{2\pi\imath/N}$.
This is parallel to the proof of unitarity of the projective
action of $PSL_2(\Z)$ in $V$ for types $(\al,\be)$; in the
case of $(\be)$, one uses here that dim$V$ is odd (when $p=0$).

When $p>0$, i.e. for  $\tilde{V}$ considered
over $\Z_{\ddot{q}}/(p^m)$, this arguments can be readily
replaced by considering the exact lists of squares
of $\lan\cdot\,,\,\cdot\ran$ upon the diagonalization
of this form in $\tilde{V}$ (in the basis of $\{e'_n\}$).
\smallskip

The extension from Part $(iv)$ of the previous
claims to the modules
$\,\varsigma_x^{\ep}\varsigma_y^{\de}(V^+_{2n+1})$
is straightforward; use (\ref{ABCT}) and (\ref{varsiabc}).

To adjust these considerations to the cases of $\checkmark$\~modules
form Part $(v)$, one needs to replace $\ddot{q}$ by its
square and use that\,
gcd$\bigl(M,\hbox{ord}(\ddot{q})\bigr)=1$
for  $g(\zeta)=\zeta^M$ (see above). Recall that here
$q,t,\ddot{q}$ must be replaced by
$q^\checkmark=q^2,t^\checkmark=t^2,\ddot{q}^\checkmark=\ddot{q}^2$.
The order of
$\tilde{T}^\checkmark$ divides $2\!N$
($\tilde{T}^\checkmark$ is always semisimple).
This gives the claims from $(i-iii)$ with $\checkmark$.

The same {\em branch cycle argument\,} gives
that $G_K^\circ$ and $G_K$ transpose the
$\varsigma_x^\ep\varsigma_y^\de$\~orbits of the modules
from $(\breve{\al}^\checkmark)$
for even $2k$, $(\breve{\be}^\checkmark)$
and those from $(\breve{\ga}^\checkmark)$
via their action from $(iv)$ for the family $(\breve{\ga})$.
This action is trivial for
$(\breve{\al}^\checkmark)$ for odd $2k$.
\sq
\smallskip

\subsection{\bf Special cases}\label{sec:special-c}
Following essentially \cite{CoH},
Example 2.4, let us calculate the automorphism $\eta$ of
$\tilde{\mathfrak{B}}_{\mathfrak{p},m}$ corresponding to
the {\em complex conjugation} in $K(\ddot{q})$, sending
$\C\ni z\mapsto \overline{z}$.
Here we can consider the whole $\b_1$,
since the complex conjugation is {\em continuous} (the
other automorphisms from Gal$(\overline{\Q}/\Q)$ are
discontinuous).

We will pick the {\em standard homotopy basis\,}.
Namely $\ga_i$ will be represented by
a loop at $b_o$ that winds once around $\al_i$
counterclockwise for $i=1,2,3$ and clockwise for
$i=0$ and winds around no other $\al_j$.
Let us assume that $\al_0<b_o< \al_1<\al_2<\al_3$
are all real. Then $\ga_0=\ga_3\ga_2\ga_1$,
$\ga_1$ and $\ga_2\ga_1$ can be chosen
to have support invariant under complex conjugation.
The composition of the paths $\ga_i$ and then $\ga_j$ is
written $\ga_j\ga_i$. One has
\begin{align}\label{gamma-cc}
&\overline{\ga_0}=\ga_0^{-1},\
\overline{\ga_1}=\ga_1^{-1},\
\overline{\ga_2\ga_1}=\ga_1^{-1}\ga_2^{-1},\\
&\overline{\ga_2}=\ga_1^{-1}\ga_2^{-1}\ga_1,\
\overline{\ga_3}=\ga_1^{-1}\ga_2^{-1}\ga_3^{-1}
\ga_2\ga_1.\notag
\end{align}
Switching to $T',A',B',C'$ from $\b_1$
corresponding to
$\ga_0,\ga_1,\ga_2,\ga_3$ and using that
$T'=A'B'C'=(A')^{-1} (B')^{-1} (C')^{-1}$
(\,$A',B',C'$ are involutions in $\b_1$),
\begin{align}\label{eta-cc}
&\overline{T'}\!=\!(T')^{-1},\
\overline{A'}\!=\!A',\
\overline{B'}\!=\!A'B'A',\
\overline{C'}\!=\!T'B'A'\!=\!T'C'(T')^{-1}.
\end{align}

Now let us substitute
$\tilde{A}=XT,\,
\tilde{B}=q^{1/4}XTY,\, \tilde{C}=T^{-1}Y,\,
\tilde{T}=q^{1/4}T$
for $T',A',B',C'$; see
(\ref{brdnrm}).
To avoid confusion with the complex conjugation
acting in matrices, we will denote the resulting
automorphism of $\b_q$ by $\eta'$; then
$\eta'(q^{1/2})=q^{-1/2}$.

Recalculating now (\ref{eta-cc}) in terms of
$X,T,Y$ we obtain that $\eta'$ acts as follows:
\begin{align}\label{eta-TXY}
&\eta'(T)=T^{-1},\
\eta'(X)=T^{-1}X^{-1}T=XT^2,\
\eta'(Y)=Y^{-1}.
\end{align}
Indeed,  $\eta'(T^{-1}Y)=TCT^{-1}=
T\eta'(Y)=YT^{-1}$, which gives the last
formula. This is the {\em DAHA-Kazhdan-Lusztig involution\,}
$\,\eta$ from Section 2.5.5 of \cite{C101} defined in terms
of $X$ instead of $Y$.
\medskip

{\sf Quotient by the center.}
Aiming at obtaining the action
of $G_K$ instead of smaller $G_K^\circ$, one can try to
consider its quotient of $\tilde{\mathfrak{B}}_{\mathfrak{p},m}$
by its center instead of enlarging
$\tilde{\mathfrak{B}}_{\mathfrak{p},m}$
to $\tilde{\mathfrak{B}}_{(p),m}$ . This is actually
what we are supposed to do for obtaining the {\em Tate module\,}
in the case $t=1$, which was discussed in the Introduction.

\begin{corollary}\label{GALMODCENTER}
In the notations from $(i)$ of the theorem,
let $\tilde{\mathfrak{B}}^\dag_{\mathfrak{p},m}$
be the quotient of $\tilde{\mathfrak{B}}_{\mathfrak{p},m}$
by its center. We assume that either $\tilde{\mathfrak{B}}$
is finite (we set $p=0=m$ in this case) or $p$ is a primitive
element modulo $N$ for $N=\ell^u$ and odd prime $\ell$.
Then $\tilde{\mathfrak{B}}^\dag_{\mathfrak{p},m}$
does not depend on the choice of the prime ideals
$\mathfrak{p}\supset (p)$
up to isomorphisms and the action of $G_{K}^\circ$ in
$\tilde{\mathfrak{B}}_{\mathfrak{p},m}$  and
$\tilde{\mathfrak{B}}_{\mathfrak{p},m}^\checkmark$
provides
a homomorphism
$G_{K}\to
\hbox{Aut}\bigl(\tilde{\mathfrak{B}}^\dag_{\mathfrak{p},m}\bigr)$
for any $\mathfrak{p}\supset (p)$ and $m\in \N$ and its
counterpart
for $\tilde{\mathfrak{B}}_{\mathfrak{p},m}^{\dag\checkmark}$.
If $\,Image(\b_1)$ is finite (and $p=0$), then $G_K$ acts
directly in $\tilde{\mathfrak{B}}^\dag\equal
\tilde{\mathfrak{B}}$/Center$(\tilde{\mathfrak{B}})$.
\end{corollary}
{\em Proof.} The consideration is essentially the same
as in the proof of Theorem \ref{GALACTION}. We
pick $\,g\in G_K$ and consider two homomorphisms
\begin{align*}
&\b_1/U^\dag\,=\,\hbox{Aut}(F^\dag_{\C}\!\to\! P^1_{\C})\toeq
\,\tilde{\mathfrak{B}}^\dag_{\mathfrak{p},m}, \\
&\b_1/U^{\dag\,\prime}=\hbox{Aut}(F^{\dag\,\prime}_{\C}
\!\to\!P^1_{\C})\toeq
\tilde{\mathfrak{B}}^\dag_{\mathfrak{p},m},
\end{align*}
where the (normal) subgroups $U^\dag$, $U^{\dag\,\prime}$
correspond to the cover $F^\dag_{\C}$ associated
with $\tilde{\mathfrak{B}}^\dag_{\mathfrak{p},m}$
and its $\,g$\~image $F^{\dag\,\prime}_{\C}$.

Then we conclude that the conjugation by $\phi_{\hat g}=
h_{\hat g}\si_g$ in GL$(\tilde{V})$ will not
change $U^\dag$. Indeed, the relations in
$\mathfrak{B}_{\mathfrak{p},m}$ will be
transformed by the automorphisms from
Gal$(\Q(q^{1/2})/\Q)$ and therefore do not depend
on the particular choice of $q^{1/2}$ modulo the
center of $\mathfrak{B}_{\mathfrak{p},m}$. Thus the whole $G_K$,
not only $G_K^\circ=$Gal$(\overline{\Q}/K^\circ)$
as in the theorem, will act
in the cover $F^\dag_{\C}\!\to\! P^1_{\C}$. This
gives the required for
$\tilde{\mathfrak{B}}_{\mathfrak{p},m}^\dag$; the passage to
$\tilde{\mathfrak{B}}_{\mathfrak{p},m}^{\dag\checkmark}$
is straightforward.
\sq
\smallskip

{\sf Dimension 4.}
let us consider the
module $V=V'_4$ from
(\ref{v4prime}) as an example.
It is irreducible for a primitive $(2\!N)$th
root of unity $q^{1/2}$ and
remains irreducible when $q^{1/2}\mapsto \check{q}^{1/2}$
for a primitive
root $\check{q}^{1/2}$ of order $2\!\check{N}$. 

Following Section \ref{sec:dim24}, we denote by
$\mathfrak{B}^\dag_N$  and $\mathfrak{B}^\dag_{\check{N}}$,
the quotients of $\,Image(\b_q)\,$ or $\, Image(\b_1)\,$
by the center. Let us also consider 
the commutator subgroup
\begin{align}\label{commsubg}
\c_N^\dag\equal
<\!T^2,X^2,Y^2\!>\!/\!<\!T^{2M},X^{2N},Y^{2N}\!>\,\simeq\,
\Z_{M}\times \Z_N^2,
\end{align}
where $M=N$ for odd $N$ and $M=N/2$ otherwise. 

The quotient of $\mathfrak{B}^\dag_N/\c_N^\dag$  is isomorphic
to $\Z_2^2$ and corresponds to the cover of
$P^1$ that is the composition
$E \simeq E'\to E'/E'_2\simeq E\to P^1$;
the first map is the multiplication by $2$ ($E'_2$ is
the group of points of the $2$nd order). The action
of $\mathfrak{B}^\dag_N$ in $K^\circ(E')$ is via translations by
$E'_2$; $T$ is the reflection $s: z\mapsto -z$ there.
Note that $E'/\{s\}\simeq P^1$ covers $P^1$ with the Galois group
$\Z_2^2$ and $6$ ramification points of
index $2$ over $\al_1,\al_2,\al_3$.

Let us briefly describe the cover $F$ of $E$ corresponding to 
the group $\mathfrak{B}_N^\dag$ and
the cover of $E'$ corresponding to the group $\c_N^\dag$.
One needs to extend the field $K^\circ(E')$ by the functions
$\phi_i=f_i^{1/N}$ for $f_i\in K^\circ(E')$ with the divisors
$(f_i)=2\,0'_i-2\,0'$, where $0',0'_i\, (i=1,2,3)\,$ are 
the points of the second order in $E'$.
Note that $f_1f_2f_3$ is a perfect
square in $K^\circ(E')$. Recall that $K^\circ=K(q^{1/4})$
(actually $q^{1/2}$ is sufficient below).  

We associate the translations $e_1,e_2$ by the points 
$0'_{1},0'_{2}\in E'_2$ respectively 
with $X$ and $Y$. Then the automorphisms
in $Gal(K^\circ(F)/K^\circ(P^1))$
corresponding to the generators $X,Y,T 
\in \mathfrak{B}_N^\dag$ (modulo the center) 
can be presented as follows. For a primitive
$(2\!N)$th root of unity $\om$, $\al,\be\in \C^*$ and
$\al/\be=\om$, one can pick $\phi_{1,2,3}$ to ensure
the relations
\begin{align*}
&X\!\begin{pmatrix}\phi_1\\ \phi_2\\ \phi_3\end{pmatrix}
\!=\!
\begin{pmatrix}\al\om\phi_1^{-1}\\ \frac{\al}{\om}\phi_3 \phi_1^{-1}\\
\phi_2\phi_1^{-1}\end{pmatrix}\!\!,
Y\!\begin{pmatrix}\phi_1\\ \phi_2\\ \phi_3\end{pmatrix}
\!=\!
\begin{pmatrix}\frac{\be}{\om}\phi_3\phi_2^{-1}\\ \be\om\phi_2^{-1}\\
\phi_1\phi_2^{-1}\end{pmatrix}\!\!,
T\!\begin{pmatrix}\phi_1\\ \phi_2\\ \phi_3\end{pmatrix}
\,=\,
\begin{pmatrix}\phi_1\\ \om^2\phi_2\\
\phi_3\end{pmatrix}\!\!,\\
&X^2\!\begin{pmatrix}\phi_1\\ \phi_2\\ \phi_3\end{pmatrix}
\!=\!
\begin{pmatrix}\phi_1\\ \om^{-2}\phi_2\\
\om^{-2}\phi_3\end{pmatrix}\!\!,
Y^2\!\begin{pmatrix}\phi_1\\ \phi_2\\ \phi_3\end{pmatrix}
\!=\!
\begin{pmatrix}\om^{-2}\phi_1\\ \phi_2\\
\om^{-2}\phi_3\end{pmatrix}\!\!,
T^2\!\begin{pmatrix}\phi_1\\ \phi_2\\ \phi_3\end{pmatrix}
\,=
\begin{pmatrix}\phi_1\\ \om^4\phi_2\\
\phi_3\end{pmatrix}\!\!.\\
\end{align*}
Check the relations $(TX)^2=1=(TY^{-1})^2=Y^{-1}X^{-1}YXT^2$ 
using 
\begin{align*}
&X^{-1}\!\begin{pmatrix}\phi_1\\ \phi_2\\ \phi_3\end{pmatrix}
\!=\!
\begin{pmatrix}\al\om\phi_1^{-1}\\ \al\om\phi_3 \phi_1^{-1}\\
\om^2\phi_2\phi_1^{-1}\end{pmatrix}\!\!,\ \,
Y^{-1}\!\begin{pmatrix}\phi_1\\ \phi_2\\ \phi_3\end{pmatrix}
\!=\!
\begin{pmatrix}\be\om\phi_3\phi_2^{-1}\\ \be\om\phi_2^{-1}\\
\om^2\phi_1\phi_2^{-1}\end{pmatrix}\!\!.
\end{align*}

\comment{
\begin{align*}
&X\!\begin{pmatrix}\phi_1\\ \phi_2\\ \phi_3\end{pmatrix}
\!=\!
\begin{pmatrix}\om^2\phi_1^{-1}\\ \phi_3 \phi_1^{-1}\\
\phi_2\phi_1^{-1}\end{pmatrix}\!\!,
Y\!\begin{pmatrix}\phi_1\\ \phi_2\\ \phi_3\end{pmatrix}
\!=\!
\begin{pmatrix}\om^{-1}\phi_3\phi_2^{-1}\\ \om\phi_2^{-1}\\
\phi_1\phi_2^{-1}\end{pmatrix}\!\!,
T\!\begin{pmatrix}\phi_1\\ \phi_2\\ \phi_3\end{pmatrix}\!=\!
\begin{pmatrix}\phi_1\\ \om^2\phi_2\\
\phi_3\end{pmatrix}\!\!,\\
&X^2\!\begin{pmatrix}\phi_1\\ \phi_2\\ \phi_3\end{pmatrix}
\!=\!
\begin{pmatrix}\phi_1\\ \om^{-2}\phi_2\\
\om^{-2}\phi_3\end{pmatrix}\!\!,
Y^2\!\begin{pmatrix}\phi_1\\ \phi_2\\ \phi_3\end{pmatrix}
\!=\!
\begin{pmatrix}\om^{-2}\phi_1\\ \phi_2\\
\om^{-2}\phi_3\end{pmatrix}\!\!,
T^2\!\begin{pmatrix}\phi_1\\ \phi_2\\ \phi_3\end{pmatrix}
\,=\,
\begin{pmatrix}\phi_1\\ \om^4\phi_2\\
\phi_3\end{pmatrix}\!\!.\\
\end{align*}
Check the relations $(TX)^2=1=(TY^{-1})^2=Y^{-1}X^{-1}YXT^2$ 
using 
\begin{align*}
&X^{-1}\!\begin{pmatrix}\phi_1\\ \phi_2\\ \phi_3\end{pmatrix}
\!=\!
\begin{pmatrix}\om^2\phi_1^{-1}\\ \om^2\phi_3 \phi_1^{-1}\\
\om^2\phi_2\phi_1^{-1}\end{pmatrix}\!\!,\ \,
Y^{-1}\!\begin{pmatrix}\phi_1\\ \phi_2\\ \phi_3\end{pmatrix}
\!=\!
\begin{pmatrix}\om\phi_3\phi_2^{-1}\\ \om\phi_2^{-1}\\
\om^2\phi_1\phi_2^{-1}\end{pmatrix}\!\!.
\end{align*}
}
We obtain the required action.
Note that $K^\circ(E')$ extended 
by the $\phi_1 \phi_2/\phi_3$, which is invariant
with respect to $X^2$ and $Y^2$, corresponds
to the homomorphism $\c_N^\dag\to \c_N^\dag/<\!X^2,Y^2\!>\simeq
\Z_{M}$.

In terms of $P^1$ and the images $\,o',\al'_1,\al'_2,\al'_3\,$
of $\,0',0'_1,0'_2,0'_3\,$, one can take 
$f_i=\frac{\al'_i-z}{z}$, 
where $z$ is the coordinate of $P^1$. Then for
$j=1,2$,
$$
f_i^{e_j}\!=\!u_{ij}\frac{f_k}{f_j},
\ u_{ij}\!=\!\frac{\al'_i\!-\!\al'_j}{\al'_i},
\hbox{\, where } 1\!\le\! i,k\!\le\! 3, \
i\!\neq\!j\!\neq\!k,\, k\!\neq\!i,
$$  
and $f_j^{e_j}=u_{ij}u_{kj} f_j^{-1}$. Setting 
$\al_1'=\al_3'-\al_2'$, we arrive at the relations
above for the action of $X,Y$ with $\om=-1$ and
$\al\!=\!-\be\!=\!\frac{\al'_1-\al'_2}{\al'_3}$.
\smallskip

\subsection{\bf Triangle groups}
As it has been already mentioned, the root system $C^\vee C_1$ 
is a natural setting for Theorem \ref{GALACTION} and its proof
(based on the Riemann Existence Theorem for $\C P^1$ with $4$
ramification points).
However $\b_1$ of type $A_1$ and its quotients can be associated
with the covers of elliptic curves, which is of obvious
importance. Moreover, the case of $A_1$ has the following
interesting relation to the {\em triangle groups\,}.
\smallskip

Given $l,m,n \in \{2,3,4,\ldots\}$, the corresponding
triangular group is
\begin{align}\label{triang}
\De(l,m,n)\!\equal
<\!a,b,c\,\mid\, a^2\!=\!b^2\!=\!c^2\!=\!1\!=\!(ab)^l\!=\!
(bc)^m\!=\!(ca)^n\!>\!.
\end{align}
It contains a normal subgroup of index $2$,
the {\em von Dyck group\,},
\begin{align}\label{triangprime}
&D(l,m,n)\!=
<\!x,y,z\,\mid\, x^l=y^m=z^n=1=xyz\!>, \hbox{ where }\\
&x=ab,\, y=bc,\, z=ca,\ \,
axa=x^{-1},\, aya=xz,\, aza=z^{-1}.\notag
\end{align}
It is commonly called a triangle group, as well as
$\De(l,m,n)$,
especially in the papers devoted to simple quotients
of triangle groups.

The groups $\De(l,m,n;r),\,D(l,m,n;r)$ are defined by
imposing the additional relation $(abac)^r=1=(xz^{-1})^r$.
This definition is related to Conjecture 5.3 from \cite{Sch}
concerning the discreteness of the images of $\De(l,m,n)$
in $PU(2,1)$. The case when $l=m=n$ is called {\em equilateral}.

\begin{proposition}\label{TRIANGB}
(i) For the modules $V=V_{2N-4k},V^{\ep,\de}_{2|k|}$
over $\Z_{\ddot q}$ of types $(\al,\breve{\ga})$ from
Theorem \ref{psl-two-inv}
and those obtained upon applying $\iota$, the group
$\tilde{\mathfrak{B}}^\dag$,
the image $\tilde{\mathfrak{B}}$
of $\b_1$ in automorphisms of $V$ divided by its
center, is a natural
quotient of $\De(2\!N,2\!N,2\!N;2\!N)$\,:
\begin{align}\label{abctotrig}
&a\mapsto \tilde{A}\!=\!XT,\,
b\mapsto \tilde{B}\!=\!q^{1/4}XTY,\,
c\mapsto \tilde{C}\!=\!T^{-1}Y,\
x\mapsto q^{1/4}Y,\\
y\mapsto &q^{1/4}T^{-1}X^{-1}T\!=\!q^{1/4}XT^2,\
z\mapsto Y^{-1}X^{-1}\!=\!q^{1/4}Y^{-1}\tau_-(X^{-1})Y. \notag
\end{align}
Namely, the orders of the images of $x,y,z$ are $2\!N$, as well as 
the order of \,Image$(abac)=YXY=\tau_+\tau_-(\tilde{Y})$.

Assuming the semisimplicity of $T$ (i.e. for $t\neq -1$),
the order of
$\,Image(abc)=\tilde{T}=q^{1/4}T=
\tilde{A}\,\tilde{B}\,\tilde{C}$ modulo
scalar matrices is \,ord$(-t)\,=\,$
$2\!N/\hbox{gcd}(2\!N,N-2|k|)$. Also, the
order of $\,Image(xz^{-1}y)=\tilde{T}^2\,$
in \,GL$(V)/\C^*\,$ equals \,ord$(t^2)=
N/\hbox{gcd}(N,2|k|)$.

(ii) In the case of odd $N$ for $(\breve{\al}^\checkmark)$
when $2k$ is even and for $(\breve{\ga}^\checkmark)$,
we substitute $q\mapsto q^\checkmark=q^2$,
$t\mapsto t^\checkmark=t^2$ in the formulas
from (\ref{abctotrig}):
\begin{align}\label{abctotrigcheck}
&a\mapsto XT^\checkmark,\
b\mapsto q^{1/2}XT^\checkmark Y^\checkmark,\
c\mapsto (T^\checkmark)^{-1}Y^\checkmark,\\
&x\mapsto q^{1/2}Y^\checkmark,\,
y\mapsto q^{1/2}(T^\checkmark)^{-1}X^{-1}T^\checkmark,\,
z\mapsto Y^{-1}X^{-1}. \notag
\end{align}
Then the orders of the images of $x,y,z$ modulo scalars
are $N$, as well as the order of \,Image$(abac)=
\tau^\checkmark_+\tau_-^\checkmark(\tilde{Y}^\checkmark)$.
Accordingly, the images of $a,b,c,x,y,z$ satisfy the
relations from $\De(N,N,N;N)$. The dimensions of the
corresponding spaces $V^\checkmark$ are
$N-2k$ for $(\breve{\al}^\checkmark)$
and $2|k|$ for $(\breve{\ga}^\checkmark)$.

Since $t^\checkmark\neq -1$ ($N$ must be odd when
$\,\checkmark$ is used),
$\,Image(abc)=q^{1/2}T^\checkmark$ modulo
scalar matrices is of order\, ord$(-t^2)\,=\,$
$2\!N/\hbox{gcd}(N,2|k|)$ and the order of
$(T^\checkmark)^2$ is $N/\hbox{gcd}(N,2|k|)$.
\end{proposition}

{\em Proof.} Calculating the orders of $\tilde{T}$ and
$\tilde{T}^2$ is straightforward. Since we
need them modulo scalar matrices, the coefficients
of proportionality do not matter here.
Concerning $\,Image(x)=q^{1/4}Y$, its spectrum
in modules of types $(\al,\ga)$
is described in (\ref{nonsymp}); see also
Theorem \ref{tnegk},$(ii)$ and
Theorem \ref{perfectposit},$(i)$. This gives\,
ord$(Y\hbox{mod }\C^*)$.
Applying conjugations and the automorphisms
$\tau_{\pm}$, which act in $V_{2N-4k}$ and
in $V_{2|k|}^+$, we see that the images
of $y,z,abac=xz^{-1}$ will have the same orders modulo $\C^*$
as $Y$. So the corresponding triangular groups are
{\em equilateral}.
The modules from $(\breve{\ga})\setminus (\ga)$ are
obtained from $V_{2|k|}^+$  by changing the signs
of $X,Y$, which does not influence their orders
modulo scalars.
The cases with $\checkmark$ are straightforward.
\sq
\smallskip

Note that we excluded in the theorem the modules
of type $(\be),(\be^\checkmark)$ and ($\al^\checkmark)$
for odd $2k$ due to the presence of the Jordan blocks
in $X,Y$. See Theorem 2.9.3 from
Section 2.9.2 of \cite{C101}.
The images of the elements
$x,y,z,xz$  become of infinite order in these cases.
Their orders in $\tilde{\mathfrak{B}}^\dag_{\mathfrak{p},m}$
and $\tilde{\mathfrak{B}}^{\dag\checkmark}_{\mathfrak{p},m}$
are
$2\!N p^m$  for $(\be)$ and
$N p^m$ for $(\al^\checkmark,\be^\checkmark)$.

\subsection{\bf The Livn\'e groups}
Interestingly, the DAHA-modules of dimension $3$ already
provide rich theory, matching the so-called
{\em Livn\'e groups}.

The images of the
triangular group in $PU(2,1)$  are extensively
studied algebraically and using the methods of complex
hyperbolic geometry. Generally, the images of the
product $\,abc\,$
(our $\tilde{T}$) do not satisfy quadratic equations.
For instance the element $\,abc\,$
in the well-studied $\De(4,4,4;5)$
has $3$ distinct eigenvalues in $PU(2,1)$; see \cite{Der}.

Our discrete groups for dim$V=3$ appeared exactly the
canonical subgroups of the Livn\'e lattices (examples of
the Mostow groups). Conjecture \ref{DISCCONJ} provides
examples of arithmetic discrete groups \,Image$(\b_1)$
in higher dimensions.


Let us demonstrate what our construction gives
for\, dim$V=3$, when it results in
subgroups of $PU(2,1)$.
It occurs only for the series $(\breve{\ga})$ and
$(\breve{\ga}^\checkmark)$,
so the group \,Image$(\b_1)$ is
never finite when \,dim$V=3$, except for the boundary
case $(\ga^\bullet)$ from (\ref{boundcasea}).

Let us consider $(\ga)$, including $(\ga^\bullet)$.
Then $k=-3/2,\,
t=q^{-3/2}$, $N=3,4,\ldots\ $. Using (\ref{epmtex}),
$$
V^+_3\!=\!\mathscr{X}/(\ep_{2}\!-\!\ep_{-1})\!=\!
\Q[q^{1/4}][X^{\pm 1}]/\bigl(
(X\!-\!q^{1/4})(X\!-\!q^{-1/4})(X\!-\!q^{3/4})\bigr).
$$
Also,
$V^+_3=\oplus_{i=1}^3 \Q[q^{1/4}](\ep'_i+\ep'_{i-3})\subset
W_{6}=\!\mathscr{X}/(\ep_{3})$
as DAHA-modules.
Here $\ep'_m$ are the images in $V^+_3$ of \
$\ep_0=1,\, \ep_1=t^{1/2}X,$
\begin{align}
\ep_{-1}& =
t^{1/2}(\frac{1-tq}{1-t^2 q}X^{-1}+
\frac{1-t}{1-t^2 q}X),\ \,
\ep_2 = t(\frac{1-tq}{1-t^2 q}X^{2}+
\frac{q-qt}{1-t^2 q}), \notag\\
\ep_{-2}&\,= \ \frac{t(1-tq)(1-tq^2)}{(1-t^2 q)(1-t^2 q^2)}\
 \bigl(X^{-2}+\ \
\frac{1-t}{1-tq^2}X^2+\ \frac{(q+1)(1-t)}{1-tq^2}\bigr),\notag\\
\ep_{3}\ &= \frac{t^{3/2}(1-tq)(1-tq^2)}{(1-t^2 q)(1-t^2 q^2)}
\bigl(X^3+
q^2\frac{1-t}{1-tq^2}X^{-1}+ q\frac{(q+1)(1-t)}{1-tq^2}X\bigr).
\notag
\end{align}

Since we calculate modulo $(\ep_2-\ep_{1})$,
it suffices to take here $\ep_i'$ for $i=1,2,0$; the exact
sums $\ep'_i+\ep'_{i-3}$
are needed if we realize $V_{3}^+$
inside $W_6$, which is $\mathscr{X}/(\ep_3)$.
Moreover, one can easily
recalculate the formulas to the basis $\{e_1',e_2',e_0'\}$
in $V_3^+$\,;\, $Y$ will remain unchanged and the matrix for $T$
from (\ref{matT-tot}) will be conjugated
by \,diag$(e_m(t^{-\frac{1}{2}})\,,\, m=1,2,0)$.

For $m=1,2,3$, one has:
\begin{align}\label{Yspectrum}
Y(\ep_m\!-\!\ep_{m-3})/(\ep_m\!-\!\ep_{m-3})=
Y(e_m)/e_m=t^{-\frac{1}{2}}q^{-\frac{m}{2}}=
t^{\frac{1}{2}}q^{\frac{3-m}{2}}.
\end{align}
The spectrum of $X$ in $V_3$ is the same as that of
$Y^{-1}$. The spectrum of $\tilde{Y}=q^{1/4}Y$
is $\{1,q^{1/2},q^{-1/2}\}$.

The action of $T$ on the vectors $\{\ep_m+\ep_{m-3},
\ep_{3-m}+\ep_{-m}\}$ for $0<m<3$ is the same as its
action on $\{\ep_m,\ep_{-m}\}$. It is
given by the matrix
\begin{align}\label{matTv3}
T_{\{\ep_m,\ep_{-m}\}}\ =\
\begin{pmatrix}
\frac{t^{1/2}-t^{-1/2}}{1-tq^m}
   & t^{\frac{1}{2}}\frac{1-q^m}{1-tq^m}
\\ t^{-\frac{1}{2}}\frac{1-t^2q^m}{1-tq^m}
   & \frac{t^{1/2}-t^{-1/2}}{1-t^{-1}q^{-m}}
\end{pmatrix}.
\end{align}
This presentation formally works even for $m=0,3$;
then $T$ acts as $t^{1/2}$, which is the
sum of the entries in either column (the
basic vectors then become proportional for such $m$).

Thus in the basis $\{\ep_m+\ep_{3-m}\,,\,m=1,2,3\}$,
\begin{align}\label{matT-tot}
T\ =\
\begin{pmatrix}
\frac{t^{1/2}-t^{-1/2}}{1-tq}
   & t^{\frac{1}{2}}\frac{1-q}{1-tq} & 0
\\ t^{-\frac{1}{2}}\frac{1-t^2q}{1-tq}
   & \frac{t^{1/2}-t^{-1/2}}{1-t^{-1}q^{-1}} &0
\\ 0 & 0 & t^{1/2}
\end{pmatrix}.
\end{align}
Its eigenvalues are $\{t^{1/2},t^{1/2},-t^{-1/2}\};$
the eigenvalues of $\tilde{T}=q^{1/4}T$ are
$\{q^{-1/2},q^{-1/2},-q\}$. The matrix $T$ has
coinciding eigenvalues and a Jordan $2$\~block when
$N=3$ at the eigenvalue $t^{1/2}=q^{-3/4}=\imath$.

Since $Y^{-1}\,=\,$diag\,$(\,t^{1/2}q^{1/2}=q^{-1/4},\,
t^{1/2}q=q^{1/4},\, t^{1/2}q^{3/2}=q^{3/4}\,)$ in the basis
above, the matrix for $\tilde{C}=T^{-1}Y=Y^{-1}T$ equals
\begin{align}\label{matC-tot}
\tilde{C}\!=\!
\begin{pmatrix}
q^{\frac{1}{2}}\frac{t-1}{1-tq}
   & q^{\frac{1}{2}}t\frac{1-q}{1-tq} & 0
\\ q\frac{1-t^2q}{1-tq}
   & q\frac{t-1}{1-t^{-1}q^{-1}} &0
\\ 0 & 0 & q^{\frac{3}{2}}t
\end{pmatrix}\!=\!
\begin{pmatrix}
q^{\frac{1}{2}}\frac{q^{-\frac{3}{2}}-1}{1-q^{-\frac{1}{2}}}
   & q^{-1}\frac{1-q}{1-q^{-\frac{1}{2}}} & 0
\\ q\frac{1-q^{-2}}{1-q^{-\frac{1}{2}}}
   & q\frac{q^{-\frac{3}{2}}-1}{1-q^{\frac{1}{2}}} &0
\\ 0 & 0 & 1
\end{pmatrix}.
\end{align}

Following \cite{Par1} and the notations
from this paper (see formula (2.4) and Theorem 1.2 there), let
\begin{align}\label{matJI}
J\,=\,
\begin{pmatrix}
0 & 0 & 1
\\ 1  & 0 & 0
\\ 0 & 1 & 0
\end{pmatrix},\ \,
I_1\,=\,
\begin{pmatrix}
1 & \tau & \overline{\tau}
\\ 0 & -1 & 0
\\ 0 & 0 & -1
\end{pmatrix}.
\end{align}
We set $I_2=JI_1J^2,\, I_2= J^2 I_1 J,\,
I_{123}=I_1I_2I_3=(JI_1)^3$; use that $J^3=1$.
Here $\tau\in \C^*$. If $|\tau|=2\cos(\pi/l)=
2\cos(\phi/2)$ for $\phi\equal 2\pi/l$
and $|\tau^2-\overline{\tau}|=2\cos(\pi/m)$,
then $\De[l;m]\equal \,<\!I_1,I_2,I_3\!>$ is
a quotient of $\De(l,l,l;m)$ for $3\le l,m\in \N$.
We send $a\mapsto I_1,b\mapsto I_2,c\mapsto I_3$.

The groups $\De[l,l,l]$ for
$$
\tau=\tau_l\equal e^{2\imath\phi/3}+e^{-\imath\phi/3}
=2e^{\imath\phi/6}\cos(\phi/2)
\for \phi\equal 2\pi/l\, (l\ge 3)
$$
are natural subgroups of the {\em Livn\'e groups\,},
which are $\,<\!J,I_1\!>.$
In this case,
$|\tau^2-\overline{\tau}|=|\tau|$ and $m=l$, so
they are of the type $\De[l;l]$. See 
(\ref{uvw-cond}) and around for further discussion.
The eigenvalues of $I_1, I_{12}=I_1I_2$\, (coinciding
with those for $I_{23}=I_2I_3,I_{31}=I_3I_1$)\, and
the eigenvalues of $I_{123}$
can be found following entry $(ii)$ in the
table after Theorem 3.1 in \cite{Par1}.
Namely, they are 
\begin{align}\label{eigenYT}
\{-1,-1,1\}_1,\
\{1,e^{\imath \phi},e^{-\imath\phi}\}_{12},\
\{e^{2\imath \phi},-e^{-\imath\phi},-e^{-\imath\phi}\}_{123}.
\end{align}
Here $I_{123}$ is parabolic at $l=6$ and semisimple otherwise. 

\begin{proposition}\label{LIVNE}
Let $q^{1/2}=e^{\imath\phi}$ for
$\phi= 2\pi/(2\!N)=\pi/N$.
\smallskip

(i) The image of $\b_1$ in
GL$(V^+_{3})$ divided by the center is conjugate to
$\De[2\!N;2\!N]$. More exactly, conjugation by a proper matrix
sends
\begin{align}\label{conjlivn}
&\tilde{A}\mapsto -I_1,\ \,
\tilde{B}\mapsto -I_2,\ \,\tilde{C}\mapsto -I_3,\ \,
\tilde{T}=q^{1/4}T\mapsto -I_{123}, \\
&\tilde{Y}=q^{1/4}Y\mapsto I_{12},\ \
\tilde{X}\tilde{T}^2=q^{1/4}XT^2\mapsto I_{23},
\ Y^{-1}X^{-1}\mapsto
I_{31}.\notag
\end{align}

Confirming Conjecture \ref{DISCCONJ} above for $V_{2|k|}^+$ of
type $(\ga)$, the sequence $2\!N=6,8,10,12,18$
from (\ref{gadiscr}) for dim$V=3$ coincides
with the even part of the list of Proposition 4.8
from \cite{Par1}. These are the only
values of $2\!N$ when the image of $\b_1$ in
GL$(V^+_{3})$ is discrete non-finite; all images
are arithmetic for such $2\!N$.
\smallskip 

(ii) Additionally, let us substitute
$q\!\mapsto\! q^\checkmark\!\!=\!q^2,
q^{1/2}\!\mapsto\! (q^{1/2})^\checkmark\!\!=\!q$ for odd $N\ge 3$.
Then the action of $\HH^\checkmark=
\HH_{q,q^{-3/2}}$ in $\hat{V}^{+\checkmark}_3$
is given by the same formulas (\ref{Yspectrum}) and (\ref{matT-tot})
upon the substitution above and this module is (remains)
irreducible due to Proposition \ref{ODDN}.
Then the conjugation equivalence from (\ref{conjlivn})
holds and the image of $\b_1$ in  GL$(V^{+\checkmark}_{3})$
modulo the center  is conjugate to $\De[N;N]$, which 
corresponds to $\phi^\checkmark\!\!=\!2\pi/N$.

Accordingly, the odd cases from
\cite{Par1}, Proposition 4.8 give that
the image of $\b_1$ in GL$(V^{+\checkmark}_{3})$
is discrete and non-finite if and only if $N=5,7,9$.
The images for $N=5,7$ are arithmetic, which matches
Conjecture \ref{DISCCONJ}; the image of $\b_1$ for $N=9$ is
discrete and non-arithmetic.
\end{proposition}

{\em Proof.} We only need to establish the equivalence
of the theory of images of triangular groups
in $PU(2,1)$ and that concerning
the images of $\b_1$ in $3$-dimensional irreducible
modules of DAHA of type $A_1$.  It suffices to observe
that $I_{123}$ satisfies the quadratic equation in
the first theory, which uses (\ref{eigenYT})
and semisimplicity of $I_{123}$ for generic real $\phi$.
Then we use the classification
of $3$-dimensional modules of $\HH\,$ and its
variant for $\HH^\checkmark\,$ from the second theory.
\sq

\smallskip
We note that our analysis of the discreteness of the
images of $\b_1$
for the lists in (\ref{gadiscr}), including
the case dim$V=3$, is similar to
the approach from \cite{Par1} and other papers devoted
to triangle groups in $PU(2,1)$. We analyze the
signature of the DAHA inner product, which is proportional
to that in \cite{Par1,Par2} in the Livn\'e case. Though
our Conjecture \ref{DISCCONJ} is {\em in any dimensions\,};
it is for the images of
$\De(2\!N,2\!N,2\!N; 2\!N)$ and $\De(N,N,N; N)$
subject to the quadratic relation for
$\tilde{T}\sim abc$ in irreducible rigid modules
of $\HH$ and $\HH^\checkmark$.
\smallskip

The {\em parabolic case\,} is interesting, which is for
$V^+_{2|k|}$ with $2\!N=6$, $k=-3/2$
and $q^{1/2}=e^{\pi\imath/3}$.
Recall that here $t^{1/2}=q^{-3/4}=\imath\,$
and $T$ has one Jordan $2$\~blocks at the eigenvalue
$t^{1/2}$;
this is the boundary case $(\ga^\bullet)$ in
(\ref{boundcasea}).
 Then the Livn\'e group $\,<\!J,I_1\!>\,$
is conjugate to the Eisenstein-Picard modular
group, which is  $PU\bigl(2,1;\, \Z[q^{1/2}]\bigr)$
supplied with the natural inner product.
See \cite{Par2} and references therein.

Here $\tilde{T}\sim I_{123}=(I_1J)^3$.
Assuming that \, $p\neq2,3$, we can calculate in this case
the image of $\mathfrak{B}^\dag_{\mathfrak{p},m}$
of $\b_q$ in automorphisms of $V_3^+$  (modulo scalars)
considered over
$\Z[q^{1/2}]/\mathfrak{p}^m$. Multiplying $T,X,Y$ by $t^{1/2}$,
we make $\Z[q^{1/2}]$
the ring of definition of $V_3^+$ (instead of $\Z[q^{1/4}]$).
Let $\mathfrak{p}\supset (p)$ be a prime
ideal in $\Z[q^{1/2}]$.
We obtain that for $q^{1/2}=e^{\pi\imath/3}$,
\begin{align}\label{livne6}
\mathfrak{B}^\dag_{\mathfrak{p},m}\ =\
PU\bigr(2,1;\, \Z[q^{1/2}]/
\mathfrak{p}^m\,\bigl).
\end{align}

Thus Theorem \ref{GALACTION} and
Corollary \ref{GALMODCENTER} supply
$\mathfrak{B}^\dag_{\mathfrak{p},m}$
and $\mathfrak{B}^\dag_{(p),m}$
with an action of $G_K^\circ$ and $G_K$,
where $K$ is the field of definition
of $\al_0\in P^1\supset \{\al_1,\al_2,\al_3\}$
and the corresponding elliptic curve $E\ni 0$ in this
theorem. Recall that $\{\al_1,\al_2,\al_3\}$ are not
assumed to be
individually defined over $K$ (only as a set). The group
$G_K$ permutes
the corresponding generators $\tilde{A},\tilde{B},\tilde{C}$
accordingly; see Part $(iv)$ of Theorem \ref{GALACTION}.

We would like to
mention here that Livn\'e in his thesis considered his groups
in connection with a certain branched cover
of degree $2$ of the universal elliptic curve.

\subsection{\bf Some perspectives}
One can try to
obtain more general groups $\De[l,m,n]$, the
images of triangle groups $\De(l,m,n)$ in $PU(2,1)$ with the
DAHA-type quadratic relation for $I_{123}$. Let us recall the
basic construction of the group $\De[l,m,n]$ (see e.g. \cite{Sch}).
For $u,v,w\in \C^*$, we begin with unitary complex reflection
\begin{align*}
I_1\,=\,
\begin{pmatrix}
1 & u & \overline{v}
\\ 0 & -1 & 0
\\ 0 & 0 & -1
\end{pmatrix},
\ \,I_2\,=\,
\begin{pmatrix}
-1 & 0 & 0
\\ \overline{u} & 1 & w
\\ 0 & 0 & -1
\end{pmatrix},
\ \,I_3\,=\,
\begin{pmatrix}
-1 & 0 & 0
\\ 0 & -1 & 0
\\ v & \overline{w} & 1
\end{pmatrix}.
\end{align*}
with respect to the Hermitian form 
$$
(e_1,e_1)\!=\!(e_2,e_2)\!=\!(e_3,e_3)=2, \
(e_1,e_2)\!=\!\overline{u},\,
(e_2,e_3)\!=\!\overline{v},\, (e_3,e_1)\!=\!\overline{w}
$$
of determinant $det=2(\Re(uvw)\!-\!|u|^2\!-\!|v|^2\!-\!|w|^2\!+\!4)$
for the basic vectors $\{e_i\}$. Namely, $I_i(z)=-z+(z,e_i)e_i\,$
for $z\in \C^3$.

The eigenfunctions of matrices $I_{12},I_{23}$ and $I_{31}$
are respectively 
$$\{1,e^{\pm\imath\al}\},\
\{1,e^{\pm\imath\be}\},\ \{1,e^{\pm\imath\ga}\}
\for 0\le \al,\be,\ga \le \pi
$$ 
if and only if 
$|u|\!=\!2\cos(\al/2)$,\, $|v|\!=\!2\cos(\be/2),\,
|w|\!=\!2\cos(\ga/2)$. Accordingly, 
the matrices $I_{12},I_{23},I_{31}$ are of orders
$l,m,n$ exactly when
$\,\al\!=\!2\nu_1\pi/l$,\, $\be\!=\!2\nu_2\pi/m$,\,
 $\ga\!=\!2\nu_3\pi/n$, where $\nu_i\in \N$,\,
$\nu_1 <l, \nu_2<m, \nu_3<n$, \,gcd$(\nu_1,l)=1$ and 
so on for $\nu_{2,3}$. 
\smallskip

Setting $q=e^{2\imath\phi}$ for $0\le \phi\le \pi$,
the condition
\begin{align}\label{uvw-cond}
3+uvw-(|u|^2+|v|^2+|w|^2)=e^{2\imath\phi}-2e^{-\imath\phi}
\end{align}
is necessary and sufficient
for $I_{123}$ to
have the following eigenvalues:
$$
\{-q^{1/4}t^{1/2},-q^{1/4}t^{1/2},q^{1/4}t^{-1/2}\}\!=\!
\{-q^{-1/2},-q^{-1/2},q\}
\hbox{\, for } t^{1/2}\!=\!q^{-3/4}.
$$
This condition is a direct calculation of
$Tr(I_{123})$; the corresponding relation
for $Tr(I_{123}^2)$ follows from (\ref{uvw-cond})
automatically due to $\phi\in\R$. Note that
$det<0$ and the Hermitian form is of signature $(2,1)$
for $0<\phi<\pi/2$.

Equation (\ref{uvw-cond}) reduces
to that with $u,v\in \R_+^*$ via the  
substitutions
$u\mapsto \ep_1u, v\mapsto \ep_2u, w\mapsto \ep_3w$
for unimodular $\ep_i$ such that $\ep_1\ep_2\ep_3=1$.
Then $\Im w=(\sin(2\phi)+2\sin(\phi))/(uv)$ and one 
arrives at a quadratic equation for $\Re w$,
which always has solutions in $\R$ for sufficiently large $u,v$;
the corresponding $I_{123}$ will be generally 
non-semisimple.

Let us impose (\ref{uvw-cond}) for generic real $\phi$
and assume furthermore that $I_{123}$ is semisimple. 
Then $T\equal -q^{-1/4}I_{123}$ will satisfy the DAHA 
quadratic relation for $t^{1/2}=q^{-3/4}$ and we can
employ Part $(ii)$ of Theorem \ref{tnegk}. We conclude that 
the set of matrices $\{-I_1,-I_2,-I_3,-I_{123}\}$
is conjugated to $\{\tilde{A},\tilde{B},\tilde{C},\tilde{T}\}$
in $V_{2|k|}^+$ for $k=-3/2$.

In terms of (\ref{uvw-cond}), the  semisimplicity of
$I_{123}$ for sufficiently general real $\phi$ 
implies that  $u=\ep_1\tau, v=\ep_2\tau, w=\ep_3\tau$
for $\tau= e^{2\imath\phi/3}+e^{-\imath\phi/3}$
and unimodular $\ep_i$ such that $\ep_1\ep_2\ep_3=1$.
We do not see how to derive this fact directly from
(\ref{uvw-cond}) and the semisimplicity of $I_{123}$.
Recall that the coincidence of $\{u,v,w\}$ with $\{\tau,\tau,\tau\}$
up to $\ep_i$ reflects the
projective $PSL_2(\Z)$\~invariance of $V_{2|k|}^+$.


Now let $\phi=\nu\pi/N$ for $N\ge 3$ and $(\nu,N)=1$.
Using Theorem 2.9.2 of \cite{C101}, we obtain
that the irreducible $\HH$\~modules of dimension $3$
for $q=e^{2\imath\phi}$ can be only
$V_3^+$ or $V_3^{+\checkmark}$
up to a possible change of signs of $T,X,Y$. This results in
$l=m=n$ and the group $\De[2\!N,2\!N,2\!N;2\!N]$ or
that with $2\!N\mapsto N$
for $\checkmark$. We conclude that $3$-dimensional
representations of DAHA of type $A_1$ for 
$q=e^{2\imath\nu\pi/N}$ with $N\ge 3$ and $\nu$ as above
generate
only (subgroups of) the Livn\'e groups;
see Proposition \ref{LIVNE}.
\smallskip

{\sf Non-rigid theory.}
The connection with triangle groups from
Proposition \ref{TRIANGB} links our work to
\cite{LLM} and quite a few papers
devoted to finite/simple quotients of triangular groups.
For instance, Corollary 1.2 from \cite{LLM}
establishes that there are infinitely many
prime numbers $p$ (their density is nonzero) such
that linear and orthogonal groups  over $\Z/(p^m)$
can be presented as quotients of $D(2\!N,2\!N,2\!N;2\!N)$,
or $D(N,N,N;N)$, which we can reach via rigid DAHA-modules
for $A_1$.

Given $V$ of
types  $(\al,\breve{\ga})$ with the standard
DAHA inner product, finding $\,p\,$ such that
$\,\tilde{\mathfrak{B}}^\dag_{\mathfrak{p},m}\,=\,$
$PU(V)\,$ for $V$ considered over
$\Z_{\ddot{q}}/\mathfrak{p}^m$ is of obvious
interest; cf. (\ref{livne6}).
Note that we need unitary groups here (not considered
in \cite{LLM}); also our
$\tilde{\mathfrak{B}}^\dag_{\mathfrak{p},m}$ are
of special kind due to the quadratic relation
for $\tilde{T}$ (which is the key in DAHA theory).

Using all, not only rigid, DAHA-modules from
Theorem 2.9.2 of \cite{C101} is quite
natural for establishing a connection with
\cite{LLM}, which is essentially based on the deformation
theory of the character varieties for triangle groups,
i.e. {\em non-rigid\,} representations. Even in type $A_1$, 
one can obtain here interesting examples of quotients of 
various triangular groups.
The action of $G_K$
can be generally defined for non-rigid representations;
it will then transform the parameters of DAHA-modules. 

Finding the transformation of these parameters under the
projective action of $PSL_2(\Z)$ (a similar and 
simpler problem) is important in the theory of 
the {\em Lam\'e and Heun\,} equations. The latter equation
corresponds to the root system $C^\vee C_1$ and is directly
related to the Painlev\'e VI equation.
These equations require DAHA at the critical center charge
$q=1$ (the $t$\~parameters can be 
arbitrary). The unitarity of the corresponding monodromy 
representations (for proper $t$ and the spectral
parameter), which are DAHA modules 
of dimension $2$ at $q=1$, and the existence of the 
action of the Fourier transform $\sigma$ there are among
the core questions in the theory of these equations. 
\smallskip

Let us mention here the $3\times 3$ Fuchsian representation
of the Painlev\'e VI equation from \cite{Boa} and related
topics considered there. This representation is connected with
our using $3$-dimensional DAHA modules, but 
this is beyond the present
paper (we study only rigid modules). Also, the monodromy
nature of DAHA modules combined with
Proposition \ref{LIVNE} provide an interpretation
of the Livn\'e groups in terms of the
monodromy of the {\em KZB} and similar local systems.

\smallskip


{\sf Compatibility at different roots of unity.}
A challenging problem concerning Theorem \ref{GALACTION}
is in establishing compatibility of
groups $\tilde{\mathfrak{B}}^\dag_{\mathfrak{p},m}$
for different roots of unity $q$. 
By contrast, given $q$ and prime $p$, the compatibility with
respect to $m$ is straightforward. One can readily  define
the action of $G_K^\circ$
and $G_K$ in
$$
\tilde{\mathfrak{B}}_{\mathfrak{p},\infty}=\varprojlim_{\,m}
\tilde{\mathfrak{B}}_{\mathfrak{p},m},\
\tilde{\mathfrak{B}}_{\mathfrak{p},\infty}^\dag=\varprojlim_{\,m}
\tilde{\mathfrak{B}}_{\mathfrak{p},m}^\dag,\
\tilde{\mathfrak{B}}_{(p),\infty}=\varprojlim_{\,m}
\tilde{\mathfrak{B}}_{(p),m}.
$$
Accordingly, the conjugation matrices
$h_{\hat{g}}$ from Theorem \ref{GALACTION}
will have entries in the $p$-adic rings
$\varprojlim_{\,m}
\Z_{\ddot{q}}/\mathfrak{p}^m$
and $\varprojlim_{\,m}
\Z_{\ddot{q}}/(p^m).$ Recall that the $p$\~adic limit is not 
needed if $\tilde{\mathfrak{B}}$ is already finite.
\smallskip

The dependence
of $\tilde{\mathfrak{B}}^\dag_{\mathfrak{p},m}$
on $N\,$ is very far from
straightforward even in case of the triangle groups.
The deformation construction from Section \ref{sec:DEFORMV}
suggests the following (general) approach, which seems
important in its own right. Let us focus on
 $V_{2N-4k}$ of type $(\al)$ from Proposition
\ref{ext-verlinde}; the module  $V'_{4|k'|}$ for
$k'=k-N/2$ is its deformation,
well defined for any nonzero $q$. Upon the substitution
$q\mapsto \check{q}\equal q^d$ for $d|N$, one obtains the action of
$\check{\HH}\equal \HH_{q^{d/2},t^{d/2}}$ in  $V'_{4|k'|}$.
The decomposition
of this module is interesting  in its own right.

This construction is connected with the passage $q\mapsto
q^\checkmark=q^2$ for odd $N$, which is used to address the
case of odd $N$ when $q^{1/2}$ is not a primitive root of unity.
However, $q\mapsto \check{q}$
cannot be associated with any Galois automorphisms
due to $d|N$ (unless $N=\ell^u$ and we take $p=\ell$).
Potentially there can be links here with Lusztig's
{\em Frobenius map\,} in quantum groups, but this is beyond
this paper.
\smallskip

Let $\check{k}=k$ mod $\check{N}$ for 
$\check{N}=N/d,\, k,\check{k}\!\in\! \Z_+/2$.
Assuming that $\check{k}<\check{N}/2$, the irreducible
$\check{\HH}$\~module $V_{2\check{N}-4\check{k}}$ is a
quotient of  $V_{2N-4k}$.
 This links $N$ and $\check{N}$, but does not result in
any immediate connections at level of the corresponding
images of $\b_1$.

For instance, one can assume here that $k<\check{N}/2$ or
that $d$ is odd and $|k'|<\check{N}/2$; then we have
that $\check{k}=k$ and, respectively, $\check{k}'=k'$. 
In the latter case, the module $V'_{2|k'|}$ remains 
irreducible for $q\mapsto \check{q}$.
\smallskip

{\sf Dimension 4.}
Under the inequality  $|k'|<\tilde{N}/2$ for odd $d$,
let us consider the
module $V=V'_4$ from
(\ref{v4prime}) as an example.
It is irreducible for a primitive $(2\!N)$th
root of unity $q^{1/2}$ and
remains irreducible when $q^{1/2}\mapsto \check{q}^{1/2}$
for a primitive
root $\check{q}^{1/2}$ of order $2\!\check{N}$. Following
Section \ref{sec:dim24}, we denote by
$\mathfrak{B}^\dag_N$  and $\mathfrak{B}^\dag_{\check{N}}$,
the quotients of $\,Image(\b_q)\,$ or $\, Image(\b_1)\,$
by the center.
Then one has the homomorphisms from (\ref{NcheckN})
\begin{align}\label{NcheckNx}
&\mathfrak{B}^\dag_N/\!<\!Z\!>\ \to\
\mathfrak{B}^\dag_{\check{N}}/\!<\!\check{Z}\!>\ ,\
X\mapsto X^d, Y\mapsto Y^d,
T\mapsto T^d,
\end{align}
for $Z=Y^2T^{-2}X^{-2}$ and $\check{Z}=Y^{2d}T^{-2d}X^{-2d}$;
there are no such homomorphisms without the
division by $<\!Z\!>$ unless $d=1,N$. Obviously such a division
(which means imposing the relation $T^2=Y^2X^{-2}$ in
$\mathfrak{B}^\dag_N$)
results in loss of information. However the real problem
here is in establishing compatibility of the maps from
(\ref{NcheckNx})
with the action of the Galois group $G_K$. We see no reasons for
the group $<\!Z\!>$ to be fixed under $G_K$ (upon a proper
conjugation in $\mathfrak{B}^\dag_N$).
This would hold if $<\!Z\!>$ were conjugated in
$\mathfrak{B}^\dag_N$ to the ramification subgroups $<\!T\!>$,
$<\!A\!>$, $<\!B\!>$ or $<\!C\!>$,
but this is not the case.

Instead of division by $<\!Z\!>$
and trying the maps from (\ref{NcheckNx}), switching to
the commutator subgroups
seems the best option here. Let
$$
\c_N^\dag\equal
<\!T^2,X^2,Y^2\!>\!/\!<\!T^{2M},X^{2N},Y^{2N}\!>\,\simeq\,
\Z_{M}\times \Z_N^2,
$$
where $M=N$ for odd $N$ and $M=N/2$ otherwise. We
follow (\ref{commsubg}) and the corresponding
part of Section \ref{sec:special-c}.

Recall that 
$\mathfrak{B}^\dag_N/\c_N^\dag\simeq \Z_2^2$ 
corresponds to the cover of
$P^1$ by $E'\simeq E$ for the composition
$E'\to E'/E'_2\simeq E\to P^1$, where
the fist map is the multiplication by $2$ and $E'_2$ is
the group of points of order $2$.
The cover of $E'$ associated with $\c_N^\dag$
corresponds to $K^\circ(E')$ extended by the functions
$f_i^{1/N}$ for $f_i\in K^\circ(E')$ with the divisors
$(f_i)=2\,0_i-2\,0$, where $0_i\, (i=1,2,3)\,$ are 
nonzero points of the second order of $E'$.

We obtain the following system of homomorphisms: 
\begin{align}\label{NcheckNxx}
&\c^\dag_N\,\to\,
\c^\dag_{\check{N}},\where
X\mapsto X^d,\, Y\mapsto Y^d,\,
T^2\mapsto T^{2d}\, \hbox{ for odd $d$},
\end{align}
which is compatible with the projective action
of $PSL_2(\Z)$ and the action of $G_K$. Thus the latter action
results in the following ``ramified" analogue of the Tate module
(though not that interesting).
For an odd prime $\ell$ and $N=\ell^u$, the absolute Galois
group  $G_K$ acts in 
$\t_{\ell}\equal\varprojlim_{\,u}\c^\dag_{\ell^u}$ by
$3\times 3$-matrices with $\ell$-adic entries.

\smallskip
The action of the absolute Galois group
is actually quite similar to
that of the projective $PSL_2(\Z)$ (more generally,
the  corresponding braid groups). Technically, we use the 
same rigidity argument for both. 
In the case of $t=1 (k=0)$, the connection is the
most direct; namely, the absolute Galois group acts via  
the projective $PSL_2(\Z)$. 
Both groups are important
ingredients in the theory of the corresponding
{\em Hurwitz spaces\,}; see e.g. Section 3.1 in \cite{Det}.  
\smallskip

The projective action of
$PSL_2(\Z)$ in $\c_N^\dag$ can be readily calculated.
Since we divide by the center, we omit the powers
of $q$ in the formulas below.
It is determined by the relations
\begin{align}\label{tau-triv}
\tau_-(X^2)=Z=\tau_+(Y^{-2}) \for
Z=Y^2T^{-2}X^{-2} \hbox{\, from
(\ref{Cformula})}.
\end{align}
For instance, let us check the first formula:
$$
\tau_-(X^2)=YXYX=Y (YXT^2) X=Y^2 (T^{-1}X^{-1})(X^{-1}T^{-1})=Z.
$$
Recall that $\tau_-(Y)=Y,\,\tau_+(X)=X,\,\tau_{\pm}(T)=T\,$.
We see that the action of $PSL_2(\Z)$
in $\c_N^\dag$ factors 
through its projection onto $PSL_2(\Z_2)$.
\smallskip

{\sf Refined Jones polynomials.}
Formula (\ref{tau-triv}) with $q$ (i.e. considered
in the whole $\b_q$) and upon its extension
to all powers of $X$ and $Y$ can be used to
calculate  the images of the 
{\em refined Jones polynomials} 
of torus knots from \cite{CJ} upon the
substitution $t=-q^{-1}$
(in the case of $A_1$). The torus knots are encoded by the
(first columns) of elements of $\ga\in PSL_2(\Z)$ and their 
refined Jones polynomials are the evaluations in $\HH$ of the 
elements
$\ga(p_n)$, considered as elements of $\HH$, 
for arbitrary symmetric Macdonald polynomials $p_n$, which
add {\em colors\,} to the theory.

For $V$ of types $(\al,\be,\ga)$ from Theorem \ref{psl-two-inv}
(and their generalizations for any root systems),
the evaluations of $\ga(p_n)\in \HH$ coincide with 
the values of $\ga(p'_n)$ at $t^{-1/2}$ for the corresponding
$q,t\,$, where $p'_n$ is the image of $p_n$ in $V$ assuming
that it is well defined (i.e. for sufficiently small $n$). 
We employ the fact that $\ga$ act in $V$ (the rigidity of $V$). 
\smallskip

The nonsymmetric polynomials $e_n$ can be used instead
of $p_n$ in this construction. For instance, 
$\ga(X)$ is needed for $n=1$, which requires a calculation
entirely within the group $\,Image(\b_q)$. In this case,
one can disregard powers of $q$ and calculate modulo the center;
the formulas from \cite{CJ} are anyway up to powers of $q,t$.
The colored case ($n>1$) is similar,
though we need to track now the powers of $q$.
For example, formulas (\ref{tau-triv})
with proper $q$-corrections provide simple expressions
for $\ga(p_n)$ for even $n$ (under the substitution
$t\!=\!-q^{-1}$). Note that $\ga(p'_n)(t^{-1/2})$ becomes a pure
power of $q$, i.e. trivial, for $t=1(k=0)$ and any $n$ (the same
holds for any root systems), since $e_n$ are monomials in this
case.

Due to Theorem \ref{GALACTION}, the Galois group
$G^\circ_K$ acts in $V$ for the modules
$V$ considered there. Hence the evaluation of $g(p'_n)$,
where $p_n'$ is the
image of $p_n$ in $V$ (assuming that $p_n$ is well defined)
can be considered as certain {\em arithmetic refined Jones 
polynomial\,} of
$g\in G_K^\circ$. It is a collection of elements from 
$\Z[\ddot{q}]/(p^m)$ for $p,m\,$ as above (trivial for
$t=1$). They are actually
from  $\Q(\ddot{q})$\, for finite $\,Image(\b_q)$; can this
be true for infinite $\,Image(\b_q)$?
\smallskip

{\sf Concluding remarks.}
$(a)\,$ A generalization of
Theorem \ref{GALACTION} to rigid modules for any root
systems, at least to {\em perfect representations\,}
from \cite{C101}, is a natural challenge.
Such modules describe the monodromy
of certain local systems connected with the {\em KZB\,}
or {\em elliptic QMBP\,}, which can be hopefully
defined over algebraic numbers. This would give an alternative
approach to the action of the absolute Galois group in the
rigid DAHA modules, without using the Riemann Existence Theorem.
The natural setting for the approach based on {\em RET\,}
is the root system $C^\vee C_1$, where quite a ramified
theory at roots of unity is expected. 
\smallskip

$(b)\,$ On the other hand, even for $A_1$ and
for the simplest nontrivial $3$-dimensional modules
of DAHA, the corresponding images of the elliptic
braid group exactly match the {\em Livn\'e groups}.  Furthermore,
the simplest nontrivial case $q=e^{2\pi\imath/3}, t=-1$
of Proposition \ref{LIVNE} results
in the Livn\'e lattice for $\phi=\pi/3$, directly related
to the Eisenstein-Picard modular group, which has interesting
algebraic and analytic theory.
Starting with an
elliptic curve over $K\subset \overline{\Q}$, we obtain an
action of Gal$(\overline{\Q}/K)$ in $p$-adic completions of
this modular group for $p\neq 2,3$ and generalize this construction
to arbitrary rigid DAHA-modules of type $A_1$.
\smallskip

$(c)\,$ One of the major facts of the DAHA theory 
is the projective action
of $PSL_2(\Z)$
in the algebra itself and its rigid irreducible representations.
This property is entirely conceptual and is directly related to
the topology of the elliptic configuration space. In the present
paper, we use such an approach for the absolute Galois group
$G_K$ instead of the projective $PSL_2(\Z)$. As an application,
certain {\em refined Jones
polynomials\,} of $g\in G_K$ can be defined in type $A_1$,
which are systems of elements from $\Z[\ddot{q}]/(p^m)$,
instead of those in \cite{CJ} for torus knots
encoded  by the corresponding matrices $\ga\in PSL_2(\Z)$.
\medskip

{\bf Acknowledgements.}
The author thanks David Kazhdan, Maxim Kontsevich, Nikita
Nekrasov and Yan Soibelman for useful discussions and Hebrew
University, IHES and SCGP (where this work was reported)
for invitations and hospitality.

\bibliographystyle{unsrt}

\end{document}